\documentclass[a4paper,11pt]{article}
\small
\normalsize

\usepackage[utf8]{inputenc}
\usepackage[T1]{fontenc}
\usepackage[english]{babel}
\usepackage{times}
\usepackage{frcursive}
\usepackage{ulem}
\usepackage{color}
\usepackage[pdftex]{hyperref}

\usepackage{theorem}

\usepackage[all]{xy}

\usepackage{stmaryrd}
\usepackage{graphicx}
\usepackage{paralist}
\usepackage{amsfonts}
\usepackage[leqno]{amsmath}
\usepackage{amssymb,mathrsfs}
\usepackage{color}
\usepackage{times}
\usepackage{enumerate,vmargin}
\usepackage{tocvsec2}
\usepackage{verbatim}
\usepackage[labelfont=bf,font=it]{caption}
\usepackage[numbers,sort&compress]{natbib}

\setmarginsrb{2.9cm}{2.6cm}{2.9cm}{1.7cm}{0cm}{0mm}{0cm}{9mm}

\hypersetup{
    unicode      = false,     
    pdftoolbar   = true,      
    pdfmenubar   = true,      
    pdffitwindow = true,      
    pdfnewwindow = true,      
    colorlinks   = true,      
    linkcolor    = blue,      
    citecolor    = red,      
    filecolor    = blue,      
    urlcolor     = green       
}

\definecolor{vio}{rgb}{.5,.1,.5}
\definecolor{grey}{rgb}{.5,.5,.5}
\definecolor{cnic}{rgb}{.9,0.5,0}
\definecolor{mycolor}{rgb}{.6,0,.8}

\newtheorem{lemma}{Lemma}[section]
\newtheorem{theorem}[lemma]{Theorem}
\newtheorem{proposition}[lemma]{Proposition}
\newtheorem{corollary}[lemma]{Corollary}

{\theorembodyfont{\rmfamily}}
{\theorembodyfont{\rmfamily}}
{\theorembodyfont{\rmfamily}}
{\theorembodyfont{\rmfamily}\newtheorem{remark}[lemma]{Remark}}
{\theorembodyfont{\rmfamily}\newtheorem{rem}[lemma]{Remark}}
{\theorembodyfont{\rmfamily}\newtheorem{definition}[lemma]{Definition}}
{\theorembodyfont{\rmfamily}\newtheorem{example}[lemma]{Example}}

\newcommand{\noi}{\noindent}

\newcommand{\lgeo}{\llbracket}
\newcommand{\rgeo}{\rrbracket}

\newcommand{\un}{{\bf 1}}

\newcommand{\cA}{\mathcal A}

\newcommand{\cD}{\mathcal D}

\newcommand{\cF}{\mathcal F}

\newcommand{\cR}{\mathcal R}

\newcommand{\bE}{\mathbf E}

\newcommand{\bP}{\mathbf P}
\newcommand{\bQ}{\mathbf Q}

\newcommand{\bm}{\mathbf m}

\newcommand{\bmu}{\boldsymbol \mu}
\newcommand{\fTheta}{\boldsymbol \Theta}

\newcommand{\gzH}{\mathbf{h}}

\newcommand{\bbU}{\mathbb U}
\newcommand{\bbZ}{\mathbb Z}

\newcommand{\bbM}{\mathbb M}
\newcommand{\bbN}{\mathbb N}

\newcommand{\bbR}{\mathbb R}
\newcommand{\bbT}{\mathbb T}
\newcommand{\bbW}{\mathbb W}
\newcommand{\boo}{\varrho}
\newcommand{\booo}{\varrho}

\newcommand{\bgg}{\mathtt b}

\newcommand{\Wid}{\bbW^*_{\! \bgg }}
\newcommand{\Wif}{\bbW^*_{\! [0, 1]}}

\newcommand{\Wbd}{\overline{\bbW}_{\! \bgg}}
\newcommand{\Wbf}{\overline{\bbW}_{\! [0, 1]}}
\newcommand{\rmo}{\mathrm{o}}
\newcommand{\gzW}{\mathcal W}

\newcommand{\ccF}{\mathscr F}
\newcommand{\ccG}{\mathscr G}

\def\cq{$\hfill \square$}
\def\cqfd{$\hfill \blacksquare$}
\def\ino{ \! \in \! }

\def\bC{\mathbf{C}}

\def\bT{\mathbf{T}}

\newcommand{\bdelta}{{\boldsymbol{\delta}}}

\newcommand{\bx}{\mathbf{x}}

\newcommand{\bt}{\mathbf{t}}

\newcommand{\ttT}{\mathtt{T}}
\newcommand{\ttR}{\mathtt{R}}

\def\bu{\mathbf{u}}

\newcommand{\ftau}{\boldsymbol{\tau}}
\newcommand{\fSigma}{\boldsymbol{\Sigma}}
\newcommand{\bbd}{\mathbf{d}}
\definecolor{orange}{rgb}{1,0.5,0}

\def\epp{\varepsilon}
\def\bba{\mathbf{a}}
\def\ttH{\mathtt{H}}
\def\MMM{\mathbf{M}}
\def\MMT{\mathbf{MT}}

\definecolor{italiangreen}{RGB}{0,140,69}
\newcommand{\green}{\color{italiangreen}}

\usepackage[refpage,noprefix]{nomencl}                          
\makeglossary
\makenomenclature                                           

\graphicspath{{.}{figures/}}
\title{ \textsc{Scaling limits of tree-valued branching random walks}}
\date{}

\author{Thomas \textsc{Duquesne}
\thanks{Sorbonne Universit\'e, Campus Pierre et Marie Curie, 
LPSM, Case courrier 158, 4, place Jussieu, 
75252 Paris Cedex 05 
France. Email: \textit{\{thomas.duquesne, robin.khanfir, shen.lin\}@sorbonne-universite.fr}}
\and Robin \textsc{Khanfir}\footnotemark[1]
\and Shen \textsc{Lin}\footnotemark[1]
\and Niccol\`o \textsc{Torri}
\thanks{{Universit\'e Paris-Nanterre, Laboratoire MODAL'X, UMR CNRS 9023 and FP2M, CNRS FR 2036, France.}
Email: \textit{ntorri@parisnanterre.fr}. N. Torri was also supported by the project Labex MME-DII (ANR11-LBX-0023-01).
}
}

\begin{document}

\maketitle

\begin{abstract}
We consider a branching random walk (BRW) 
taking its values in the $\mathtt{b}$-ary rooted tree $\bbW_{\! \mathtt{b}}$
(i.e.~the set of finite words written in 
the alphabet $\{ 1, \ldots, \mathtt{b} \}$, with $\mathtt{b}\! \geq \! 2$). 
The BRW is indexed by a critical Galton--Watson tree conditioned to have $n$ {vertices}; 
its offspring distribution is aperiodic and is in the domain of attraction of a $\gamma$-stable law, $\gamma \ino (1, 2]$. The jumps of the BRW are those of a 
{nearest-neighbour} null-recurrent random walk on $\bbW_{\! \mathtt{b}}$ 
(reflection at the root of $\bbW_{\! \mathtt{b}}$ 
and otherwise: probability $1/2$ to move closer 
to the root of $\bbW_{\! \mathtt{b}}$ and probability $1/(2\mathtt{b})$ 
to move away from it to one of the $\mathtt{b}$ sites above). 
We denote by $\cR_{\mathtt{b}} (n)$ the range of the
BRW in $\bbW_{\! \mathtt{b}}$ which is the set of all sites in $\bbW_{\! \mathtt{b}}$ visited by the BRW. We first prove a law of large numbers for $\# \cR_{\mathtt{b}} (n)$ and  
we also {prove} that if we equip 
$\cR_{\mathtt{b}} (n)$ (which
is a random subtree of $\bbW_{\! \mathtt{b}}$) with its graph-distance
$d_{\mathtt{gr}}$, then there exists a scaling sequence 
$(a_n)_{n\in \bbN}$ satisfying $a_n \! \rightarrow \! \infty$ such that 
the metric space $(\cR_{\mathtt{b}} (n), a_n^{-1}d_{\mathtt{gr}})$, equipped with its normalised empirical measure, converges to the reflected Brownian cactus with $\gamma$-stable branching mechanism: namely, a random compact real tree that is a
variant of the Brownian cactus introduced by N.~Curien, J-F.~Le Gall and G.~Miermont in \cite{CuLGMi13}.

\medskip

\noi
\textbf{Keywords} $\, $ Branching random walks $\cdot$ Galton--Watson tree $\cdot$ Scaling limit $\cdot$ Superprocess $\cdot$ Brownian snake $\cdot$ Brownian cactus
$\cdot$ Real tree

\medskip

\noi
\textbf{Mathematics Subject Classification} $\, $ 60J80  $\cdot$ 60G50 $\cdot$ 60G52 $\cdot$ 60F17

\end{abstract}

\section{Introduction}
\label{introsec}

 Since the seventies, branching random walk (BRW) is an area of research that is intensively studied and is linked to travelling wave solutions of 
semi-linear partial differential equations (FKPP) or various models of statistical mechanics (Generalized random energy model, Mandelbrot's cascades, Gaussian free field): 
we refer to the book of Z.~Shi \cite{Shi15} for an overview of this topic; we also refer to the works of S.~Gouëzel, I.~Huerter, S. Lalley and T.~Sellke
\cite{LalSel97, LalHue00, Lal06, LalGou13} for the study of BRW in hyperbolic spaces and to 
T.~Liggett \cite{Lig96} and to I.~Benjamini and S.~Müller \cite{BenMue12} for branching random walks on trees.
In most of the previous works, BRWs are indexed by an infinite supercritical 
Galton--Watson tree (GW-trees) and questions focus on various survival events or extremal behaviours of  BRWs.

In this paper, we consider instead a BRW that takes its values 
in the $\mathtt{b}$-ary tree {$\bbW_{\! \mathtt{b}}$} and that is indexed by a critical 
Galton--Watson tree conditioned to have $n$ vertices. 
The jumps of the BRW are those of a 
{nearest-neighbour} null-recurrent RW on $\bbW_{\! \mathtt{b}}$. Namely,  
at the root of $\bbW_{\! \mathtt{b}}$ (that is denoted by $\varnothing$), it is reflected and elsewhere,  
with probability $1/2$, it jumps {to the neighbour} closer to the root of $\bbW_{\! \mathtt{b}}$ and with probability $1/(2\mathtt{b})$, it moves further from the root of $\bbW_{\! \mathtt{b}}$ and it jumps to one of the $\mathtt{b}$ sites above. 
We study the range $\cR_{\mathtt{b}} (n)$ of {this} BRW when $n\! \rightarrow \! \infty$. 
More precisely, we first show that $\frac{_1}{^n} \# \cR_{\mathtt{b}} (n) $ converges {in probability} to a constant; this law of large numbers is the analogue of the results due to J-F. Le Gall and L.~\cite{LGLi15, LGLi16} who treat the cases of $\bbZ^d$-valued BRWs that are indexed by critical GW-trees conditioned to have $n$ vertices.  
We then prove that $\cR_{\mathtt{b}} (n)$, seen as a subtree of the $\mathtt{b}$-ary tree, converges, when it is suitably rescaled, to a continuum random tree called the  reflected Brownian cactus with $\gamma$-stable branching mechanism: namely, a random compact real tree that is a
variant of the Brownian cactus introduced by N.~Curien, J-F.~Le Gall and G.~Miermont in \cite{CuLGMi13}. 

This limit theorem for $\cR_{\mathtt{b}} (n)$ is related to earlier works on scaling limits of the range 
of tree-valued critical or near-critical biased random walks (RWs): in particular we refer to D.~\cite{Du05} who deals with {near-critical biased} RWs on $\mathtt{b}$-ary trees, to Y.~Peres and O.~Zeitouni \cite{PerZei08} 
who show that the distance to the root of a critical {biased} RW in a Galton--Watson environment is diffusive, to A.~Dembo and N.~Sun \cite{DeSu12} who study the cases of critical biased RWs on $N$-type GW-trees,  
to E.~Aïdékon and L.~de Raphélis \cite{AidRap17} who improve Y.~Peres and O.~Zeitouni's result and who show that the range of the same RW 
converges when suitably rescaled, to a variant of the Brownian CRT, and to X.~Chen and 
G.~Miermont \cite{ChMi17} who show that rescaled {Brownian} bridges and loops in hyperbolic spaces converge to the Brownian CRT. Their work is based on a previous results due to P.~Bougerol and T.~Jeulin \cite{BoJe99}. Independently, A.~Stewart shows in his PhD Thesis \cite{MR3697600} that the rescaled simple RW bridges on a $d$-regular tree ($d\geq 3$) converge to the Brownian CRT.

Let us describe more precisely the results that we obtain. 
We consider a (rooted and ordered) Galton--Watson tree $\tau$ with offspring distribution $\mu$ that satisfies the following: we fix $\gamma \ino (1, 2]$ and we assume 
\begin{equation} 
\label{hyposta} 
(\mathbf{H}) \quad \left\{ 
\begin{array}{l}
\!\! (\mathrm{H}_1):  \; \sum_{k\in \bbN} k\mu(k)= 1 , \\
\!\! (\mathrm{H}_2): \; \textrm{$\mu$ is aperiodic {(namely, $\mu$ is not supported by a proper subgroup of $\bbZ$),}} \\
\!\! (\mathrm{H}_3):\textrm{Either $\gamma\! = \! 2$ and $\sum_{k\in \bbN} k^2\mu(k) \! < \! \infty$, or $\gamma\ino (1, 2)$ and $\mu$ is in } \\ 
\; \; \; \; \; \; \; \; \, \textrm{the domain of attraction of a $\gamma$-stable law. }
\end{array} \right.
\end{equation}
Note that $(\mathrm{H}_1)$ implies that a.s.~the total number of vertices $\# \tau$ is finite; $ (\mathrm{H}_2)$ 
implies that for all large enough integers $n$, $\bP (\# \tau \! = \! n)\! >\! 0$. We 
translate $ (\mathrm{H}_3)$ into the following assertion: let $(L_n)_{n\in \bbN}$ be an i.i.d.~sequence of 
$(\{\! -1\} \! \cup \! \bbN) $-valued {random variables} such that $\bP (L_n\! = \! k)\! = \! \mu (k+1)$, 
$k\! \geq \! -1$ and let $X$ be a real {random variable} whose law is spectrally positive 
$\gamma$-stable; it is characterised by its Laplace exponent: $\log \bE [\exp (-\lambda X)]\! = \! \lambda^\gamma$, $\lambda \ino [0, \infty)$. Then under $ (\mathrm{H}_3)$, there exists a nondecreasing $\tfrac{\gamma -1}{\gamma}$-regularly varying sequence $(a_n)_{n\in \bbN}$ such that
\begin{equation}
\label{Favoriten}
\frac{a_n}{n} \big( L_1 + \ldots + L_n -n  \big) 
\overset{\textrm{(law)}}{\underset{n\rightarrow \, \infty}{-\!\!\! -\!\!\! -\!\!\! \longrightarrow}} X \; .
\end{equation}
As we see below, 
under $(\mathbf{H})$, $\tau$ behaves regularly when it is conditioned to be large: namely, we see 
that suitably rescaled versions of $\tau$ under $\bP (\, \cdot \, | \, \# \tau\! = \! n)$ converge in distribution 
when $n\! \rightarrow \! \infty$.

Conveniently, we view $\tau$ as a family tree whose ancestor is the root and where 
siblings are ordered by birth-rank. The depth-first exploration of $\tau$ is 
the sequence of vertices $(u_k)_{0\leq k < \# \tau}$ that is defined recursively as follows: 
$u_0$ is the root and for all $k\ino \{ 0, \ldots, \#\tau \! -\! 2\}$, let $v$ be the most recent 
ancestor of $u_k$ having at least one unexplored child (note that possibly $v\! = \! u_k$); 
then $u_{k+1}$ is the unexplored child of $v$ with least birth-rank. Our first result is 
the following law of large numbers for the size of the range of the $\bbW_{\! \bgg}$-valued $\tau$-{indexed} critical branching random walk (see Sections \ref{empirisec} and \ref{treeBRWsec} for more precise definitions).
\begin{theorem}
\label{main1} Let $\tau$ be a Galton--Watson tree with offspring distribution $\mu$ that satisfies $(\mathbf{H})$. Recall that $(u_k)_{0\leq k < \# \tau}$ stands for the depth-first exploration of $\tau$. Conditionally given $\tau$, let $(Y_{v})_{v\in \tau}$ be a $\bbW_{\! \bgg}$-valued $\tau$-indexed critical branching random walk starting at $Y_{\mathtt{root}}\! = \! \varnothing$. Then, there exists a constant $c_{\mu, \bgg}\ino (0, \infty)$ that only depends on $\mu$ and $\bgg$ such that 
\begin{equation}
\label{main1eff}
\forall \epp \ino (0, \infty), \quad \lim_{n\rightarrow \infty} \bP \Big(  \tfrac{1}{n}\max_{1\leq k\leq n}  
\big|\# \big\{Y_{u_l}; 0 \! \leq \! l \! < \! k\big\} \, -\, c_{\mu, \bgg} \, k \big| \, >\epp \, \Big|\,  \# \tau \! = \! n \Big) = 0 \; , 
\end{equation}
In particular, for all $\epp \ino (0, \infty)$, we get
$ \lim_{n\rightarrow \infty}\bP \big(\,  \big| \frac{1}{n} \# \cR_\bgg -c_{\mu, \bgg} \big| \! >\! \epp \, \big| \, \#\tau \! = \! n \big) \! = \! 0$, where $\cR_{\bgg}$ stands for the range $\{ Y_v; v\ino \tau \}$ of the branching random walk.   
\end{theorem}

  Let $\tau$, $(Y_v)_{v\in \tau}$ and $\cR_{\bgg} \! = \! \{ Y_v; v\ino \tau \}$ be as in 
Theorem \ref{main1}. Observe that $\cR_{\bgg}$ is a subtree of $\bbW_{\! \bgg}$. Our second main result is a limit theorem for rescaled versions of the metric spaces $(\cR_\bgg, d_{\mathtt{gr}})$ where 
$ d_{\mathtt{gr}}$ stands for the graph-distance. To state it, let us first recall a limit theorem for $(\tau,  d_{\mathtt{gr}})$. 
Set $H_k (\tau)\! = \! d_{\mathtt{gr}} (\mathtt{root}, u_k)$ for all $k\ino \{ 0, \ldots, \# \tau \! -\! 1 \}$, that is the height process of $\tau$. Note that 
$H(\tau)\! =\! \big(H_k(\tau)\big)_{k\in \{ 0, \ldots, \# \tau - 1 \}}$ entirely codes $\tau$. Then, Theorem 3.1 in D.~\cite{Du03} asserts the following: 
\textit{assume $(\mathbf{H})$ as in (\ref{hyposta}) and let $(a_n)$ be as in \eqref{Favoriten}; then there is a nonnegative continuous process $H=(H_{s})_{s\in [0, 1]}$ such that }
\begin{equation}
\label{Loos}
\big( \tfrac{1}{a_n} H_{\lfloor ns \rfloor}(\tau) \big)_{\! s\in [0, 1]} \quad \textit{under} \quad  \bP \big( \, \cdot \, \big| \, \# \tau \! = \! n \big) \; \underset{n\rightarrow \infty}{- \!\!\! - \!\!\! - \!\!\! \longrightarrow } \; H, 
\end{equation}
\textit{weakly on $\bC([0, 1], \bbR)$}. When $\gamma\! \! = \! 2$, $H$ is the normalised Brownian excursion and this result is due to Aldous (see {Theorem} 23 in Aldous \cite{Al93}). When $\gamma \ino (1, 2)$, $H$ is the normalised excursion of the $\gamma$-stable height process that is a local-time function of a 
$\gamma$-stable spectrally positive L\'evy process. 

The metric space that is the limit of $(\tau,  \tfrac{1}{a_n} d_{\mathtt{gr}})$ as $n\! \rightarrow \! \infty$ 
is derived from the normalised excursion of the $\gamma$-stable height process $H$ as follows: for all 
$s_1, s_2 \ino [0, 1]$, 
we set
\[d_H (s_1, s_2) \! = \! H_{s_1}+ H_{s_2} -2\min_{s_1\wedge s_2 \leq s \leq 
s_1 \vee s_2} H_s.\] We easily check that a.s.~$d_H$ is a pseudo-metric on $[0, 1]$. We introduce the relation $\sim_H $ on $[0, 1]$ by setting $s_1 \sim_H s_2$ {if and only if} $d_H (s_1, s_2)\! = \! 0$; 
clearly, $\sim_H$ is an equivalence relation and the normalised $\gamma$-stable L\'evy tree is taken as the quotient space $ T_H \! = \! [0, 1 ] / \! \sim_H $, equipped with the distance induced by $d_H$ that we keep denoting $d_H$. We denote by $p_H \! : \! [0, 1] \rightarrow T_H$ the canonical projection. Note that $p_H$ is continuous; therefore $T_H$ is compact and connected{.} Moreover, $T_H$ is a real tree, namely, a metric space such that all pairs of points are joined by a unique simple arc that turns out to be a geodesic (see Definition \ref{errtredef} for more details). We set {$r_H\! := \! p_H (0)$} that is viewed as the root of $T_H$ and we equip $T_H$ with the measure $\mu_H$ that is the image of the Lebesgue measure on $[0, 1]$ via $p_H$, namely, $\int_{T_H} \! f \, d\mu_H\! = \! \int_0^1 \! f(p_H(s)) \, ds$, for all 
continuous $f\! : \! T_H \! \rightarrow \! \bbR$. The convergence (\ref{Loos}) 
then implies the following one. 
\begin{equation}
\label{Strobl}
\big( \tau, \tfrac{1}{a_n} d_{\mathtt{gr}} , \mathtt{root}, \tfrac{1}{n} \mathtt{m}  \big) 
\quad \textrm{under} \quad  \bP \big( \, \cdot \, \big| \, \# \tau \! = \! n \big) \; \underset{n\rightarrow \infty}{- \!\!\! - \!\!\! - \!\!\! \longrightarrow } \; (T_H, d_H, r_H, \mu_H) 
\end{equation}
where $\mathtt{m}\! = \! \sum_{v\in \tau} \delta_v$ stands for the counting measure on $\tau$. Here the convergence holds weakly on the space $\bbM$ of isometry classes of pointed measured compact metric spaces equipped with the Gromov--Hausdorff--Prokhorov distance 
$\bdelta_{\mathtt{GHP}}$ that makes it a Polish space, as proved in {Theorem 2.5 of R.~Abraham, J-F.~Delmas and  P.~Hoscheit \cite{AbDeHo13}}. (See (\ref{defGHP}) for a precise definition of $\bdelta_{\mathtt{GHP}}$ and see (\ref{rappGHP}) for more details.) {For more details on L\'evy trees see J-F.~Le Gall and Y.~Le Jan \cite{LGLJ98} and D.~and  J-F.~Le Gall \cite{DuLG02, DuLG05} (see also Section \ref{statreesc}). }

The limit of rescaled versions of the metric spaces $(\cR_\bgg, d_{\mathtt{gr}})$ is constructed as 
follows: as proved in {D.~and  J-F.~Le Gall}~\cite{DuLG05} (Lemma 6.4 p.~600, that is recalled in Lemma \ref{Bregenz}), conditionally given $H$, there exists a Hölder-continuous centered Gaussian process $\sigma \ino T_H \! \longmapsto  \! W_\sigma \ino \bbR$ whose covariance is characterised by 
$\bE \big[ \big| W_{\sigma_1} \! - W_{\sigma_2} \big|^2 \big| \, H \big] \!  = \! 
d_H (\sigma_1, \sigma_2)$, for all $\sigma_1, \sigma_2 \ino T_H$. Then, we set 
\begin{equation*}
\forall \sigma_1, \sigma_2\ino T_H, \quad d_{H, W} ( \sigma_1, \sigma_2) = |W_{\sigma_1}  |+ | W_{\sigma_2}| -2 \!\!\!\! \!\! \!   \min_{ \quad \sigma \in \lgeo \sigma_1, \sigma_2 \rgeo} \!\! \!\! \!\! \! |W_\sigma| , 
\end{equation*}    
where $\lgeo \sigma_1, \sigma_2 \rgeo $ is the unique geodesic that joins $\sigma_1$ to 
$\sigma_2$ in $T_H$. In Lemma \ref{snkmtrcL}, we prove that $d_{H, W}$ is a pseudo-metric on $T_H$;  
we then define the equivalence relation $\sim_{H, W} $ on $T_H$ by setting $\sigma_1 \sim_{H, W}  \sigma_2$ {if and only if}
$d_{H, W} ( \sigma_1, \sigma_2) \! = \! 0$ and we denote by $T_{H, W}\! = \! T_H / \! \sim_{H, W}$ the quotient metric space and we keep denoting by 
$d_{H, W}$ the resulting metric; we denote by $\pi_{H, W}\! : \! T_H \! \rightarrow \! T_{H, W}$ the canonical projection that is continuous. Thus $T_{H, W}$ is compact and connected, and 
$(T_{H, W}, d_{H, W})$ is a real tree (see Section \ref{Propcact}, especially Proposition \ref{Dornbirn}, for more properties of $T_{H, W}$). 
It turns out that this kind of spaces has been introduced in {N.~Curien, J.-F.~Le Gall and G.~Miermont \cite{CuLGMi13} (see also J.-F.~Le Gall \cite{LG15} for a different purpose)}; 
they coined the name Brownian cactus, so we call $(T_{H, W}, d_{H,W})$ the normalised reflected Brownian cactus with $\gamma$-stable branching mechanism.  
We next set $r_{H, W}\! = \! \pi_{H, W} (r_H)$ that is viewed as the root of $T_{H, W}$ 
and we equip $T_{H, W}$ with the measure $\mu_{H, W}$ that is the image 
of $\mu_H$ via $\pi_{H, W}$: namely, $\int_{T_{H, W}} \! f \, d\mu_{H, W}\! = \! \int_{T_H} \! f(\pi_{H, W}(\sigma)) \, \mu_H (d\sigma) $, for all 
continuous $f\! : \! T_{H, W} \! \rightarrow \! \bbR$. Our second result is the following limit theorem. 
\begin{theorem}
\label{main2} Let $\tau$ be a Galton--Watson tree with offspring distribution $\mu$ that satisfies $(\mathbf{H})$ as in (\ref{hyposta}). We denote by $\mathtt{m}\! = \! \sum_{v\in \tau} \delta_v$ the counting measure on $\tau$ and we denote by $d_{\mathtt{gr}}$ the graph distance on $\tau$. 
Conditionally given $\tau$, let $(Y_{v})_{v\in \tau}$ be a $\bbW_{\! \bgg}$-valued $\tau$-indexed critical branching random walk starting at $Y_{\mathtt{root}}\! = \! \varnothing$. 
We denote by $\cR_{\bgg}\! = \! \{ Y_v; v\ino \tau\}$ the range of $Y$, by $\bm^\bgg_{\mathtt{occ}}\! = \! \sum_{v\in \tau} \delta_{Y_v}$ the occupation measure of $Y$ and by $d_{\mathtt{gr}}$ the graph distance on $\cR_{\bgg}$. 
We denote by $(T_H, d_H, r_h, \mu_H)$ the normalised $\gamma$-stable L\'evy tree and by 
$(T_{H, W} , d_{H, W}, r_{H, W}, \mu_{H, W})$ the normalised reflected Brownian cactus with $\gamma$-stable branching mechanism as defined above. Let $(a_n)_{n\in \bbN}$ be as in (\ref{Favoriten}). 
Then, the following limit holds weakly on $(\bbM, \delta_{\mathtt{GHP}})^2$ 
\begin{eqnarray}
\label{Mattighofen}
\lefteqn{ \Big( \big( \tau, \tfrac{1}{a_n} d_{\mathtt{gr}}, \mathtt{root}, \tfrac{1}{n} \mathtt{m} \big) \, , \, \big( \cR_{\bgg} , \tfrac{1}{\sqrt{a_n}} d_{\mathtt{gr}}  , \varnothing, \tfrac{1}{n} 
\bm^\bgg_{\mathtt{occ}}\big) \Big)
\; \textrm{under $\; \bP \big(  \cdot  \big| \, \# \tau \! = \! n \big) $ } } \hspace{40mm} 
                    \nonumber\\
 & &\overset{\textrm{(law)}}{\underset{n\rightarrow \infty}{-\!\!\! -\!\!\!-\!\!\! -\!\!\!  -\!\! \! \longrightarrow}} \; 
\Big( \big( T_H, d_H, r_H, \mu_H \big) \, , \, \big(T_{H, W} , d_{H, W}, r_{H, W}, \mu_{H, W} \big) \Big).
\end{eqnarray}
Moreover, denote by 
$\bm^\bgg_{\mathtt{count}}\! = \! \sum_{x\in \cR_\bgg} \delta_x$ the counting measure on $\cR_\bgg$; then, for all $\epp \ino (0, \infty)$, \begin{equation*}
 \lim_{n\rightarrow \infty} \bP \Big(  \,  d^{_{\, (n)}}_{{\mathtt{Prok}}} \big( \tfrac{c_{\mu, \bgg}}{n}  \bm^\bgg_{\mathtt{occ}} , \tfrac{1}{n} \bm^\bgg_{\mathtt{count}}  \big)  \! >\! \epp \, \Big| \, \# \tau \! = \! n \Big)\! = \! 0\; ,
\end{equation*}
where $c_{\mu, \bgg}$ is as in Theorem \ref{main1} and where $ d^{_{\, (n)}}_{^{\mathtt{Prok}}}$ stands for the Prokhorov distance on the space of finite measures on $(\cR_\bgg, \tfrac{1}{\sqrt{a_n}}d_{\mathtt{gr}})$. It implies that the following limit holds jointly with (\ref{Mattighofen}):
$$ \big( \cR_{\bgg} , \tfrac{1}{\sqrt{a_n}} d_{\mathtt{gr}}  , \varnothing, \tfrac{1}{n} 
\bm^\bgg_{\mathtt{count}}\big) \; \textrm{under $\; \bP \big(  \cdot  \big| \, \# \tau \! = \! n \big) $ }
\overset{\textrm{(law)}}{\underset{n\rightarrow \infty}{-\!\!\! -\!\!\!-\!\!\! -\!\!\!  -\!\! \! \longrightarrow}} \; 
 \big( T_{H, W} , d_{H, W}, r_{H, W}, c_{\mu, \bgg}. \mu_{H, W} \big) . $$

\end{theorem}

Theorem \ref{main1} is the analogue of Proposition 5 and Theorem 7 in J-F.~Le Gall and L.~\cite{LGLi16} established for $\bbZ^d$-valued BRWs. Let us mention that our strategy of proof is similar: we define a specific invariant shift for infinite tree-valued BRWs and we use the subadditive ergodic theorem; the constant $c_{\mu, \bgg} $ is interpreted as the probability that the invariant BRW visits its starting point only once.  

The proof of Theorem \ref{main2} is distinct from that of Theorem \ref{main1}. On one hand, it relies on general arguments on weak limits of random metrics (see Proposition \ref{submetric}). As an application of these results, we prove that the range of critical biased RWs 
on $N$-type supercritical GW-trees converges to the tree coded by a reflected Brownian motion (see Corollary \ref{corodiff}). This result is derived from a much more difficult 
result due to A.~Dembo and N.~Sun \cite{DeSu12} that asserts that the distance from the root of the RW converges, when suitably rescaled, to a reflected Brownian motion. The same idea allows to recover previous  scaling limits for the range of RWs on supercritical GW-trees such as in D.~\cite{Du05} (critical biased RWs on $\mathtt{b}$-ary trees) or in E.~Aïdékon and L.~de Raphélis \cite{AidRap17} (biased RWs and RWs in random environment on single-type GW-trees). We refer to the end of Section \ref{singRW} for more details.    
On the other hand, the proof of Theorem \ref{main2} uses limit theorems for discrete snakes that have been obtained by S.~Janson and J-F.~Marckert in \cite{JanMar05} in the Brownian case and by C.~Marzouk \cite{Mar20} in the stable cases.

\subsection*{Organization of the paper}

The paper is organised as follows. In Section \ref{empirisec}, we introduce notations on trees with an infinite line of ancestors that constitute a natural state-space for invariant tree-valued BRWs. 
In Section~\ref{trbrwdfsc}, we define the kinds of BRWs that we study; Section \ref{memesec} is devoted to metric properties of the range of the so-called free BRWs. In Section \ref{zdczze}, we state a coupling for BRWs that is a key {argument} in the proof of Theorem \ref{main1}. In Section \ref{brwestsec}, we prove estimates that are  used mostly to prove that $c_{\mu ,\bgg}\! >\! 0$. Theorem \ref{main1} is proved in Sections \ref{lawinvsec} and \ref{Pfmain1}. In Section \ref{psdmsec}, we prove general convergence results for random metrics. In Section \ref{singRW} we apply these results to get scaling limits for the range of RWs on $N$-type GW-trees {(see Corollary \ref{corodiff})}.  In Section \ref{snkmtrsec} we introduce snake metrics and we prove specific results. 
In Sections \ref{statreesc} and \ref{rfbrsnsc}, we recall definitions and properties on stable L\'evy trees and L\'evy snakes. Section \ref{Propcact} is devoted to basic properties of reflected Brownian cactuses. 
Theorem \ref{main2} is proved in Section \ref{Pflimth}.

\paragraph*{Acknowledgements} The authors would like to thank the referee and the editor for their helpful suggestions.

\section{Tree with a possibly infinite line of ancestors.}
\label{empirisec}
\paragraph*{Words.}  Recall that $\bbN$ stands for the set of nonnegative integers $\{ 0, 1, 2 , \ldots \}$ 
and that $\bbN^*\! = \! \bbN \backslash \{ 0\}$. Let $A$ be a set with more than two elements that is viewed as an alphabet. We denote by 
$\bbW_{\! A}$ the finite words written with alphabet $A$: namely, 
\begin{equation} 
\label{Wdesun}
\bbW_{\! A}\! = \! \bigcup_{n\in \bbN} A^n\; .
\end{equation}
Here, $A^0$ is taken as $\{ \varnothing \}$, $\varnothing$ being the empty word.  
Let $u\! = \! (a_1, \ldots , a_n)\ino \bbW_{\! A}$ be distinct from $\varnothing$. We set $|u|\! = \! n$ that is 
the \textit{height} of $u$, with the convention that $|\varnothing |\! = \! 0$.  
We next set $\overleftarrow{u}\! = \! (a_1, \ldots, a_{n-1})$ that is interpreted as the \textit{parent of $u$} 
(if $n\! = 1$, then $\overleftarrow{u}\! = \! \varnothing$). 
More generally for all $p\ino \{ 1, \ldots, n\}$, we set 
$u_{| p}\! = \! (a_1, \ldots , a_p)$, with the convention: $u_{| 0}\! = \! \varnothing$. 
For all $v\! = \! (b_1, \ldots, b_m)\ino \bbW_{\! A}$, we set 
$u\ast v\! = \! (a_1, \ldots , a_n, b_1, \ldots,  b_m)$ that is 
the \textit{concatenation} of $u$ with $v$, with the convention that 
$\varnothing \ast u\! = \! u \ast  \varnothing\! = \! u$. We shall also define the most recent common ancestor of 
$u$ and $v$ in $\bbW_{\! A}$ as $u\wedge v \! = \! u_{| p} \! =\!  v_{| p}$ where $p\! = \! \max \{ k\ino \bbN \, :  u_{| k}\! =\!  v_{| k} \}$.  
We shall consider three cases: 
\smallskip

\begin{compactenum}
\item[$\bullet$] $A=\bbN^*$; in that case we use the notation $\bbU\! := \! \bbW_{\bbN^*}$, 
the letter $U$ being for Ulam.  
\item[$\bullet$] $A= \{ 1, \ldots, \bgg\}$, $\bgg$ being an integer $\geq 2$; we use the notation 
$\bbW_{\! \bgg}\! := \! \bbW_{\! \{ 1, \ldots , \bgg\}}$; $\bbW_{\! \bgg}$ is the \textit{$\bgg$-ary tree}.
\item[$\bullet$] $A= [0, 1]$; we call $\bbW_{[0, 1]}$ the \textit{free tree}. 
\end{compactenum}

\begin{definition}
\label{ort}
\textit{Rooted ordered trees} can be viewed as subsets $t\! \subset \! \bbU$ that satisfy the following. 
\begin{compactenum}

\smallskip

\item[$(a)$] $\varnothing \ino t$. 

\smallskip

\item[$(b)$] If $u\ino t\backslash \{ \varnothing\}$, then $\overleftarrow{u}\ino t$. 

\smallskip

\item[$(c)$] For all $u\ino t$, there exists $k_u (t)\ino \bbN$ such that $u \! \ast \! (i)\ino t$ 
{if and only if} $1\! \leq \! i \! \leq \! k_u (t)$.  
\end{compactenum}

\smallskip

\noi
We denote by $\bbT$ the set of rooted ordered trees. \cq 
\end{definition}
\noi
The quantity $k_u (t)$ is interpreted as the \textit{number of children of $u$} and $u \! \ast \! (i)$ is the\textit{ $i$-th child of $u$}, $1\! \leq \! i \! \leq \! k_u (t)$. 
If $k_u (t)\! = \! 0$, then there is no child stemming from $u$ and assertion $(c)$ is empty.  
We next set {the shift of $t$ at $u$ by} $\theta_u t\!  = \! \{ v\ino \bbU: u \! \ast \! v\ino t \}$ that is also a rooted ordered tree: it is viewed as the \textit{subtree of the descendants stemming from $u$}.  
Unless otherwise specified, all the random variables that are mentioned in this paper are defined on the same probability space $(\Omega, \cF, \bP)$. 

\begin{definition}
\label{GWtreedef}
We equip $\bbT$ with the sigma-field $\ccF(\bbT)$ generated by the sets $\{ t\ino \bbT\! : \! u\ino t\}$, $u\ino \bbU$. A \textit{Galton--Watson tree with offspring distribution $\mu$} (a \textit{GW($\mu$)-tree}, for short) is a $(\ccF, \ccF(\bbT))$-measurable r.v.~$\tau\! :\! \Omega \! \rightarrow \! \bbT$ that satisfies the following. 
\begin{compactenum}

\smallskip

\item[$\bullet$] $k_\varnothing (\tau)$ has law $\mu$. 

\smallskip

\item[$\bullet$] For all $k\! \geq\!  1$ such that $\mu(k)\! >\! 0$, the subtrees $\theta_{(1)} \tau, \ldots , \theta_{(k)} \tau$ under $\bP (\, \cdot \, | k_\varnothing (\tau)\! = \! k)$ are independent with the same law as $\tau$ under $\bP$. \cq 
\end{compactenum}
\end{definition}

\noi
Recall that $\tau$ is a.s.~finite if and only if $\mu$ is critical or subcritical: $\sum_{k\geq 1} k\mu (k)\! \leq \! 1$.

\paragraph*{Bilateral words.}
We shall consider branching random walks seen from the spatial and genealogical position of a tagged individual. 
To that end, it is convenient to introduce ordered trees that are rooted possibly 
at a negative generation and we also introduce their local limits that may have an infinite line of ancestors. 
It is therefore convenient to introduce words indexed by possibly negative numbers.

To simplify, we set $\bbZ_-\! = \! \bbZ \backslash \bbN^*$ and $\overleftarrow{\bbZ}\! = \! 
\bbZ  \cup  \{ -\infty \}$.    
{An} $A$-\textit{bilateral word} is a (possibly infinite) sequence $u\! =\!  (a_k)_{k_1 <k\leq k_2}$ where 
$k_1 \ino \overleftarrow{\bbZ}$ and $k_2 \ino \bbZ$ are such that $k_2 \! \geq \! k_1$ and where 
$a_k \ino A $ for all $k_1 \! <\!  k \! \leq \! k_2$. We denote by $\overline{\bbW}_{\! A}$ the set of bilateral words, including the empty word denoted by $\varnothing$ 
(if $k_1\! \geq  \! k_2$, then we agree on $u\! = \! \varnothing$). 
If $u\neq \varnothing$, we introduce the following notation. 
\begin{equation}  
\label{extremdef}\begin{aligned}
& |u|_-\! =\!  k_1 , \,\text{be the \textit{depth} of $u$,}\\
&|u|\! = \! k_2, \, \text{be the \textit{relative height} of $u$ (note that it may take negative values),}\\
&\overleftarrow{u}\! = \! (a_k)_{k_1 < k\leq k_2-1} \, \text{be the \textit{parent} of $u$,}\\
&u_{|_{(l_1, l_2]}} \! =\!  (a_k)_{k_1 \vee l_1 <k \leq k_2 \wedge l_2}\,,\\
&\mathtt{end} (u)= a_{k_2} \,,
\end{aligned}
\end{equation}
for all $l_1\ino \overleftarrow{\bbZ} $, $l_2 \ino \bbZ$ such that $l_1\leq l_2$. To simplify,  we also set 
\begin{equation}
\label{simplifff}
\forall l\ino \bbZ, \quad u_{| l }=  u_{|(-\infty , l] } . 
\end{equation}
Let us stress that a bilateral word has at most a finite number of letters $a_k$ indexed by positive indices, $k\ino \mathbb N^*$, while it can have infinitely many letters indexed by negative indices, $k\ino \overleftarrow{\bbZ}\cap \bbZ_-$. {Note that $\varnothing$ has neither relative height nor depth. Note that if $l_1\! \leq \! k_1$ and $k_2 \leq l_2$, then $u |_{(l_1, l_2]}\! = \! u$.}

\smallskip

\noi
$\bullet$ \textit{Shift.}  For all $l \in \bbZ$, we denote the $l$-\textit{shift operator} $\varphi_{l}\! : \! \overline{\bbW}_{\! A} \! \rightarrow \!  \overline{\bbW}_{\! A}$ by:  
\begin{equation*}  
\varphi_{l} (u) = (a_{k+l})_{k_1 -l < k\leq k_2 -l }, 
\end{equation*}
Note that $|\varphi_{l} (u) |_- \! = \! |u|_- \! - l$, that $|\varphi_{l} (u)|\! = \! |u| - l$ and that 
$\varphi_l (\varnothing)\! = \! \varnothing$. Clearly, $\varphi_l \circ  \varphi_{l^\prime} \! = \! \varphi_{l+ l^\prime}$, $\varphi_l$ is bijective and $\varphi_{l} \circ  \varphi_{-l}$ is the identity map.

\smallskip

\noi
$\bullet$ \textit{Concatenation.}  For all $u\! =\!  (a_k)_{k_1 <k\leq k_2} \ino \overline{\bbW}_{\! A}$ and all $v\! = \! (b_k)_{1\leq k\leq k_3} \ino \bbW_{\! A}$, we define 
 \begin{equation} 
\label{concbila}
u \! \ast \! v\!  = \! (c_l)_{k_1< k\leq k_2+k_3} \quad \textrm{where} \quad c_k 
= \left\{ 
\begin{array}{ll}
a_k &  \textrm{If $k_1 \! < \! k \! \leq \! k_2$,} \\
b_{k-k_2} &  \textrm{If $k_2 \! < \! k \! \leq \! k_2+ k_3$.}
\end{array} \right.
\end{equation}
The bilateral word $u \! \ast \! v$ is the concatenation of a bilateral word $u$ on the left with 
a null depth word $v$ on the right. Note that $|u\ast v|_-\! = \! |u|_-$, that  
$|u\ast v|\! = \! |u|+ |v|$ and that $u\ast \varnothing\! = \! u$. 

\smallskip

\noi
$\bullet$ \textit{Convergence in $\overline{\bbW}_{\! A}$.} Assume that $(A, d_A)$ is a Polish space. 
We equip $\overline{\bbW}_{\! A}$ with the following local convergence. 

\smallskip

 \begin{compactenum}
\item[] \textit{Let $u^{_{(p)}}_{^{\! }}\! \ino \overline{\bbW}_{\! A}$, $p\ino \bbN$; the sequence of words $u^{_{(p)}}_{^{\! }}$ converges to $u$ if $|u^{_{(p)}}_{^{\! }}|_-\! \rightarrow \! |u|_-$ in $\overleftarrow{\bbZ}$
and if for all $l\ino \bbN$ and for all $\epp \ino (0, 1)$, there exists $p_{l, \epp}\ino \bbN$ such that for all $p\! \geq \! p_{l, \epp} $, $|u^{_{(p)}}_{^{\! }} \! |\! = \! |u|$, $(-l)\! \vee \! |u^{_{(p)}}_{^{\!  }} \! |_-\! = \! (-l) \! \vee \! |u|_-\, $ and   
$\, \max_{(\! -l)  \vee |u|_- \! \leq k\leq |u| } d_A \big(u^{_{(p)}}_{^{k}} \! , u^{_{\!}}_ {^{k}} \big) < \epp $. 
}

\smallskip

\end{compactenum}
It is easy to see that this convergence corresponds to a Polish metric and we equip $\overline{\bbW}_{\! A}$ with the corresponding Borel sigma-field. 
Note that the shifts operators $\varphi_l$ are homeomorphisms with respect to local convergence.    
If $A \! = \! \bbN^*$, we shall use the notation $\overline{\bbU} \! := \! \overline{\bbW}_{\bbN^*}$. 
If $A\! = \! \{ 1, \ldots, \bgg\}$, we shall use the notation $\overline{\bbW}_\bgg \! := \! \overline{\bbW}_{\! \{ 1, \ldots, \bgg\}}$. Note that in these  cases, $\bbN^*$ and $\{ 1, \ldots, \bgg\}$ are equipped with the discrete topology.

\paragraph*{Trees with infinite line of ancestors.} 
\begin{definition}
\label{Wsubtree} A non-empty subset $R \! \subset \! \overline{\bbW}_{\! A}$ is a subtree of $\overline{\bbW}_{\! A}$ if it satisfies the following. 

\smallskip

\begin{compactenum}

\item[$(a)$] There exists $|R|_-\ino \overleftarrow{\bbZ}$ such that $|v|_-\! = \! |R|_-$ for all $v\ino R\backslash \{ \varnothing\}$.

\smallskip

\item[$(b)$] For all $v\ino R \backslash \{ \varnothing\}$ 
such that $|v| \! >\! |R|_-$, $\overleftarrow{v}\ino R$.

\smallskip

\item[$(c)$] For all $u$ and $v$ in $R$, there exists $l\ino \bbZ$ such that $u_{| l} \! = \! v_{| l}$ (see \eqref{simplifff}).   \cq 
\end{compactenum}
\end{definition}
We call $|R|_-$ the \textit{depth of the subtree $R$}.  
If $|R|_-\! >\! -\infty$, then $(b)$ implies that the empty word $\varnothing$ is an element of $R$ 
 and $|R|_-+1\! = \! \min_{v\in R \backslash \{ \varnothing \}} |v|$, 
 $(c)$ is always fulfilled since for all $l \! \leq \! |R|_-$ and all $v\in R$, we get $v_{| l}\! = \! \varnothing$. If $|R|_-\! =\! -\infty$, then $\varnothing \! \notin \! R$. 
 
\smallskip 
 
\noi
$\bullet$ \textit{Common ancestor.}
Recall from (\ref{simplifff}) the notation $u_{| l}$. Let $R\! \subset \! \overline{\bbW}_{\! A}$ be a subtree as in Definition \ref{Wsubtree}.  We define the \textit{common ancestor} of 
$u, v\ino R$ by 
\begin{equation}  
\label{comancdef}
u\wedge v  = \! u_{| b(u,v)} \!\! =  v_{| b(u,v)} \quad \textrm{where} \quad  
b(u,v)\! :=\! \max \{ l\ino \bbZ : \, u_{| l} \!\! = \!  v_{| l} \}, 
\end{equation}
that is well-defined thanks to Definition \ref{Wsubtree} $(c)$.

\smallskip 
 
\noi
$\bullet$ \textit{Graph distance on subtrees of $\overline{\bbW}_{\! A}$.}
A subtree $R \! \subset \! \overline{\bbW}_{\! A}$ as in Definition \ref{Wsubtree} corresponds to the following graph-tree: its set of vertices is $R$ and its set of edges is $\{ \{ v, \overleftarrow{v}\}; v\ino R  :  |v| \! >\! |R|_- \}$. We easily observe that the graph distance $d_{\mathtt{gr}}$ on $R$ is given by 
\begin{equation}
\label{graphsubt}
\forall x,y\ino R, \quad d_{\mathtt{gr}}(x,y)= |x| + |y| -2 |x\wedge y| \; .
\end{equation}

\smallskip

\smallskip 
 
\noi
$\bullet$ \textit{Ordered trees with a possibly infinite line of ancestors.} 
We next extend Definition \ref{ort} to ordered trees with a possibly infinite line of ancestors as follows. 
\begin{definition}
\label{deftreeee} Recall the notation $\overline{\bbU}\! := \! \overline{\bbW}_{\bbN^*}\! $. 
A subset $t \! \subset \! \overline{\bbU}$ is an ordered tree (with a possibly infinite line of ancestors) if it is a subtree of $\overline{\bbU}$ as in Definition \ref{Wsubtree} satisfying 
$(a)$, $(b)$ and $(c)$ with $A\! = \! \bbN^*$ 
and if it furthermore satisfies the condition that any word has a finite number of children, that is, 

\smallskip

\noi
$(d)\, $  $\forall u\ino t$, $k_u (t)\!\!  : = \! \# \{ v\ino t\! : \! \overleftarrow{v}\! = \! u\} \! < \! \infty$ and $\{ 1, \ldots , k_u (t)\} \! =\!  \{ \mathtt{end} (v); v\ino t \! \! : \!   \overleftarrow{v}\! = \! u\}$ if $k_u (t) \! \geq \! 1$.

\smallskip

\noi 
We shall consider that the singleton $\{ \varnothing\}$ is the only tree with one point. 
We denote by $\overline{\bbT}$ the set of ordered trees. \cq  
\end{definition}

\noi
For all $k\ino \overleftarrow{\bbZ}$, 
we set $\bbT_k \! =\! \{ \{ \varnothing \} \} \cup \{ t\ino \overline{\bbT}: |t|_- \! = \! k \}$. Note that $\bbT_0\! = \! \bbT$, where $\bbT$ is as in Definition \ref{ort} and observe that $\varphi_{-k} (\bbT) \! = \! \bbT_k$ {when $k>-\infty$}.

\smallskip

\noi
$\bullet$ \textit{Lexicographical order and successor of a vertex.} Let $t\ino \overline{\bbT}$. By Definition \ref{Wsubtree} $(c)$, the vertices of $t$ are totally ordered by the \textit{lexicographical order} $\leq_{t}$ that is formally defined as follows. Let $u,v\ino t\backslash \{ \varnothing\}$; recall from (\ref{comancdef}) the definition of $b(u,v)$ and from (\ref{extremdef}) the definition of $\mathtt{end} (\cdot)$; then,   
\begin{equation}
\label{lextredef}
u \leq_{t} v  \quad \textrm{{if and only if}} \quad \mathtt{end} (u_{| b(u,v) +1}) \leq  \mathtt{end} (v_{| b(u,v) +1})\; .
\end{equation}
Note that $\leq_{t}$ actually depends on $t$: it is not defined on the whole set of bilateral words $\overline{\bbU}$ but only on $t$ (indeed, to define $\leq_{t}$, branching points have to be well-defined, which requires possibly infinite words to share a prefix). If $\varnothing \ino t$, then $\varnothing$ is the $\leq_{t}$-least element of $t$. We denote by $<_t $ the strict order associated with $\leq_{t}$. We also introduce the following related notation: for all $u \ino t$, the \textit{successor} $\mathtt{scc} (u)$ of $u$ is defined as the $\leq_{t}$-least element of $\{ v \ino t \! : \! u \! <_{t} \! v\}$ if this set is not empty, otherwise we simply take $\mathtt{scc} (u)\! = \! u$.

\smallskip

\noi
$\bullet$ \textit{Subtree.} Let  $t\ino \overline{\bbT}$ {and $u\in t$}. The subtree $\theta_u t$ stemming from $u$ is defined as follows. 
\begin{equation}
\label{curshitr} 
\theta_u t\! = \! \{ v\ino \bbU: u\ast v \ino t \} , 
\end{equation}
where we recall from (\ref{concbila}) the definition of the concatenation $\ast$ of a bilateral word on the left with a null-depth word on the right. Note that $|\theta_u t|_-\! = \! 0$, namely: $\theta_u t \ino \bbT$, where $\bbT$ is as in Definition \ref{ort}.

\vspace{-3mm}

\paragraph*{Pointed labelled trees.} To deal with branching random walks, we introduce labelled trees where the label of a vertex is viewed as its position in space. More precisely, let $(E, d_E)$ be a Polish metric space. 
We define the space of \textit{pointed $E$-labelled trees} as follows. For all $k\ino \overleftarrow{\bbZ}$, we set: 
$$ \bbT^\bullet_k  (E) \! = \! \{ \partial \} \! \cup \! \big\{  \bt\! =\!  \big(t, \boo \, ; \bx \! =\!  (x_v)_{v\in t} \big)\! : t\ino \bbT_k, \, \boo\ino t, \, x_v\ino E, v\ino t  \big\} \quad  \textrm{and} \quad \overline{\bbT}^\bullet\!  (E)\! =  \!\! \bigcup_{ k\in   \overleftarrow{\bbZ}} \! \bbT^\bullet_k  (E), $$
where $\partial$ stands for a cemetery point. 
Here, the label of $v\ino t$ is $x_v\ino E$ that is viewed as the spatial position of $v$. 
If there is no label, we simply write $\bbT^\bullet_k$ and $\overline{\bbT}^\bullet$.

\smallskip

\noi
$\bullet$ \textit{Shift operator on labelled trees.} Shift operators act naturally on the space of pointed $E$-labelled trees as follows: 
let $l\ino \bbZ$; we set  $\varphi_l (\partial)\! = \! \partial$ and for all 
$\bt\! = \! (t,\boo; \bx)\! \in \overline{\bbT}^\bullet\!  (E)\backslash \{\partial \}$, 
we set 
\begin{equation*}
\varphi_l (\bt)\! = \! (\varphi_l (t),\varphi_l (\boo); \bx \!  \circ \!  \varphi_{-l}  ), \quad \textrm{where} \quad 
 \bx \!  \circ \!  \varphi_{\! -l} \! = \! (x_{\varphi_{\! -l} (w)})_{w\in \varphi_l (t)}. 
 \end{equation*}

\noi
$\bullet$ \textit{Truncation.} We next define a natural truncation procedure for pointed labelled trees along the line of ancestors of the distinguished point.  

\vspace{-3mm}

\begin{definition}
\label{truncdef}
Let $ \bt\! =\!  \big(t, \boo \, ; \bx \! =\!  (x_v)_{v\in t} \big) \ino \overline{\bbT}^\bullet \!   (E)$. 
Let $p\ino \overleftarrow{\bbZ}$ and $q\ino \bbZ \! \cup \! \{ \infty \}$ be such that $p\! < \!  q$. We define the following.  
\begin{compactenum}

\smallskip

\item[$-$] If $(p,q] \cap  \big( |\boo|_{_{\! -}}, |\boo| \big]  \! = \! \emptyset$, then we set $[\bt ]^q_p\! = \!  \partial$.  
We also set $[\partial ]^q_p\! = \!  \partial$.

\smallskip

\item[$-$]If $(p,q] \cap  \big( |\boo|_{_{\! -}}, |\boo| \big]  \! \neq \! \emptyset$, 
then we set $[\bt]_p^q\! = \!  (t^\prime, \boo^\prime ; \bx^\prime)$
where 
$$ \boo^\prime \! =\!  \boo_{| ( p, q]} , \quad t^\prime \! = \! \big\{  v_{| (p,q]};\,  v\ino t :  v_{| p}\! =
\! \boo _{| p} \big\} \; \textrm{and} \;  x^\prime_{v^\prime}\! = \! x_v, $$
where $v\ino t$ is
such that $v^\prime \! = \!  v_{| (p,q]}$ and $v_{| p} \! =
\! \boo_{| p } $ (recall notation $v_{| p}$ from (\ref{simplifff})). 
\end{compactenum}

\smallskip

\noi
We simply set $[\bt]_p$ instead of $[\bt]_p^\infty$. If $\bt \ino  \bbT^\bullet_{k}  (E) $, then note that $[\bt]_p^q\ino  \bbT^\bullet_{k\vee p}  (E)$. We use a similar notation for pointed trees without label.  \cq 
\end{definition}

\noi
$\bullet$ \textit{Local convergence on $\overline{\bbT}^\bullet \!   (E)$.}
For all $\bt \! = \! (t,\boo; \bx)$, $\bt^\prime \!= \! (t^\prime, \boo^\prime; \bx^\prime)$ in 
$\overline{\bbT}^\bullet \! (E)$, we first set 
$$ \Delta (\bt, \bt^\prime)\! = \! \un_{\{ (t, \boo)\neq (t^\prime, \boo^\prime)\}} + 
\un_{\{ (t, \boo)= (t^\prime, \boo^\prime) \}} \max_{v\in t} \big(1\! \wedge \! d_E (x_v, x^\prime_v) \big) $$
with $\Delta (\bt, \partial)\! = \! 1$ and $\Delta (\partial, \partial)\! = \! 0$. We easily check that $\Delta$ is a metric on $\overline{\bbT}^\bullet \! (E)$. Then, we define the local convergence as follows. 

\smallskip

\begin{compactenum}
\item[] \textit{Let $\bt_n\! = \! (t^{(n)}\! , \boo^{(n)} ; \bx^{(n)}) \ino \overline{\bbT}^\bullet\! (E)$, $n\ino \bbN$; the sequence $(\bt_n)_{n\geq 0}$ is said to \textit{converge locally} to $\bt\! = \! (t, \boo; \bx) \ino  
 \overline{\bbT}^\bullet\!  (E)$ if for all $\epp \ino (0, 1)$ and all $q\ino \bbN$, there exists 
 $n_{q, \epp} \ino \bbN$ such that for all integers $n\! \geq \! n_{q, \epp}$, 
 $\Delta \big( [\bt^{(n)}]_{-q}^q,  [\bt]_{-q}^q \big) \! < \! \epp$ (when there is no label, it simply means that $[\bt^{(n)}]_{-q}^q \! =\!   [\bt]_{-q}^q\, $).}
\end{compactenum}

\smallskip

\noi
Local convergence corresponds for instance to the following metric: 
\begin{equation*}
\forall \bt, \bt^\prime \ino  \bbT^\bullet\! (E) , \quad \bdelta_{\mathtt{loc}} \big(\bt, \bt^\prime \big)= \sum_{q\in \bbN} 2^{-q-1} \Delta \big(   [\bt]_{-q}^q, [\bt^{\prime}]_{-q}^q \big), 
\end{equation*}
with $\bdelta_{\mathtt{loc}} (\bt, \partial)\! = \! 1$ and $\bdelta_{\mathtt{loc}}(\partial, \partial)\! = \! 0$. We easily check that $(\bbT^\bullet\! (E), \bdelta_{\mathtt{loc}})$ is Polish and we note that shift operators are isometries. 

\vspace{-3mm}

\begin{definition}
\label{ggloupin} Let $(t,\boo)\ino \overline{\bbT}^\bullet $ be distinct from $\partial$. 
\begin{compactenum}

\smallskip

\item[$(a)$] We set $\mathtt{cent} (t,\boo)\! = \! \varphi_{|\boo|} (t, \boo)$ that is the 
\textit{centering} map: it shifts trees so that their distinguished point is at relative height $0$.

\smallskip

\item[$(b)$] We next set $\mathtt{scc} (t,\boo) \! = \! \mathtt{cent} (t, \mathtt{scc} (\boo))$ where we recall that 
$ \mathtt{scc} (\boo)$ stands for the vertex of $t$ coming next in the lexicographical order as defined by (\ref{lextredef}). 
We call $\mathtt{scc} (\cdot)$ the \textit{successor map}.

\smallskip

\item[$(c)$] Observe that there is a unique pointed tree $(t^\prime, \boo^\prime)\ino \overline{\bbT}^\bullet$ and a unique one-to-one map $\psi \! : \! t^\prime \! \rightarrow \! \{\boo_{| p}; p\! \leq \! |\boo| \} \cup  \{ v\ino t\! : \!  \boo \! <_t \! v \}$ such that $\psi (\varrho^\prime) \! = \! \varrho$, that is increasing with respect to the lexicographical order and 
that preserves the relative height; we set $[(t, \boo ) ]^+\! = \! (t^\prime, \boo^\prime )$ 
that is called the \textit{right-part} of $(t,\boo)$. {Intuitively, the difference between $(t',\rho')$ and $\{\boo_{| p}; p\! \leq \! |\boo| \} \cup  \{ v\ino t\! : \!  \boo \! <_t \! v \}$ is that $(t',\rho')$ respects the convention $(c)$ in Definition \ref{ort} we have imposed on rooted ordered trees.} We also set 
\begin{equation*}
\mathtt{scc}^+ (t,\boo) \! = \! [\mathtt{scc} (t,\boo) ]^+\; .
\end{equation*} 
The map $\mathtt{scc}^+ (\cdot)$ is called the \textit{right-successor}. 
\end{compactenum}

\smallskip

\noi
By convenience, we set $\mathtt{cent} (\partial)\! =\!  \mathtt{scc} (\partial) \! = \! \mathtt{scc}^+ (\partial)\! = \! \partial$. \cq 

\end{definition}

\noi
Note that 
\begin{equation}
\label{cikdnd}
[\mathtt{cent} (t,\boo)]^+\! = \! \mathtt{cent} \big( [(t, \boo)]^+\big) \quad \textrm{and} \quad \mathtt{scc}^+ (t, \boo)\! = \! \mathtt{scc}^+\!  \big( [(t, \boo)]^+\big)\; .
\end{equation}
Let us state a technical result about the continuity of the maps $\mathtt{cent}, [\cdot]^+, \mathtt{scc}$ and $\mathtt{scc}^+$.

\begin{lemma}
\label{consucc} The maps $\mathtt{cent} (\cdot)$ and $[\, \cdot \, ]^+$ are 
locally continuous and the maps $\mathtt{scc} (\cdot) $ and  $\mathtt{scc}^+ (\cdot) $ are locally continuous at the pointed trees $(t,\boo)$ such that $\mathtt{scc} (\boo) \! \neq \! \boo$. 
\end{lemma}
\noi
\textbf{Proof.} Since shift-operators are $\bdelta_{\mathtt{loc}}$-isometries, 
the continuity of $\mathtt{cent}$ follows from the continuity of 
$(t,\boo)\! \rightarrow \! |\boo|$ that is a direct consequence of the definition of local convergence; $[\, \cdot\, ]^+$ is locally continuous because $[[\bt]_{-q}^q]^+\!\! = \! [[\bt]^+]_{-q}^q $ and to complete the proof, it is then sufficient to prove that $(t,\boo) \!  \mapsto \!(t,\mathtt{scc}(\boo))$ is locally continuous at trees such that $\mathtt{scc} (\varrho) \! \neq \! \varrho$. 

To that end, let $(t_n, \boo_n)\! \rightarrow \! (t, \boo)$ locally in $ \overline{\bbT}^\bullet\!\! $. Set $p_0 \! = \! |\mathtt{scc} (\boo)|$ and suppose  that $\boo\! \neq \! \mathtt{scc} (\boo)$, which implies $|\boo \! \wedge \! \mathtt{scc} (\boo)|\! = \! p_0\! -\! 1$. Let $p\ino \bbN$ be such that 
$-p\! < \! p_0\! -\! 1 \! <\! |\boo|+1 \! < \! p$; there is $n_p\ino \bbN$ such that for all $n\! \geq \! n_p$, $(t^\prime_n, \boo_n^\prime)\! := \! [(t_n, \boo_n)]_{^{-p}}^{_{p}}$ is equal to $(t^\prime, \boo^\prime)\! := \! [(t, \boo)]_{^{-p}}^{_{p}}$; thus, the successor of  $\boo^\prime_n$ in 
$t^\prime_n$ is equal to the successor of $\boo^\prime$ in $t^\prime$. Consequently, 
$[(t_n, \mathtt{scc} (\boo_n))]_{^{-p}}^{_{p}}\! = \! (t^\prime_n, \mathtt{scc} (\boo_n^\prime))\! = \! (t^\prime, \mathtt{scc} (\boo^\prime))\! = \! [(t, \mathtt{scc} (\boo))]_{^{-p}}^{_{p}}$, which entails the desired result. \cqfd

\vspace{-3mm}

\paragraph*{Infinite pointed Galton--Watson trees.}
For all $t\ino \overline{\bbT}$ and for all $u\ino t$, recall from (\ref{curshitr}) the definition of the subtree $\theta_u t \ino \bbT$.  

\vspace{-3mm}

\begin{definition}
\label{GWptdef} Let $(r(j,k))_{k\geq j\geq 1}$ be a probability measure 
on the octant $\{ (j,k)\ino (\bbN^*)^2\! : \! j\! \leq \! k\}$. 
Let $\mu$ be a probability measure on $\bbN$. 
Let $\ftau^*\! = \! (\tau^*\! , \varrho)\! :\!  \Omega \! \rightarrow \!   \overline{\bbT}^\bullet\! $ be a 
Borel-measurable random pointed tree such that a.s.~$|\varrho|\! = \! 0$ and $|\tau^*|_-\!\! = \! -\infty$. We introduce the following notation. 
$$\forall p\ino \bbN, \quad  \booo (p) = \varrho_{| \, ]-\infty , -p]}, 
\quad \mathtt{Sp} \! = \! \{  \varrho (p) \, ; p\ino \bbN\} \; \, \textrm{and} \; \,
\partial \mathtt{Sp}\! = \! \big\{ u\ino \tau^* \backslash \mathtt{Sp}: \overleftarrow{u}\ino \mathtt{Sp} \big\} .$$  
Then, $\ftau^*$ is an \textit{infinite pointed Galton--Watson tree with offspring distribution $\mu$ and dispatching measure $r$} if the 
r.v.~$S\! :=\! \big( \mathtt{end} (\varrho (p)),  k_{\varrho (p +1)} (\tau^*)\big)_{\! p\in \bbN}$ are i.i.d.~with law $r$ and if conditionally given $S$, 
the subtrees $(\theta_{u} \tau \, ; \, u\ino \partial \mathtt{Sp})$ are independent GW($\mu$)-trees.

We shall deal with the following special cases that are well-defined if $m_\mu \! := \! \sum_{k\in \bbN} k\mu (k) \! < \! \infty$.
\begin{compactenum}

\smallskip

\item[$(i)$] If $r(j,k)\! = \! \mu (k)/m_\mu$, for all $k\! \geq \! j\! \geq \! 1$, then we say that $\ftau^*$ is an \textit{infinite pointed GW($\mu$)-tree} (an IPGW($\mu$)-tree for short). 

\smallskip

\item[$(ii)$] If $r(j,k)\! = \! \un_{\{ j=1\}}\overline{\mu} (k)$, for all $k\! \geq \! j\! \geq \! 1$, 
where $\overline{\mu} (k)\! = \! m_{\mu}^{-1}\! \sum_{l\geq k}\mu (l)$ for all $k\! \geq \! 1$, then 
 we say that $\ftau^*$ is the \textit{right part of an infinite pointed GW($\mu$)-tree} (an IPGW$^+ $($\mu$)-tree for short). 
\end{compactenum}

\smallskip
\noi
Note that if $\ftau^*$ is an IPGW($\mu$)-tree, then $[\ftau^*]^+\! $ is an IPGW$^+ $($\mu$)-tree. \cq 
\end{definition}
IPGW-trees are related to GW-trees via the \textit{many-to-one} principle (or the \textit{one-point decomposition of GW-trees}) that asserts 
the following: let $\mu$ be a probability distribution on $\bbN$ such that $m_\mu \! : = \! \sum_{k\in \bbN} k\mu (k)\! < \! \infty$. 
Let $\tau$ be a GW($\mu$)-tree and let $\ftau^*$ be an IPGW($\mu$)-tree as in Definition \ref{GWptdef}. Then for all Borel-measurable functions $F \! : \! \bbN \times \overline{\bbT}^\bullet \! \rightarrow \! [0, \infty)$,   
\begin{equation}
\label{sizebiased}
\bE \Big[\sum_{v\in \tau} F \big( |v| ; \varphi_{|v|} (\tau, v)  \big) \Big]= \sum_{p\geq 0} m_\mu^p\, \bE\big[ F (p\, ;  [\ftau^*]_{-p})\big] \; , 
\end{equation}
where we recall from SL{Definition} \ref{truncdef}  the notation $ [\ftau^*]_{-p}$ for the pointed 
tree $\ftau^*$ truncated above the ancestor of $\booo$ at generation $-p$. 
Based on this identity, the following proposition  
shows that IPGW trees are local limits of critical GW-trees conditioned to be large and seen from a uniformly chosen vertex. This result is part of the folklore; its proof derives from (\ref{sizebiased}) and it is left to the reader. 

\vspace{-2mm}

\begin{proposition}
\label{GWptGWloc} Let $\mu$ be a probability distribution on $\bbN$ such that $\sum_{k\in \bbN} k\mu (k)\! = \! 1$. We assume that $\mu $ is aperiodic. 
Let $\tau$ be a GW($\mu$)-tree; let $\bu$ be uniformly distributed on the set of vertices of $\tau$. Let $\ftau^*$ be an IPGW($\mu$)-tree as in Definition \ref{GWptdef}. 
Then 
$$\mathtt{cent} (\tau, \bu) \quad  \textrm{under} \quad 
\bP (\, \cdot \, | \, \# \tau \! = \! n) \;  \underset{n\rightarrow \infty}{-\!\!\! -\!\!\! \longrightarrow} \; \ftau^* $$

\vspace{-2mm}

\noi
weakly on $\overline{\bbT}^\bullet$ with respect to local convergence. 
\end{proposition}

\vspace{-1mm}

\noi
We use the previous proposition to prove the following one. 

\vspace{-2mm}

\begin{proposition}
\label{succinv}  Let $\mu$ be a probability measure on $\bbN$ such that $\sum_{k\in \bbN} k\mu (k)\! = \! 1$. Recall $\mathtt{scc} (\cdot)$ 
and $\mathtt{scc}^+ (\cdot)$ from {Definition \ref{ggloupin}}. Then, the law of IPGW($\mu$)-trees (resp.~IPGW$^+$($\mu$)-trees) is preserved by $\mathtt{scc} (\cdot)$ (resp.~by $\mathtt{scc}^+ (\cdot)$).
\end{proposition}

\vspace{-2mm}

\noi
\textbf{Proof.} Let us first mention that a different proof of the result for $\mathtt{scc}^+ (\cdot)$ is given in Proposition 2 in Le Gall and L.~\cite{LGLi16}.  Then, note that the result for $\mathtt{scc}^+ (\cdot)$ is implied by the result for $\mathtt{scc} (\cdot)$ by (\ref{cikdnd}) and since the {right-part} of a IPGW($\mu$)-tree is an IPGW$^+ $($\mu$)-tree.
Let us prove the result for $\mathtt{scc} (\cdot)$. Let $\ftau^*\!= \! (\tau^*\! , \varrho)$ be an IPGW($\mu$)-tree and let $\tau$ be a GW($\mu$)-tree; let $\bu$ be uniformly distributed on the set 
of vertices of $\tau$; denote by $v_*$ the last vertex of $\tau$ with respect to the lexicographical order. Set $\bu^\prime\! = \! \mathtt{scc} (\bu)$ if $\bu \! \neq \! v_*$ and $\bu^\prime\! = \! \varnothing$ if $\bu\! = \! v_*$: clearly, $\bu^\prime$ is uniformly distributed on $\tau$ and $\bP ( \mathtt{scc} (\bu)\! \neq \! \bu^\prime | \, \# \tau\! = \! n)\! = \! 1/ n$.   
Let $F\! :\! \overline{\bbT}^\bullet \! \rightarrow \! [0, \infty)$ be locally continuous and bounded. By Definition \ref{ggloupin} $(b)$, 
$\mathtt{scc} (\tau, \bu)\! =\! \mathtt{cent} (\tau, \mathtt{scc}(\bu))$. Thus, 
\[
\big| \, \bE [ F(\mathtt{scc} (\tau, \bu)) | \# \tau \! =\!  n \big] \! -\!  \bE \big[ F( \mathtt{cent} (\tau, \bu^\prime))|\# \tau \! = \! n  \big] \, \big| \!  \leq \! 2 \lVert F \rVert_\infty /n.
\]
First suppose that $\mu$ is aperiodic. 
Then the previous inequality and Proposition \ref{GWptGWloc} imply that 
\[
\lim_{n\rightarrow \infty}\bE \big[ F(\mathtt{scc} (\tau, \bu)) |\# \tau \! = \! n  \big]\! = \! \lim_{n\rightarrow \infty}\bE \big[ F( \mathtt{cent} (\tau, \bu^\prime))|\# \tau \! = \! n  \big]\! = \! \bE [ F( \ftau^*)].
\]
Moreover, since it is clear that a.s.~$\mathtt{scc} (\varrho) \! \neq \! \varrho $, 
Lemma \ref{consucc} and Proposition \ref{GWptGWloc} entail 
\[
\lim_{n\rightarrow \infty}\bE \big[ F( \mathtt{scc} (\tau, \bu))|\# \tau \! = \! n  \big]\! = \! \bE [ F(\mathtt{scc} (\ftau^*))],
\]
which completes the proof when $\mu$ is aperiodic. 

 Let us consider a general $\mu$. For all $\epp \ino (0, 1)$, we set $\mu_\epp\! := \! \epp \delta_1 + (1\! -\! \epp) \mu$. Namely, $\mu_\epp$ is a critical aperiodic offspring distribution. Let $\ftau^*_\epp$ be an infinite GW($\mu_\epp$)-tree. We easily check that $\ftau^*_\epp\! \rightarrow \! \ftau^*$ locally as $\varepsilon \! \rightarrow \! 0$. The local continuity of $\mathtt{scc}(\cdot) $ (Lemma \ref{consucc}) entails the desired result. \cqfd

\section{Tree-valued branching random walks.}
\label{treeBRWsec}

\subsection{Definitions.}
\label{trbrwdfsc}
Let $E$ be a (Polish) space of labels. Let $(q(y,dy^\prime))_{ y\in E}$ be a transition kernel and 
let $\varpi$ be a Borel probability measure on $E$. 
For all pointed tree $ \bt\! =\! (t, u)\ino \overline{\bbT}^\bullet\! \! $, we define the law $Q_{\varpi, \bt}$ on $\overline{\bbT}^\bullet \! (E)$ of the $q$-branching random walk with genealogical tree $t$ and such that $\varpi$ is the law of the spatial position of the distinguished individual $u$. 

To that end, we first assume that $t$ is finite and we introduce the following notation: for any 
$v, w \ino t$, we denote by $\lgeo v, w\rgeo$ the set of vertices on the shortest path joining $v$ to $w$ in the tree $t$. If $w\! \neq \! v$, then we denote by $\overleftarrow{w}^v$ the unique $v^\prime\ino \lgeo v, w \rgeo $ that is at graph-distance $1$ from $w$; 
we call $\overleftarrow{w}^v$ the \textit{$v$-parent of $w$}.

Then, the r.v.~$\fTheta\! =\! (t,u ; (Y_v)_{v\in \tau})\! :  \Omega \! \rightarrow \!  \overline{\bbT}^\bullet\!  (E)$ has law 
$Q_{\varpi, \bt}$ if 
\begin{equation}
\label{Ybbs}
\textrm{
the joint law of $(Y_v)_{v\in t}$ is} \quad \,    \varpi (dy_u) \!\! \!  \prod_{v\in t\backslash \{ u\}} \!\!\! q( y_{\overleftarrow{v}^u}, dy_v),
\end{equation}
The definition \eqref{Ybbs} can be extended to the case in which $\bt$ is infinite. For such purpose it is enough to note that 
$[\fTheta]_{p}^{q}$ has law $Q_{{\varpi},{[(t,\varrho)]_p^q}}$ 
(see Definition \ref{truncdef} for the truncation $[\cdot ]_p^q$). 

\smallskip
 
Let us note that 
\begin{equation*} 
\bt_n \underset{n\rightarrow \infty}{-\!\!\! -\!\!\! \longrightarrow} \bt \; \, \textrm{locally in $\overline{\bbT}^\bullet$} \quad \Longrightarrow \quad 
 Q_{\varpi, \bt_n}  \underset{n\rightarrow \infty}{-\!\!\! -\!\!\! \longrightarrow} Q_{\varpi,\bt}  
\; \, \textrm{weakly on $\overline{\bbT}^\bullet \! (E)$.}
\end{equation*}

 When $\varpi\! = \! \delta_y$ for some $y\ino E$, we simply write 
$Q_{y, \bt}$ instead of $Q_{\delta_y, \bt}$. Since $y\mapsto q(y, A)$ is Borel-measurable for all Borel subsets $A$ of $E$, it is easy to check that $y\mapsto  Q_{y, \bt}$ 
is also Borel measurable and that $Q_{\varpi, \bt}\! = \! \int_E \! \varpi (dy) \,  Q_{y, \bt}$.  

As an immediate consequence of the definition, we also get the following: fix $l \ino \bbZ$ and set $(t^\prime, u^\prime) \! = \! 
 (\varphi_l (t), \varphi_l (u))$. 
 \begin{equation}  
\label{atschift}
\textrm{If $\fTheta= (t,u; (Y_v)_{v\in t})$ has law $Q_{\varpi, (t,u)}$, 
then $\big( t^\prime, u^\prime ; 
(Y_{\varphi_{-l} (w)} )_{w\in t^\prime } \big)$ has law $Q_{\varpi, (t',u')} $.}
\end{equation}

\begin{definition}
\label{defimodel}
We fix $\bt\! = \! (t,u) \ino \overline{\bbT}^\bullet$. We shall consider mostly the four following cases. 
\begin{compactenum}

\smallskip

\item[$(i)$] $E \! = \! \bbW_{\! \bgg}$, the $\bgg$-ary tree equipped with the local convergence; $\varpi \! = \! 
\delta_{\varnothing}$ and $q(x,dy)\! =\!  p^+_\bgg (x,dy)$ where for all measurable $f: \bbW_{\! \bgg} \! \rightarrow \! [0, \infty)$, 
\begin{equation*} 
\int_{\bbW_{\! \bgg}} \!\!  \!\! p_\bgg^+ (x, dy) \, f(y) = \left\{ 
\begin{array}{ll}
\tfrac{1}{2} f(\overleftarrow{x}) + \tfrac{1}{2\bgg} \sum_{1\leq i\leq \bgg} \, f \big(  x\!  \ast \!  (i) \big) &  
\textrm{if $x\! \neq \! \varnothing$ } \\
\tfrac{1}{\bgg} \sum_{1\leq i\leq \bgg}f((i)) &  \textrm{if $x\! = \! \varnothing$.}
\end{array} \right.
\end{equation*}
We denote by $Q^{_{+\bgg}}_{^{\bt}}$ the law of the $\bbW_{\! \bgg}$-valued branching random walk with transition kernel $p_{^\bgg}^{_+}$ and with "initial" position $\varnothing$ in $\bbW_{\! \bgg}$.

\smallskip

\item[$(ii)$] $E \! = \! \bbW_{[0, 1]}$, the free tree equipped with the local convergence; $\varpi \! = \! 
\delta_{\varnothing}$ and $q(x,dy)\! =\!  p^+(x,dy)$ where for all measurable $f: \bbW_{[0, 1]} \! \rightarrow \! [0, \infty)$,
\begin{equation*} 
\int_{\bbW_{[0, 1]}}  \!\! \!\! \!\!  \!\! p^+ (x, dy) \, f(y) = \left\{ 
\begin{array}{ll}
\tfrac{1}{2} f(\overleftarrow{x}) + \tfrac{1}{2}\int_0^1 \! ds   \, f \big( x\!  \ast \!  (s) \big) &  
\textrm{if $x\! \neq \! \varnothing$ } \\
\int_0^1 \! ds   \, f \big( (s) \big)  &  \textrm{if $x\! = \! \varnothing$.}
\end{array} \right.
\end{equation*}
We denote by $Q^{_{+}}_{^{\bt}}$ the law of the $\bbW_{\! [0, 1]}$-valued branching random walk with transition kernel $p^{+}$ and with "initial" position $\varnothing$ in $\bbW_{\! [0, 1]}$.

\smallskip

\item[$(iii)$] $E \! = \! \Wid \! := \!  \{x\ino {\overline{\bbW}_\bgg}: |x|_{\! -} \! = \! -\infty \}$, equipped with the local convergence; we fix $\rmo \ino \Wid$ and we take 
$\varpi= \delta_{\rmo}$ and 
$q(x,dy)\! = \! p_\bgg (x, dy)$ where for all $x\ino { \Wid}$ and for all measurable $f:  \Wid \! \rightarrow \! [0, \infty)$, 
\begin{equation} 
\label{brwTinfd}
\int_{\Wid} \!\!  \!\! p_\bgg (x, dy) \, f(y) = 
\tfrac{1}{2} f(\overleftarrow{x}) + \tfrac{1}{2\bgg} \sum_{1\leq i\leq \bgg} \, f \big( x\!  \ast \!  (i) \big)  \; .
\end{equation}
We denote by $Q^{_{\bgg}}_{^{\rmo, \bt}}$ the law of the $ \Wid$-valued branching random walk with transition kernel $p_{\bgg}$ and with "initial" position $\rmo$ in $ \Wid$.

\smallskip

\item[$(iv)$] $E \! = \! \Wif\! := \{ x\ino \Wbf:  |x|_{\! -} \! = \! -\infty \}$, equipped with the local convergence; we fix $\rmo \ino \Wif$; we take 
$\varpi \! = \!  \delta_{\rmo}$ and 
$q(x,dy)\! =\!  p(x,dy)$, where for all $x\ino \overline{\bbW}_{[0, 1]}$ and for all measurable $f:\Wif\! \rightarrow \! [0, \infty)$, 
\begin{equation*} 
\int_{\Wif}  \!\! \!\! \!\!  \!\! p(x, dy) \, f(y) = \tfrac{1}{2} f(\overleftarrow{x}) + \tfrac{1}{2}\int_0^1 \! ds   \, f \big( x\!  \ast \!  (s) \big).  
\end{equation*}
We denote by $Q^{_{\! }}_{^{\rmo, \bt}}$ the law of the $\Wif$-valued branching random walk with transition kernel $p$ and with "initial" position $\rmo$ in $\Wif$: we shall refer to this branching random walk as the \textit{free branching random walk}. \cq 
\end{compactenum}
\end{definition}

\begin{rem}
Note that if $\fTheta\! = \! (t, u; (X_v)_{v\in t})$ has law $Q^{_\bgg}_{^{\rmo ,\bt}}$ 
(or $Q^{_{\!}}_{^{\rmo ,\bt}}$) then $|\fTheta|\! := \! (t, u; (|X_v|)_{v\in t})$ is a $\bbZ$-valued branching random walk whose spatial motion is that of the simple symmetric 
random walk. Similarly, if $\fTheta$ has law $Q^{_{+\bgg}}_{^{\bt}}$ (or $Q^{_+}_{^{\bt}}$), then $|\fTheta|$ is an $\bbN$-valued branching random walk whose spatial motion is that of the simple symmetric 
random walk reflected at $0$. \cq 
\end{rem}

\begin{definition}
\label{tracage}
We define the \textit{$\bgg$-contraction map} $\Phi_\bgg \! : \! \overline{\bbW}_{[0, 1]}\!  \rightarrow \!  \overline{\bbW}_{\! \bgg}$ as follows. For all $r\ino (0, \infty)$, we set 
$\lceil r \rceil \! =\! \min \{ k \ino \bbZ  \colon  r \! \leq \! k \}$ and by convenience we take $\lceil 0 \rceil \! = \! 1$.  
Then, for all 
$x \! = \! (a_k)_{|x|_- < k \leq |x|} \ino  \overline{\bbW}_{[0, 1]}$, we define 
\begin{equation*}
\Phi_\bgg (x)= \big(\lceil \bgg a_k \rceil)  \big)_{|x|_- < k \leq |x|} \in \, \overline{\bbW}_{\! \bgg} \; .
\end{equation*}
The $\bgg$-contraction map is measurable and  preserves the depth and the relative height of words.\cq 
\end{definition}
\begin{remark}
\label{ooobvvv}
Note that $\Phi_\bgg$ transforms respectively the kernel $p(x,dy)$ into $p_\bgg (x,dy)$ and the kernel $p^+(x,dy) $ into $p_{\bgg}^+(x, dy)$.   
It naturally extends to $\overline{\bbW}_{[0, 1]}$-labelled pointed trees as follows: if $\Theta \! = \! (t,\varrho ; (x_v)_{v\in t}) \ino \overline{\bbT}^\bullet \! (\overline{\bbW}_{[0, 1]})$, we set $\Phi_\bgg (\Theta)\! = \!  (t,\varrho ; (\Phi_\bgg (x_v))_{v\in t})$. Clearly, $\Phi_\bgg (\Theta)$ is a $\overline{\bbW}_{\! \bgg}$-labelled pointed tree and we easily check that the map is measurable. 
Moreover, if $\fTheta\! = \! (t, \varrho; (X_v)_{v\in t})$ has law $Q^{_+}_{^{\bt}}$ (resp.~$Q^{_{\!}}_{^{\rmo , \bt}}$) then 
$\Phi_\bgg (\fTheta)$ has law $Q^{_{+\bgg}}_{^{\bt}}$ (resp.~$Q^{_{\bgg}}_{^{\Phi_\bgg(\rmo), \bt}}$).  \cq 
\end{remark}

\subsection{Metric properties of the range of free branching random walks.}
\label{memesec}
  We gather basic facts 
about the range of a free branching random walk in terms of the heights of the spatial positions in $\Wif$. Let $ \bt\! = \! (t,u) \ino \overline{\bbT}^\bullet$ such that $|t|_{ -}\! = \! -\infty$.
For all $w\ino t$, we recall that $\lgeo v, w \rgeo$ is the shortest path (with respect to the graph-distance) that joins vertex $v$ to vertex $w$. Moreover, we set $\lgeo v, w \lgeo \, = \! \lgeo v, w \rgeo \backslash \{ w \! \} $, 
$\rgeo v , w \rgeo\! = \! \lgeo v, w \rgeo \backslash \{ v\!  \} $ and $\rgeo v, w \lgeo  = \! \lgeo v, v \rgeo 
\backslash \{ v, w \!  \}$; we also denote by $\rgeo\!  -\! \infty, v\rgeo$ the lineage of $v$: namely, 
$\rgeo \! -\! \infty, v\rgeo \! =\!  \{ v_{| l}\, ; l \! \leq \! |v| \}$.

We next decompose a free branching random walk by first describing the heights of the vertices as a $\bbZ$-valued branching random walk $(\mathbf h_v)_{v\in t}$ and then explaining how to embed 
$(\mathbf h_v)_{v\in t}$ randomly in $\mathbb{W}^*_{[0, 1]}$. 
More specifically, for all $v\ino t$, let $\gzH_v \ino \bbZ$ be such that 
 \begin{equation}
  \label{ZBrRW}
\forall v\ino t, \quad | \gzH_v \! -\! \gzH_{\overleftarrow{v}}| = 1 \; \quad \textrm{and} \quad \inf_{\; \, w\in \lgeo  v ,u \rgeo} \! \! \!\!\gzH_w\;  \underset{^{|v|\rightarrow -\infty}}{-\!\!\!-\!\!\!-\!\!\! -\!\!\! \longrightarrow} -\infty . 
\end{equation}
Let $(U_v)_{v\in t}$ be a family of independent r.v.~that are uniformly distributed on $[0, 1]$. With $v\ino t$, we associate a spatial position in $\Wif$ as follows. Since $\gzH$ takes arbitrary negative 
values on the lineage of $v$, for all integers $k \! \leq \! \gzH_v$,
\begin{equation}
\label{sdfljbcv}
\textrm{there is a unique $v(k) \! \in \,  \rgeo \! -\! \infty, v\rgeo$ such that} \quad  \gzH_{v(k)}\! = \!\!\!\!\!\!\!\!  \min_{ \quad w\in \lgeo v(k), v\rgeo} \!\!\!\!\! \!\gzH_w = \, k \, > \gzH_{\overleftarrow{v(k)}} \; .
\end{equation}
Namely, $ \{ \overleftarrow{v(k)}; k\! \leq \! \gzH_v \}$ is the set of vertices in the lineage of $v$ where $\gzH$ reaches a new infimum. In particular, note that $v\! =\!  v(\mathbf{h}_v)$. Then, we set 
\begin{equation}
\label{fgorth}
\forall v\ino t, \quad X_v = (U_{v(k)})_{k\leq \gzH_v} \;.
\end{equation} 
By construction $|X_v|\! = \! \gzH_v$, $v\ino t$. 
Then, $(t,u; (X_v)_{v\in t})$ is a $\Wif$-valued branching random walk satisfying 
\begin{equation*}
 X_v = \left\{  \begin{array}{ll}
\!\! \! X_{\overleftarrow{v} } \ast (U_{v})  &  
\textrm{if $\; |X_v| \! -\! |X_{\overleftarrow{v}}| \! = \! 1$,} \\
\quad  \overleftarrow{X_{\overleftarrow{v} } } \!\!\!\!\!  &  
\textrm{if $\; |X_v| \! - \! |X_{\overleftarrow{v}}| \! = \! -1$.} 
\end{array} \right.
\end{equation*}
We view $(X_v)_{v\in t}$ as a version of free branching random walk conditionally given the relative heights  $(\gzH_v)_{v\in t}$ of the spatial positions. 
The following proposition is a key point to analyse the metric of the range of free branching random walks.
\begin{proposition}
\label{founoiyr}
Let $ \bt\! = \! (t,u) \ino \overline{\bbT}^\bullet$, $(\gzH_v, U_v, X_v)_{v\in t}$ be as above. 
Then, $\{ X_v; v\ino t\}$ is a subtree of $\Wif$ as in Definition \ref{Wsubtree} and if we denote by $d_{\mathtt{gr}}$ it graph-distance, we get  
\begin{equation}
\label{clecle}
\textrm{$\bP$-a.s.~for all $v, w\ino t$,} \qquad 
d_{\mathtt{gr}} 
(X_v,X_w) = |X_v|+|X_w| -2\!\!\!\!  \min_{\; \, v^\prime\in \lgeo v, w\rgeo } \!\!\!\!  |X_{v^\prime}| \; .
\end{equation}
\end{proposition}
\noi
\textbf{Proof.} We fix $v, w\ino t$. By 
(\ref{graphsubt}) we known that $d_{\mathtt{gr}} (X_v, X_w)\! = \! |X_v|+ |X_w| \! -\! 2|X_v \wedge X_w|$. Since $t$ is countable, to prove (\ref{clecle}) it is enough to show that a.s.~$|X_v \wedge X_w|  \! = \! \min_{v^\prime\in \lgeo v, w\rgeo }   |\gzH_{v^\prime}|$. We prove it in two steps.

\textit{Step 1.} Recall from (\ref{sdfljbcv}) the definition of the $v(k)$, $k\! \leq \! \gzH_v$ and define similarly the $w(k)$, $k\! \leq \! \gzH_w$. We then set 
$$ k= \max \big\{ \ell \! \leq \! \gzH_v: v(\ell)\!  \in \, \rgeo \! -\! \infty , v\! \wedge \! w \rgeo \big\}  \quad \textrm{and} \quad  m= \max \big\{ \ell \! \leq \! \gzH_w: w(\ell)\!  \in \, \rgeo \! -\! \infty , v\! \wedge \! w \rgeo \big\} .$$
Let us prove that $k\! = \! \min_{\lgeo v, v\wedge w \rgeo} \gzH$. 
If $k\! < \! \gzH_v $, then $v(k+1) \! \in\, \rgeo v \wedge w , v \rgeo$ and $\overleftarrow{v}(k+1)\ino \lgeo v\wedge w  , v\rgeo$; then $\gzH_{\overleftarrow{v}(k+1)}\! = \!  \gzH_{v(k+1)} \! -\! 1=\!  k \! = \! \min_{\lgeo v, v\wedge w \rgeo}  \gzH$. 
If $k\! =\! \gzH_v$, then, by \eqref{sdfljbcv}, $\gzH_v\! = \! \gzH_{v(k)}\! = \! \min_{\lgeo v, v(k)\rgeo }\gzH$ and since $v(k)\! \in \rgeo -\infty, v\wedge w \rgeo $, we also get $k\! = \! \min_{\lgeo v, v\wedge w \rgeo}  \gzH$. Similarly, we get $m\! = \! \min_{\lgeo w, v\wedge w \rgeo}  \gzH$. By definition (\ref{fgorth}), 
\begin{equation}
\label{suuppin}
X_v= X_{v(k)} \ast \big( U_{v(k+1)}, \ldots, U_{v(\gzH_v)} \big) \quad \textrm{and} \quad X_w= X_{w(m)} \ast \big( U_{w(m+1)}, \ldots, U_{w(\gzH_w)} \big)
\end{equation}
with the observation that $X_v\! = \! X_{v(k)}$ (resp.~$X_w\! = \! X_{w(m)}$) if $k\! = \! \gzH_v$ (resp.~if $m\! = \! \gzH_w$). Without loss of generality, we can assume that $k\! \leq \! m$. Then, $v(k\! -\! i)\! = \! w(k\! -\! i)$, for all $i\ino \bbN$ and 
\begin{equation}
\label{uuppin}
X_{w(m)}= X_{v(k)} \ast \big( U_{w(k+1)}, \ldots, U_{w(m)} \big) \; .
\end{equation}

\textit{Step 2.} We conclude the proof by proving that 
\begin{equation}
\label{factorrr}
X_v \wedge X_w\! = \! X_{v(k)}\; .
\end{equation}
Suppose first that $k\! < \! m $ and that $k\! <\! \gzH_v$. Then $w(k+1) \! \in \, \rgeo \!  -\! \infty, v \wedge w \rgeo $ but $v(k+1) \! \in\,  \rgeo v\wedge w , v \rgeo$. Thus, $w(k+1)\! \neq \! v(k+1)$ and since the $(U_{v^\prime})_{v^\prime \in t}$ are independent with a diffuse law, a.s.~$U_{w(k+1)} \! \neq \! U_{v(k+1)}$. By (\ref{suuppin}) and (\ref{uuppin}), we get (\ref{factorrr}).     
If $k\! < \! m $ and $k\! =\! \gzH_v$, then $X_v\! = \! X_{v(k)}$ and we immediately get (\ref{factorrr}) by (\ref{suuppin}) and (\ref{uuppin}). 

Next, suppose that $k\! = \! m$. Thus $X_{w(m)}\! = \! X_{v(k)}$ by (\ref{uuppin}); if $k\! = \! m\! < \! \gzH_v \wedge \gzH_w$, then $v(k+1) \! \in \, \rgeo  v\wedge w , v\lgeo$ and $w(k+1)\! = \! w(m+1)  \! \in 
\, \rgeo v\wedge w , w \rgeo$; therefore $w(k+1) \! \neq \! v(k+1)$, which implies that a.s.~$U_{w(k+1)} \! \neq  \! U_{v(k+1)}$ and (\ref{factorrr}) consequently. 
If $k\! = \! m$ and $k\! = \! \gzH_v$ (resp.~$m\! = \! \gzH_w$), then $X_v\! = \! X_{v(k)}\! = \! X_{w(m)}$ (resp.~$X_{w}\! = \! X_{w(m)}\! = \! X_{v(k)}$) and (\ref{suuppin}) and (\ref{uuppin}) also entail (\ref{factorrr}). 
This completes the proof of (\ref{factorrr}). \cqfd 

\medskip 
We shall use Proposition \ref{founoiyr} under the following form. 
\begin{corollary}
\label{freeran} Let $t^\prime\ino \bbT$ be a finite rooted ordered tree. Let $(t^\prime , \varnothing ; (Y_{v})_{v\in t^\prime})$ be a $\bbW_{\! [0, 1]}$-valued branching random walk with law $Q^+_{\mathbf{t}^\prime}$ as in Definition \ref{defimodel} $(ii)$. Then, a.s.~for all $v, w\ino t^\prime $, we get $d_{\mathtt{gr}} 
(Y_v,Y_w)\!  =\!  |Y_v|+|Y_w| \! -\! 2  \min_{\; \, v^\prime\in \lgeo v, w\rgeo } |Y_{v^\prime}| $. 
\end{corollary}
\noi
\textbf{Proof.} For all $k\ino \bbZ_-$ we set $\varrho_k\! = \! (a_n)_{n\leq k}$, where $a_n\! = \! 1$ for all $n\ino \bbZ_-$, and we also define  
$t\! = \! \{ \varrho_0   \! \ast \!  u\, ; \, u\ino t^\prime\}\cup \{ \varrho_k; k\ino \bbZ_-\}$. Then, 
$t\ino \overline{\bbT}^\bullet$, $|t|_-\! = \! -\infty$ and $t^\prime\! =\! \theta_{\! \varrho_0} t$. Let 
$(\gzH^\prime_v)_{v\in t^\prime}$ be distributed as an $\bbN$-valued branching random walk whose initial position $\gzH^\prime_\varnothing$ is 
$0$ and whose transition kernel $q(y, dy^\prime)$ is that of the simple symmetric random walk on $\bbN$ reflected at $0$: namely,  $q(0, dy^\prime)\! = \! \delta_1$ and $q(y, dy^\prime)\! = \! \tfrac{_1}{^2} (\delta_{y-1} (dy^\prime) +\delta_{y+1} (dy^\prime))$, for all integers $y^\prime \! \geq \! 1$ (see (\ref{Ybbs})). 
We next set $\gzH_{\varrho_0 \ast v}\! = \! \gzH^\prime_{v}$, for all $v\ino t^\prime$ and $\gzH_{\varrho_k}\! = \! k$, for all 
$k\ino \bbZ_-$. Clearly $\gzH$ satisfies (\ref{ZBrRW}). We assume that $(U_{v})_{v \in t}$ is independent from $(\gzH_{v})_{v \in t}$ and we define $(X_{v})_{v \in t}$ as in (\ref{fgorth}). For all $v\ino t^\prime$, we finally set 
$Y_v\! = \theta_{X_{\varrho_0}} X_{\varrho_0 \ast v}$. Then, it is easy to see that $(t^\prime, \varnothing, (Y_v)_{v\in t^\prime})$ has law $Q^+_{\mathbf{t}^\prime}$. Since $\gzH_{\varrho_0\ast v }\! \geq \! 0$, for all $v\ino t^\prime$, we get $d_{\mathtt{gr}} (X_{\varrho_0\ast v}, X_{\varrho_0\ast w})\! = \! d_{\mathtt{gr}} (
Y_{v}, Y_{ w})$ and $|X_{\varrho_0 \ast v}|\! = \! |Y_v|$, for all $v,w\ino t^\prime$, which implies the desired result by (\ref{clecle}). \cqfd

We next consider subranges of free branching random walks. More precisely, let $ \bt\! = \! (t,u) \ino \overline{\bbT}^\bullet$ be such that $|t|_{\! -}\! = \! -\infty$; for all $v\ino t$, let $\gzH_v\ino \bbZ$ satisfy (\ref{ZBrRW}); 
let $(U_v)_{v\in t}$ be a family of independent r.v.~that are uniformly distributed on $[0, 1]$ and let $(X_v)_{v\in t}$ be derived from $(U_v)_{v\in t}$ as specified in (\ref{sdfljbcv}) and (\ref{fgorth}). 
Let $\bba\subset t$ and $v_0\ino \bba$ be such that 
\begin{equation}
\label{quasisbt}
\forall v\ino \bba\backslash \{ v_0\}, \; \textrm{$\overleftarrow{v} \ino \bba$, $\; v$ is a descendent of $v_0$ (namely $v\! = \! v_0 \! \ast \! (\theta_{v_0} v)$) and $\gzH_v\! \geq \! \gzH_{v_0}$.}
\end{equation}
Observe that for all $v\ino \bba$, {$X_{v_0}$ is a prefix of $X_{v}$}  and it makes sense to define ``the subrange'':
\begin{equation*}
\cR (\bba)= \big\{ \theta_{X_{v_0}} \! X_v \, ; \, v\ino \bba \big\} 
\end{equation*}
that is a subtree of $\Wif$ whose elements have a null depth. 
We define the following.  
\begin{equation}
\label{dpsuebf}
\forall v, w\ino t, \quad d(v,w) \! =\!  \gzH_v + \gzH_w \! -\!  2\min_{\lgeo v,w\rgeo} \gzH \; .
\end{equation}
Note that $\lgeo v, w\rgeo \! \subset \! \bba$; thus, the pseudo-metric $d$ on $\bba\! \times \! \bba$ only depends on $\bba$ and on $(\gzH_v\! -\! \gzH_{v_0} )_{v\in \bba}$. We define the relation $\sim$ on $\bba$ by setting  
$v\sim w$ if and only if  $d(v,w)\! = \! 0$ and we introduce 
\begin{equation}
\label{ttltoutim}
T(\bba)= \bba \, / \! \sim , \quad \mathtt{proj} \! : \! \bba \rightarrow T(\bba),  \textrm{the canonical projection}, \quad r= \mathtt{proj} (v_0) 
\end{equation}
and we keep denoting $d$ the (true) metric induced by $d$ on $T(\bba)$. If $ \bt\! = \! (t,\varrho) \ino \overline{\bbT}^\bullet$ and $(\gzH_v, U_v, X_v)_{v\in t}$ are as in Proposition \ref{founoiyr}, then (\ref{clecle}) implies that 
$(T(\bba) , d)$ is a graph-tree that is isometric to the subtree $\cR (\bba)$. More precisely, let $x\ino T(\bba)$ and let $v, w\ino \bba$ such that $\mathtt{proj}(v)\! = \! \mathtt{proj}(w) \! 
= \! x$; by (\ref{clecle}), we get a.s.~$X_v\! = \! X_w$ and it makes sense to set $Z_x\! = \! \theta_{X_{v_0}}X_{v}$. Then (\ref{clecle}) asserts that  
$Z: T(\bba) \! \rightarrow \cR (\bba)$ is an isometry:
\begin{equation}
\label{dkfv}
\forall x,y\ino T(\bba), \quad d(x,y)= |Z_x|+ |Z_y|\! -\! 2|Z_x \! \wedge \! Z_y|= d_{\mathtt{gr}} (Z_x, Z_y) \; .
\end{equation}
Thus, the graph-metric of the subtree $\cR(\bba)$ only depends on $\bba$ and on $(\gzH_v\! -\! \gzH_{v_0})_{v\in \bba}$. Next, the conditional law of $\cR(\bba)$ given $T(\bba)$ is characterized as follows. 
For all $x\ino T(\bba)\backslash \{ r\}$, let $V_x$ be the unique real number of $[0, 1]$ such that 
\begin{equation}
\label{applSof}
Z_x \! = \! \overleftarrow{Z}_x \! \ast \! (V_x) \; .
\end{equation}
We easily check the following.  
\begin{equation}
\label{reprzcs}
\textrm{Conditionally given $T(\bba)$, the $V_x$ are i.i.d.~$[0, 1]$-uniform r.v.}
\end{equation}
Recall from Definition \ref{tracage} the $\bgg$-contraction map $\Phi_\bgg \! : \! \overline{\bbW}_{[0, 1]}\!  \rightarrow \!  \overline{\bbW}_{\! \bgg}$. 
Then, first note that $|Z_x\wedge Z_y| \! \leq \! |\Phi_\bgg (Z_x) \wedge \Phi_\bgg (Z_y)|$, where common ancestors 
are taken respectively 
in $\bbW_{[0, 1]}$ and in $\bbW_{\! \bgg}$. Moreover, (\ref{reprzcs}) implies that 
 $\bP (|\Phi_\bgg (Z_x) \!  \wedge \! \Phi_\bgg (Z_y)| \! -\!  | Z_x \! \wedge \! Z_y| \! \geq \!  k) \! \leq \! \bgg^{-k}$, for all $k\ino \bbN$. 
This inequality combined with the argument of the proof of Corollary \ref{freeran} implies 
 the following lemma that {will be} used in Theorem~\ref{main2}. 
\begin{lemma}   
\label{contrace} 
 Let $t\ino \bbT$ be a finite rooted ordered tree. Let $(t , \varnothing ; (Y_{v})_{v\in t})$ be a $\bbW_{\! [0, 1]}$-valued branching random walk with law $Q^+_{\mathbf{t}}$ as in Definition \ref{defimodel} $(ii)$. Recall the $\bgg$-contraction map $\Phi_\bgg \! : \! \bbW_{\! [0, 1]}\!  \rightarrow \! \bbW_{\! d}$ from Definition \ref{tracage} and recall that $(t, \varnothing,
 ( \Phi_\bgg (Y_v))_{v\in t})$ has law $Q^{+ \bgg}_{\mathbf{t}}$ as in Definition \ref{defimodel} $(i)$. We denote 
 the graph distance on $\bbW_{\! [0, 1]}$ and $\bbW_{\! \bgg}$ in the same way by $d_{\mathtt{gr}}$. Then, for all $v, w\ino t$, there exists an $\bbN$-valued r.v.~$G_{v,w}$ such that
\begin{equation}
\label{zip}
2G_{v,w}\! = \! d_{\mathtt{gr}} (X_v, X_w) \! -\! d_{\mathtt{gr}} (\Phi_\bgg (X_v), \Phi_\bgg(X_w))  \quad \textrm{and} \quad \bP \big( G_{v,w} \geq k \big) \leq \bgg^{-k}, \;  k\ino \bbN. 
\end{equation} 
\end{lemma}

\subsection{A coupling between $\bbW_{\! \bgg}$- and $\bbW^{*}_{\! \bgg}$-valued branching random walks.}
\label{zdczze}
This section is devoted to the proof of the following proposition. 
\begin{proposition}
\label{coouver} Let $t\ino \bbT$ be a finite rooted ordered tree as in Definition \ref{ort} (namely, $|t|_-\! = \! 0$). We set $\bt\! = \! (t, \varnothing)$. For all $n\ino \{ 0, \ldots , \#t \! -\! 1\}$, we set $t_n\! = \! \{ v\ino t: v \leq_t  v_n \}$ where $v_n$ is the $n$-th smallest vertex of $t$ with respect to the lexicographical order $\leq_t$.     
Then, there exists 
two branching random walks 
$\fTheta_\bgg \! =\! (t,\varnothing; (Y_v)_{v\in t})$ and $\fTheta^+_\bgg \! =\! (t,\varnothing; (Y^+_v)_{v\in t})$ that satisfy the following. 
\begin{compactenum}

\smallskip

\item[$(a)$] $\fTheta^+_\bgg$ is a $\bbW_{\! \bgg}$-valued branching random walk with law 
$Q^{+\bgg}_\bt$ as in Definition \ref{defimodel} $(i)$.

\smallskip

\item[$(b)$] $\mathbf{o}_\bgg\! := \! Y_{\varnothing}$ is a $\bbZ_-$-indexed sequence of independent r.v.~that are uniformly distributed on $\{1, \ldots, \bgg \}$ and conditionally given $\mathbf{o}_\bgg$, 
$\fTheta_\bgg$ is a $\Wbd$-valued branching random walk with law 
$Q^{\bgg}_{\mathbf{o}_\bgg , \bt}$ as in Definition \ref{defimodel} $(iii)$. 

\smallskip

\item[$(c)$] For all $c\ino \bbN^*$, there exists an event $B_c$ such that $\bP (B_c)\! \leq \! 2\bgg^{-c} (\#t)^3$ and such that on $\Omega\backslash B_c$, 
\begin{equation}
\label{vouasx}
\forall n\ino \{ 0, \ldots, \#t \! -\! 1\},  \; \, 
\big| \# \{ Y_v\, ; v\ino t_n\} \! -\# \{ Y^+_v ; v\ino t_n\} \big| \leq \#  \{ v\ino t\! : |Y^+_v| \! \leq \! c+1\} \; .
\end{equation}
\end{compactenum}
\end{proposition}

\paragraph*{Overlap of independent trees.} We first prove a result concerning the overlap of independent trees that are randomly embedded in the {$\bgg$-ary} tree $\bbW_{\! \bgg}$. 
More precisely, let $(\ttT_u, r (u))$, $u\ino S$ be a finite family of rooted graph-trees (not necessarily ordered) 
equipped with their graph-distance $d_{\mathtt{gr}}$. 
To simplify notation, we set $|x|\! = \!  {d_{\mathtt{gr}}} \,(r(u), x)$, for all $x\ino \ttT_u$.  
Let $x,y\ino \ttT_u$; recall that $\lgeo x,y\rgeo$ stands for the shortest path joining $x$ to $y$ 
and that $x\wedge y$ is the most recent common ancestor of $x$ and $y$ in $\ttT_u$ rooted at $r(u)$. Namely, $\lgeo r(u), x\rgeo \cap \lgeo r(u), y\rgeo \! = \! \lgeo r(u), x\wedge y\rgeo$. Let  $x\ino \ttT_u \backslash \{ r(u)\}$; for all $i\ino \{ 1, \ldots ,|x|\}$ we denote by $x_i$ the unique ancestor of $x$ at height $i$. Let $V^u_x$, $x\ino \ttT_u$, $u\ino S$ be independent uniform r.v.~on $\{ 1, \ldots , \bgg\}$. 
Then, we define the random word:
$$Z^u_x \! =\!  \big( V^u_{x_1} , V^u_{x_2}, \ldots ,  V^u_{x_{|x|-1} } , V^u_x \big) $$
We also set $Z^u_{r(u)}\! = \! \varnothing$ and we introduce the following random subsets of the $\bgg$-ary tree $\bbW_{\! \bgg}$:
\begin{equation*} 
\forall u \ino S, \quad \ttR_u := \big\{ Z^u_x ; \, x\ino \ttT_u \big\}\; . 
\end{equation*}
For all $c\ino \bbN^*$, we also set $M_c \!= \! \sum_{u\in S} \# \{ x\ino \ttT_u : |x| \! \leq \! c \}$. 
\begin{lemma} 
\label{overlap} We keep the notation from above. Let $c\ino \bbN^*$. Let $w_u\ino \bbW_{\! \bgg}$, $u\ino S$. Then, there exists an event $A_c$ of probability 
$\bP(A_c)\leq \bgg^{-c} M_c^2$ such that on $\Omega \backslash A_c$, 
\begin{equation}
\label{ooverlap}
0 \leq \sum_{u\in S} \! \# \ttR_u  \; - \# \Big( \bigcup_{u\in S} w_u \! \ast \! ( \ttR_u ) \Big)   \leq M_c \; .
\end{equation}
\end{lemma}
\noi
\textbf{Proof.} Let $u$ and $u^\prime$ be distinct elements of $S$; let $x\ino \ttT_u$ and $y\ino \ttT_{u^\prime}$ 
be such that $|x|\! = \! |y| \! = \! {c}$. We introduce the event 
$A(x,y)\! = \! \big\{ (w_u \! \ast \! Z^u_x ) \!  \wedge \! (w_{u^\prime}\!  \ast \! Z^{u^\prime}_y ) \ino \{ w_u \! \ast \!  Z^u_x , w_{u^\prime} \! \ast \! Z^{u^\prime}_y \} \big\}$. Note that on $A(x,y)$, $w_{u} \! \wedge \! w_{u^\prime} 
\ino \{ w_{u} , w_{u^\prime} \}$; so without loss of generality, we can suppose that $|w_{u} | \! \leq \! |w_{u^\prime}  |$ and we easily check that $\bP (A(x,y) | Z^{u^\prime}_y )\! = \! \bgg^{-c}$ and thus $\bP (A(x,y))\! =\! \bgg^{-c}$. 
We then set for all $u\ino S$, 
$$ \ttT_u  (c)= \{ x\ino \ttT_u: |x|\! =\! c\} \quad \textrm{and} \quad A_c = \bigcup \big\{ A(x,y)\, ; \; x \ino \ttT_u  (c), y \ino \ttT_{u^\prime} (c), u,u^\prime \ino S, \; \textrm{distinct}  \big\} .$$
Thus, $\bP (A_c) \leq \bgg^{-c}\sum_{u,u^\prime \in S, u\neq u^\prime} \# \ttT_{u} (c)  \# \ttT_{u^\prime} (c) \leq \frac{_1}{^2} \bgg^{-c}(\sum_{u\in S} \# \ttT_{u} (c))^2 \leq \bgg^{-c}M_c^2$.  

The first inequality in (\ref{ooverlap}) is true everywhere on $\Omega$.
For the second inequality in (\ref{ooverlap}), we argue deterministically on $\Omega \backslash A_c$. 
Let $u\ino S$. We first set $R_u (c) \! = \! \{ w_u \ast Z^u_x; x\ino \ttT_u (c) \} $. 
Let $u^\prime \ino S$ be distinct from $u$; suppose that 
$W\ino R_u (c)$ and that $ W^\prime \ino R_{u^\prime} (c)$; by definition of $A_c$, $W\wedge W^\prime \! \notin \! \{ W, W^\prime\}$ and $\theta_W  (w_u \ast \ttR_u)  \cap \theta_{W^\prime}   (w_{u^\prime} \ast \ttR_{u^\prime})  \! = \! \emptyset$. Moreover, if $W, W^\prime\ino R_u (c)$ are distinct, 
since $|W|\! = \! |W^\prime|\! = \! c+ |w_u|$, we also get $W\wedge W^\prime \! \notin \! \{ W, W^\prime \}$ 
and thus $\theta_W  (w_u \ast \ttR_u)  \cap \theta_{W^\prime}   (w_{u} \ast \ttR_{u})  \! = \! \emptyset$. 
Consequently, on $\Omega \backslash A_c$, the subsets $\theta_W (w_u \ast \ttR_u) $, $W\ino R_u (c)$, 
$u\ino S$, are pairwise disjoint.

Now observe that $(w_u \ast \ttR_u)  \backslash (\bigcup_{W\in R_u(c)} \theta_W (w_u  \! \ast \ttR_u))  \! = \! 
\{ w_u \! \ast \! Z^u_x ;  |x| \! < \! c\}$ and note that $\# \big( \{ w_u \! \ast \! Z^u_x ;  |x| \! < \! c\} \big) \! \leq \! \# \big( \{ x\ino \ttT_u  \! : |x| \! < \! c  \}\big)$. 
Thus, on $\Omega \backslash A_c$, we get 
\begin{eqnarray*}
 \Big( \sum_{u\in S} \# \ttR_u  \Big)- M_c & =&  \sum_{u\in S}\!\! \big( \# \ttR_u  -\#  \{ x\ino \ttT_u  : |x| \! < \! c  \} \big) \\
 & \leq &  \sum_{u\in S}\!\! \big( \# (w_u \ast \ttR_u)  -\#   \big( \{ w_u \! \ast \! Z^u_x ;  |x| \! < \! c\} \big)   \big) \\
& \leq &  \sum_{u\in S} \# \Big( \!\!\!\!\!\!  \bigcup_{ \quad W\in R_u (c)}  \!\!\!\!\!\!  \theta_W (w_u \! \ast \! \ttR_u)   \Big) \\
&  =&  \#\Big(  \bigcup_{u\in S} \!\!\!\!\!\!  \bigcup_{\quad W\in R_u (c)} \!\!\!\!\!\!   \theta_W (w_u \! \ast \! (\ttR_u) )   \Big) 
 \leq \#  \Big( \bigcup_{u\in S}  w_u \! \ast \! (\ttR_u )\Big), 
 \end{eqnarray*}
that implies the desired result. \cqfd

\paragraph*{Path coupling.} We first state the following elementary coupling.
\begin{lemma}
\label{couplpl} Let $K\! : \! \bbN \! \rightarrow \! \bbN$ be such that $K(0)\! = \! 0$ and $K(2p+1) \! = \! K(2p+2) \! = \! 2p+2$, for all $p\ino \bbN$. Let $(\gzH_n)_{n\in \bbN}$ be a $\bbZ$-valued simple symmetric random walk such that a.s.~$\gzH_0 \! = \! 0$. 
For all $n\ino \bbN$, we set $I_n \! = \! \inf_{0\leq k\leq n} \gzH_k$ and $\gzH^+_n\! = \! \gzH_n + K(-I_n)$. Then, 
$\gzH^+$ is {an} $\bbN$-valued simple symmetric random walk reflected at $0$. 
\end{lemma}
\noi
\textbf{Proof.} For all $n\! \geq \! 1$, set $\xi_n \! = \! \gzH_n \! -\! \gzH_{n-1}$; the r.v.~are i.i.d.~and uniform on $\{ 1, -1\}$. 
Observe that $\gzH^+_n \! = \! \gzH_n\! -\! I_n +\un_{\{ I_n \; \mathrm{odd} \}}$. It is easy to check that $\gzH^+_{n+1}\!- \! \gzH^+_n \!= \!  \xi_{n+1} \un_{\{ \gzH^+_n \geq 1\}}+ \un_{\{ \gzH^+_n = 0 \}}$, which implies the desired result. \cqfd 

\medskip

The next lemma state the branching random walk version of this coupling.

\begin{lemma}
\label{coucoup} Let $t\ino \bbT$ be a rooted ordered tree as in Definition \ref{ort} (namely, $|t|_-\! = \! 0$). Let $(t, \varnothing; (\gzH_u)_{u\in t})$ be a $\bbZ$-valued branching random walk
whose spatial motion is that of a simple symmetric random walk on $\bbZ$ and whose initial position is $\gzH_\varnothing\! = \! 0$.
Recall the function $K$ from Lemma \ref{couplpl} and set 
\begin{equation}
\label{dfvsdkj}
\forall u\in t, \quad I_u= \min_{u^\prime\in \lgeo \varnothing, u \rgeo} \gzH_{u^\prime} \quad \textrm{and} \quad \gzH^+_u= \gzH_u + K(-I_u) \; .
\end{equation}
Then, $(t, \varnothing; (\gzH^+_u)_{u\in t})$ is a branching random walk
whose spatial motion is that of {an} $\bbN$-valued simple symmetric random walk reflected at $0$ and whose initial position is $\gzH^+_\varnothing\! = \! 0$. 
\end{lemma}
\noi
\textbf{Proof.} We denote by $q_+(x,dy)$ the transition kernel of the $\bbN$-valued simple symmetric random walk reflected at $0$: namely, $q_+(0, dy)\! = \! \delta_1 (dy)$ and $q_+(x, dy)\! =\! \frac{_1}{^2} (\delta_{x-1} (dy) + \delta_{x+1} (dy))$, if $x\! \geq \! 1$. For all $n\ino \bbN$, we set $t_{| n}\! = \! \{ u\ino t \! :\!  |u| \! \leq \! n \}$ and we assume the following property 

\smallskip

\begin{compactenum}
\item[$(P_n):$] \textit{the spatial motion of 
the branching random walk $\fTheta_n\! = \! (t_{|n}, \varnothing; (\gzH^+_u)_{u\in t_{| n}})$ is that of a $\bbN$-valued simple symmetric random walk reflected at $0$ with initial position 
$\gzH^+_\varnothing\! = \! 0$}.
\end{compactenum}

\smallskip

We next set $S\! = \! \{ u\ino t \! : \! |u| \! = \! n+1\}$ and $\xi_u \! = \! \gzH_u \! -\! \gzH_{\overleftarrow{u}}$ for all $u\ino S$; the r.v.~$(\xi_u)_{u\in S}$ are independent from $\fTheta_n$ and they are 
also i.i.d.~and uniform on $\{ 1, -1\}$; by definition of $\gzH^+\! $, we also get 
$\gzH^{+}_u\! = \! \gzH^{+}_{^{\overleftarrow{u}}} + \xi_u \un_{^{\!  \{\gzH^{+}_{{\overleftarrow{u}}} \geq 1 \}}}  +\un_{^{\! \{\gzH^{+}_{ \overleftarrow{u}} =0 \}}}$. This entails that 
the r.v.~$\gzH^+_{u}$, $u\ino S$ are conditionally independent given $\fTheta_n$ and that the conditional law of $\gzH^+_{u}$ is $q_+ (\gzH^{+}_{\overleftarrow{u}}, dy)$. This shows that $(P_n)$ implies $(P_{n+1})$, which recursively proves the lemma since $(P_0)$ holds true trivially. \cqfd

\paragraph*{The coupling of {$\fTheta_{\! \bgg}$ and $\fTheta^+_{\! \bgg}$.}}
 Let $t\ino \bbT$ be a finite rooted ordered tree as in Definition \ref{ort} (namely, $|t|_-\! = \! 0$). We shall use the notation $\bt\! = \! (t, \varnothing)$ for the pointed tree. 
We denote by $\varrho\! = \! (\ldots, 1, \ldots, 1)$ the infinite sequence of $1$ indexed by $\bbZ_-$. 
For all $l\ino \bbZ_-$, $\varrho_{| l}$ is the infinite sequence of $1$ indexed by the integers $\! \leq \! l$. We then set 
$$t^*\! = \! \big\{ \varrho_{| l}\, ; \; l\ino \bbZ_- \big\} \cup \big\{ \varrho \! \ast \!  u\, ; \; u\ino t \big\}\; .$$
Then we define $\bt^*\! := \! (t^*\! , \varrho)\ino \overline{\bbT}^\bullet$ and we observe that $|t^*|_-\! = \! - \infty$.  

Let $(t^*\! , \varrho \, ; (\gzH^*_v)_{v\in t^*}\! )$ be a $\bbZ$-valued branching random walk whose transition kernel is that of a simple symmetric random walk on $\bbZ$ and whose ``initial'' position is 
$\gzH^*_{^\varrho}\! = \! 0$. Let $(U^*_v)_{v\in t^*}$ be independent r.v.~that are uniformly distributed on $[0, 1]$ and that are independent from $(\gzH^*_v)_{v\in t^*}$. We define $(X^*_v)_{v\in t^*}$ as in (\ref{sdfljbcv}) and (\ref{fgorth}): namely, $X^*_v\! = \! (U^*_{v(k)})_{k\leq \gzH^*_v}$, where $v(k) \! \in \,  \rgeo \! -\! \infty, v\rgeo$, $\gzH_{v(k)}^*\! = \!  k$ and where $\{ \overleftarrow{v} (k); k\! \leq \! \gzH_v^*\}$ is the set of vertices of $ \,  \rgeo \! -\! \infty, v\rgeo$ where $\gzH^*$ reaches a new strict infimum. We set $\mathbf{o}\! = \! X^*_{^\varrho}$ that is a $\bbZ_-$-indexed sequence of mutually independent r.v.~that are uniformly distributed on $[0, 1]$. 
We easily see that conditionally given $\mathbf{o}$, $\fTheta^*\! = \! (t^*\! , \varrho\, ; (X^*_v)_{v\in t^*}\! )$ has law $Q_{\mathbf{o}, \bt^*}$ as in Definition \ref{defimodel} $(iv)$. 
Then we set the following 
$$ \forall u \ino t , \quad U_u \! = \! U^*_{\varrho \ast u}, \; X_u \! = \! X^*_{\varrho \ast u}, \; \gzH_u \! = \! \gzH^*_{\varrho \ast u}, \; Y_u \! = \! \Phi_\bgg (X_u),  $$
and we also set $ \fTheta\! = \! (t, \varnothing; (X_u)_{u\in t} )$ and $ \fTheta_{ \bgg}\! = \! (t, \varnothing; (Y_u)_{u\in t} )$, where we recall the definition of the $\bgg$-contraction map $\Phi_\bgg$ from Remark \ref{tracage}. 
Clearly, conditionally given $\mathbf{o}$, $ \fTheta$ has law $Q_{\mathbf{o}, \bt}$ as in Definition \ref{defimodel} $(iv)$ and by Remark \ref{tracage}, the conditional law of 
$\fTheta_{ \bgg}$ is $Q^{_\bgg}_{^{\Phi_\bgg(\mathbf{o}) , \bt}}$ as in Definition \ref{defimodel} $(iii)$. Namely, $\fTheta_{ \bgg}$ is distributed as in Proposition \ref{coouver} $(b)$.

Let $(U^\prime_{u})_{u\in t }$ be i.i.d.~r.v.~that are uniformly distributed on $[0, 1]$. We suppose that $(U^\prime_{u})_{u\in t }$ is independent from $(U^*_v)_{v\in t^*}$ (and thus from $(U_u)_{u\in t}$) and  from $(\gzH^*_v)_{v\in t^*}$. 
We next construct $\fTheta_{\bgg}^+$ thanks to $(U^\prime_{u}, U_u, \gzH_u)_{u\in t}$ as follows: 
For all $u\ino t$, we define 
$I_u$ and $\gzH^+_u$ by (\ref{dfvsdkj}) in Lemma \ref{coucoup} and we introduce the following notation. 
\begin{equation}
\label{cuolp}
 \forall p\ino \bbN, \; \,  S_p \! = \! \big\{ u\ino t \! : -p \! = \! I_u \! < \! I_{\overleftarrow{u}} \big\}, \; \,  S \! = \!  \bigcup_{p\in \bbN} S_p   , \;  \,  U^+_u \! = \! \left\{  \begin{array}{ll}
 U^\prime_u & \textrm{if $u\ino \bigcup_{p\in \bbN}S_{2p+1}$, } \\
 U_u & \textrm{otherwise}.
 \end{array} \right.
\end{equation}
The (disjoint) union $S$ is clearly a finite set, since $t$ is a finite tree.
Then, $(U^+_u)_{u\in t}$ are i.i.d.~$[0, 1]$-uniform r.v.~that are independent from $(\gzH^+_u)_{u\in t}$. Recall from Lemma \ref{coucoup} that $(t, \varnothing; (\gzH^+_u)_{u\in t})$ is a branching random walk
whose spatial motion is that of {an} $\bbN$-valued simple symmetric random walk reflected at $0$ and whose initial position is $\gzH^+_\varnothing\! = \! 0$. In particular, 
 $\gzH_u^+\! \geq  \! 0$ for all $u\ino t$. We then define $(X^+_u)_{u\in t}$ and $(Y^+_u)_{u\in t}$ as follows. 
 \begin{compactenum}
 
 \smallskip
 
 \item[$\bullet$] If $\gzH^+_u\! = \! 0$, then we set $X^+_u\! = \! Y^+_u \! = \! \varnothing$.

 \smallskip
 
 \item[$\bullet$] If $\gzH^+_u\! \geq \! 1$, then  for all $i\ino \{ 1, \ldots , \gzH^{+}_{u} \}$, we denote by $u(i)$ is the unique $v \ino \lgeo \varnothing , u\rgeo$ such that 
$\gzH^{+}_{{\overleftarrow{v}}}\! < \! i \! = \! \gzH^{+}_{v}\! = \! \min_{\lgeo v, u \rgeo } \gzH^+ $ and we set 
\begin{equation}
\label{cuolpp}X^+_u \! = \! (U^{_+}_{{\! u(1)}}, \ldots , U^{_+}_{{\!\!  u(\gzH^+_u)}}) \quad \textrm{and} \quad Y_u^+ = \Phi_\bgg (X^+_u) \; .
\end{equation}
\end{compactenum}
We set $\fTheta^+\! = \! (t, \varnothing; (X^+_u)_{u\in t})$ and $\fTheta^+_{\bgg}\! = \! (t, \varnothing; (Y^+_u)_{u\in t})$. 
Then, observe that $\gzH^+_u\! = \! |X^+_u|$ and that $X^+_u\! = \! X^+_{\overleftarrow{u}} \! \ast \! (U^+_u)$ if $\gzH^+_u\! = \! \gzH^+_{\overleftarrow{u}}+1$ and note that 
$X^+_u$ is the parent of $X^{_+}_{^{\overleftarrow{u}}} $ if 
$\gzH^{_+}_{^u}\! = \! \gzH^{_+}_{^{\overleftarrow{u}}} \! -\! 1$. Thus, it proves that $\fTheta^+$ as law $Q^+_{\bt}$ as in Definition \ref{defimodel} $(ii)$ and by Remark \ref{tracage},  $\fTheta^+_{\bgg}$ is distributed as in Proposition \ref{coouver} $(a)$ (namely its law is $Q^{+\bgg}_{\bt}$ as in Definition \ref{defimodel} $(i)$).

 \paragraph*{Proof of Proposition \ref{coouver}.} 
We keep the previous notations.  
We fix $n\ino \{ 0, \ldots , \#t \! -\! 1\}$. Recall that  
 $t_n\! = \! \{ v\ino t: v \leq_t  v_n \}$ where $v_n$ is the $n$-th smallest 
 vertex of $t$ in the lexicographical order $\leq_t$ on $t$. Note that $t_n \ino \bbT$, namely it is a rooted ordered tree with null depth. Recall from (\ref{cuolp}) the definition of $S$. For all $u\ino S$, we set 
 $$ \bba_u = \big\{ v\ino t_n : v\! := \! u \ast \theta_u v  \;  \textrm{ and} \; \min_{w\in \lgeo u, v\rgeo}\gzH_w \! =\! \gzH_u \big\} \subseteq t
  \quad \textrm{and} \quad \cR (\bba_u) \! = \!  \big\{ \theta_{X_u}X_{v} ; v\ino \bba_u \big\} \subseteq\Wbf \,.$$
Note that $\bba_u$ depends on $n$ and that if $u\! \notin \! t_n$, then $\bba_u$ and $\cR(\bba_u)$ are empty. 
If $u\ino t_n$, then the previous definitions make sense because $\bba_u$ satisfies (\ref{quasisbt}) (with $u=v_0$) and $X_u$ is necessarily 
a prefix of $X_v$ since $v$ is a descendent of $u$ and $I_v\! = \! I_u$. 
Observe that $\cR (\bba_u)$ is a subtree of $\Wbf$ as in Definition \ref{Wsubtree}.

  Recall from (\ref{dpsuebf}) the notation $d(v,w) \! =\!  \gzH_v + \gzH_w \! -\!  2\min_{\lgeo v,w\rgeo} \gzH$, for all $v, w\ino t$,
that is a pseudo-metric on $t$. We denote by $v \! \sim \! w$ the equivalence relation relation $d(v,w)\! = \! 0$ and as in (\ref{ttltoutim}), we set for all $u\ino S$, 
 \begin{equation*}
T(\bba_u)= \bba_u / \!\! \sim , \quad \mathtt{proj} \! : \! \bba_u \rightarrow T(\bba_u),  \textrm{the canonical projection}, \quad r(u)= \mathtt{proj}(u) \; .
\end{equation*}
Note that $T(\bba_u)$ only depends on $(\gzH_{v}\! -\! I_u)_{v\in \bba_u}$ and that
$\cR(\bba_u)$ only depends on $(\gzH_{v}\! -\! I_u ; U_v)_{v\in \bba_u}$.
Let $x\ino T(\bba_u)$; 
by (\ref{clecle}) in Lemma \ref{founoiyr}, for all $v, w\ino \bba_u$ such that $\mathtt{proj}  (v)\! = \! \mathtt{proj} (w) \! 
= \! x$, we get $X_v\! = \! X_w$ and it makes sense to set $Z^u_x\! = \! \theta_{X_{u}}X_{v}$; then, 
$Z^u \! : \! T(\bba_u) \! \rightarrow \! \cR (\bba_u)$ is an isometry given in (\ref{dkfv}). 
As in (\ref{applSof}), for all $x\ino T(\bba_u)\backslash \{ r(u) \} $ we denote by $V^u_x$ the unique real number of $[0, 1]$ such that $Z^u_x \! = \! \overleftarrow{Z}^u_x \ast (V^u_x)$. 

\smallskip

Then by   
(\ref{reprzcs}), for all 
$u\ino S$, conditionally given $\gzH$, the $V^u_x$, $x\ino T(\bba_u)$ are i.i.d.~$[0, 1]$-uniform r.v.~Since the subsets $(\bba_u)_{u\in S}$ 
are pairwise disjoint and since $\cR(\bba_u)$ only depends on {$(\gzH_{v}-I_u ; U_v)_{v\in \bba_u}$}, the $(\cR(\bba_u))_{u\in S}$ are conditionally independent given $\gzH$. Thus, 
conditionally given $\gzH$, the $V^u_x$, $x\ino T(\bba_u)\backslash \{ r(u) \}$, $u\ino S$, are i.i.d.~$[0, 1]$-uniform r.v.

    Next observe that 
\begin{equation} 
\label{ralibb}
\{ X_v; v\ino t_n \}= \bigcup_{u\in S}  X_u \! \ast \! (\cR( \bba_u)) {\green \,.}
\end{equation}
Therefore, using \eqref{cuolp} and noting that $(X_u)_{u\in S}$ only depends on $\gzH$ and on $(U^*_{\varrho_{| l}})_{l\in \bbZ_-}$ that is independent from $(U_u)_{u\in t}$, 
we conclude that, conditionally given $\gzH$ and $(X_u)_{u\in S}$, the $V^u_x$, $x\ino T(\bba_u)\backslash \{ r(u) \}$, $u\ino S$, are i.i.d.~$[0, 1]$-uniform r.v. 
From the coupling defined in (\ref{cuolp}) and (\ref{cuolpp}) we also derive easily the following. 
\begin{equation} 
\label{rarefll}
\{ X^+_v; v\ino t_n \} \! = \!\! \bigcup_{u\in S}  X^\prime_u \! \ast \! (\cR( \bba_u)), \; \textrm{where 
for all $u$ in $S$,} \; X^\prime_u:=  \left\{  \begin{array}{ll}
 U^\prime_u & \textrm{if $u\ino \bigcup_{p\in \bbN}S_{2p+1}$, } \\
\varnothing & \textrm{otherwise}.
 \end{array} \right.
\end{equation}
Since the $(U^\prime_u)_{u\in t}$ are independent from $\gzH$ and from $(U_u)_{u\in t}$, 
conditionally given $\gzH$ and $(U^\prime_u)_{u\in S}$, the r.v.~$V^u_x$, $x\ino T(\bba_u)\backslash \{ r(u) \}$, $u\ino S$, are i.i.d.~$[0, 1]$-uniform. Denote by $\ccG$ the sigma-field generated by $\gzH$ and by $(X_u, U^\prime_u)_{u\in S}$. Therefore, we have proved the following. 
\begin{equation} 
\label{fournigl}
\textrm{Conditionally given $\ccG$, the $V^u_x$, $x\ino T(\bba_u)\backslash \{ r(u) \}$, $u\ino S$, are i.i.d.~$[0, 1]$-uniform r.v. }
\end{equation}
Then we set 
$$ \forall u \ino S, \quad \mathtt{R}_u \! = \!  \Phi_\bgg (\cR (\bba_u)) \quad \textrm{and} \quad 
\forall c \ino \bbN, \quad M_c(n) = \sum_{u\in S} \# \{ x\ino T(\bba_u ) : d(r(u), x) \! \leq \! c \} \; .$$
Note that $M_c (n) \!  \leq \! \# t_n \! \leq \! \# t$. 
By (\ref{fournigl}), conditionally given $\ccG$ we can apply Lemma \ref{overlap} to $\mathtt{T}_u \! := \! T(\bba_u)$ and $w_u \!: = \! \Phi_\bgg (X_u)$, $u\ino S$, to get an event $A_c(n)$ such that  $\bP (A_c(n)) \leq \bgg^{-c}( \#t)^2$ and such that on $\Omega \backslash A_c(n)$, 
$$ 0 \leq\sum_{u\in S} \# \mathtt{R}_u \; - \# \big\{  Y_v; \, v\ino t_n  \big\}  \overset{\textrm{by (\ref{ralibb})}}{=}  
\sum_{u\in S} \# \mathtt{R}_u \; - \# \Big(  \bigcup_{u\in S} \Phi_\bgg(X_u) \! \ast \! (\ttR_u)  \Big) \leq M_c(n) \; .$$ 
Similarly, conditionally given $\ccG$, we apply Lemma \ref{overlap} to $\mathtt{T}_u \! := \! T(\bba_u)$ and to $w_u \! := \! \Phi_\bgg (X^\prime_u)$, $u\ino S$, to get an event $A^{\prime}_c (n)$ such that  $\bP (A^{\prime }_c(n)) \leq \bgg^{-c}( \#t)^2$ and such that on $\Omega \backslash A^\prime_c (n)$, 
$$ 0 \leq\sum_{u\in S} \# \mathtt{R}_u \; - \# \big\{  Y^+_v; \, v\ino t _n \big\}  \overset{\textrm{by (\ref{rarefll})}}{=}  
\sum_{u\in S} \# \mathtt{R}_u \; -\# \Big(  \bigcup_{u\in S} \Phi_\bgg(X^\prime_u) \! \ast \! (\ttR_u)  \Big) 
\leq M_c (n)\; .$$ 
Then, we set $B_c(n) \! = \!  A_c (n)\cup A^\prime_c(n)$; thus, $\bP (B_c(n))\! \leq \! 2 \bgg^{-c} (\# t)^2$ and on $\Omega \backslash  B_c (n)$, we get $ | \# \{ Y_v ;  v\ino t_n\} -\# \{ Y^+_v ; v\ino t_n\}| \!   \leq \! M_c  (n)$. 
Now observe that if $x\ino T(\bba_u)$, then $\gzH^+_x \! = \! |X^+_x| \! = \! |X^\prime_u|+ |Z^u_x|\! \leq  \! 1+ d(r(u), x)$. 
This implies $M_c (n) \! \leq  \! \# \{ v\ino t ; |Y^+_v| \! \leq \! c +1 \} $. Thus, on  $\Omega\backslash B_c (n)$, we get 
$$   \big| \# \{ Y_v\, ; v\ino t_n\} -\# \{ Y^+_v ; v\ino t_n\} \big| \leq \# \{ v\ino t ; |Y^+_v| \! \leq \! c +1 \}  \; .$$
We completes the proof of Proposition \ref{coouver} by taking $B_c\! = \! \bigcup_{0\leq n < \# t} B_c (n)$. \cqfd

\subsection{Estimates.}
\label{brwestsec}
The goal of this section is to establish Proposition \ref{ntgngbck} below.  To that end, we first state preliminary estimates. 
Recall that $\bbW^*_{\! \bgg} \! = \!  \{x\ino \overline{\bbW}_\bgg \! : \! |x|_{-} \! = \! -\infty \}$. We denote by $(Y_n)_{n\in \bbN}$ the canonical process on the space $(\Wid)^{\bbN}$ equipped with product topology (that is Polish) and with the corresponding Borel sigma-field. In this section, let us denote by $P_{ \! \rmo}$  (instead of  $Q^{_{\bgg}}_{^{\rmo, \mathbb N}}$) the law of a Markov chain on $\Wid$ with transition kernel $p_\bgg$ as defined in (\ref{brwTinfd}) and whose initial position is $\mathrm{o}$. The following result only contains {some} standard estimates. 
\begin{lemma}
\label{dlczidec} Let $\rmo \ino \Wid$ be such that $|\rmo| \! =\!  0$.
For all $k, p\ino \bbN$, we set  
\begin{equation}
\label{zmvkb}
 \rmo(p)\! = \! \rmo_{|]-\infty, -p]}, \; \,\mathtt{Sp}\! = \! \{ \rmo (p); p\ino \bbN\}, \; \, \mathbf{r}_0\! = \! 0\; \, \textrm{and} \; \, \mathbf{r}_{k+1}\! = \! \inf \{ n \! >\! \mathbf{r}_k : Y_n \ino \mathtt{Sp} \}.
\end{equation}
The following holds true. 
\begin{compactenum}

\smallskip

\item[$(i)$]  $P_{\! \rmo}$-almost surely for all $k\ino \bbN$, $\mathbf{r}_k \! < \! \infty$ and there exists $Z_k\ino \bbN$ such that $\rmo(Z_k)\! = \! Y_{\mathbf{r}_k}$. Moreover, $(Z_k)_{k\in \bbN}$ is an $\bbN$-valued birth-and-death Markov chain whose transition probabilities $(\rho (p,q))_{p,q\in \bbN}$ 
are given as follows: for all $p\in \bbN^*$,  
\begin{equation}
\label{fgfhfj}
\rho(p, p+1) \! = \! \tfrac{1}{2 } , \quad   \rho(p, p\! -\! 1)\! = \! \tfrac{1}{ 2\bgg}, \quad  \rho (p, p)\! =\! 
\tfrac{\bgg - 1
}{2\bgg} \; ,
\end{equation}
and $\rho(0, 0)\!= \! \rho (0, 1) \! = \! 1/2$. Then, $Z$ is transient which implies that almost surely $|\rmo \!  \wedge \! Y_n| \! \rightarrow \! -\infty$ and that $(Y_n)_{n\in \bbN}$ is transient under $P_{\! \rmo}$.

\smallskip

\item[$(ii)$] For all $x,y\ino \Wid$, we set $G_\bgg (x,y)\! = \! \sum_{n\in \bbN}P_x\big(Y_n\! = \! y \big)$. Then, for all $p\ino \bbN$, 
\begin{equation}
\label{greenoj}
G_\bgg (x,y)= \tfrac{2\bgg }{\bgg - 1} \bgg^{|x\wedge y|-|y|} \; .
\end{equation}

\item[$(iii)$] For all $y\ino \bbW^*_{\! \bgg}$, we set $\ttH_y\! = \! \inf\{n\ino \bbN\! : \! Y_n \! = \! y\}$, with the convention that $\inf \emptyset \! = \! \infty$. 
For all $s\ino [0, \infty)$, we set $g (s) \!   =  \! -\!  \log \big(1\! -\! \sqrt{1\! -\! e^{-2s}} ) $. Then, for all $x\ino \Wid$, $E_x [ e^{-s(1+\ttH_{\overleftarrow{x}}) }] \! = \! \exp (  -g (s))$. 
Moreover, there exists a constant $C\ino (1, \infty)$ such that  
\begin{equation}
\label{gournihusj}
\textrm{$P_\rmo$-a.s.~for all $k\ino \bbN$ and for all $s\ino [0, 1]$,} \quad E_{\rmo} \big[ e^{-s (\mathbf{r}_{k+1} -\mathbf{r}_k ) } \big| Z \big] \geq  e^{-C\sqrt{s}} \; .
\end{equation}
\end{compactenum} 
\end{lemma}
\noi
\textbf{Proof.} {Since under $P_{\! \rmo}$, $(|Y_n|)_{n\in \bbN}$ is distributed as a simple symmetric random walk on $\bbZ$, we easily see that $P_{\! \rmo}$-a.s.~for all $k\ino \bbN$, $\mathbf{r}_k \! < \! \infty$.  
The strong Markov property at the stopping times $\mathbf{r}_k$ implies that $(Z_k)_{k\in \bbN}$ is an $\bbN$-valued birth-and-death Markov chain whose 
transition probabilities $\rho (p,q)$ are given by (\ref{fgfhfj}) which easily implies that 
$Z_k\rightarrow \! \infty$ which entails $(i)$. 

Let us prove $(ii)$. 
For all $x \ino \Wid$, we introduce the stopping times 
$\ttH_x\! = \! \inf \{ n\ino \bbN\! :\!  Y_n \! = \! x \}$ and 
$\ttH^\circ_x \! = \! \inf \{ n\ino \bbN^*\! : \! Y_n \! = \! x\}$. First observe that for all $x, y\ino \Wid$, $P_{\! x} (\ttH_y \! < \! \infty)\! = \! P_{\! x\wedge y} (\ttH_y \! < \! \infty)$. By adapting the argument of $(i)$, the height of the process $Y$ restricted to $\rgeo -\infty, y \rgeo$ is a birth-and-death process with transition $\rho (\cdot, \cdot)$ and the Gambler's ruin estimate implies that $P_{\! x\wedge y} (\ttH_y \! < \! \infty ) \! = \! \bgg^{|x\wedge y| -|y|}$. Similarly, we get $P_{\! y} (\ttH^\circ_{y}\! = \! \infty )\! =\! (\bgg\! -\! 1)/2\bgg$.  
The Markov property at resp.~$\mathtt{H}_y$ and  $\ttH^\circ_y$ implies resp.~that 
$G_\bgg (x, y) \! =\!   P_{\! x}  (\ttH_{y} \! < \! \infty)G_\bgg (y,y) $ and that  
$G_\bgg (y,y)\! = \! 1+  P_{\! y}  (\ttH^\circ_{y} \! < \! \infty) \, G_\bgg (y,y) $, which entails (\ref{greenoj}).

Let us prove $(iii)$. Note that under $P_{\! x}$, $(|Y_n|)_{n\in \bbN}$ is distributed as a simple symmetric 
random walk on $\bbZ$ with initial position $|x|$. Then $\ttH_{\overleftarrow{x}}$ is the first time $(|Y_n|)_{n\in \bbN}$ reaches $|x| \! -\! 1$; therefore, the law of $\ttH_{\overleftarrow{x}}$ does not depend on $x$ and we get 
$E_x [ \exp (-s(\ttH_{\overleftarrow{x}}+1) )] \! = \! \exp (  -g (s))$ by well-known arguments. 
Next, let us work conditionally given $Z$; we set $x \! =\!  Y_{\mathbf{r}_k +1}$; on the event $\{Z_k \! = \! Z_{k+1}\}$, $\overleftarrow{x}\! = \! \rmo (Z_k)$ and $\mathbf{r}_{k+1} \! -\! \mathbf{r}_k \!- \! 1$ is the time at which 
the shifted random walk $(Y_{\mathbf{r}_k +1 +n})_{n\in \bbN}$ returns for the first time to $\overleftarrow{x}$; 
thus, by Markov at  $\mathbf{r}_k$ and by the previous argument, $P_{\rmo}$-a.s.~on the event $\{Z_k\! =\!  Z_{k+1} \}$, $E_{\rmo} [ \exp (\!-s (\mathbf{r}_{k+1} \! -\! \mathbf{r}_k ) )| Z ]\! = \! \exp (\! -g(s))$. If $Z_k\! \neq \! Z_{k+1}$, then $\mathbf{r}_{k+1} \! -\! \mathbf{r}_k \! = \! 1$. Consequently, 
$E_{\rmo} [ \exp (\! -s (\mathbf{r}_{k+1} \! -\! \mathbf{r}_k ) ) | Z ] \! =\!  \exp 
( \! -g(s) \un_{\{Z_k= Z_{k+1} \} }\!   - \! s \un_{\{Z_k\neq Z_{k+1} \}})$,
which implies (\ref{gournihusj}) because there exists $C\ino (1, \infty)$ such that $g(s) \! \leq C \sqrt{s}$ for all $s\ino [0, 1]$ and since $s\! \leq \! \sqrt{s}$ for all $s\ino [0, 1]$. }\cqfd

\begin{proposition}
\label{quantrans} 
Let $\rmo \ino \Wid$ such that $|\rmo|\! = \! 0$. Recall that $P_{\! \rmo}$ stands for the canonical law of $\Wid$-valued Markov chains with transition kernel $p_\bgg$ defined by (\ref{brwTinfd}) and with initial position $\rmo$. 
Let $f\! : \! [0, \infty) \! \rightarrow \!  [0, 1]$ be such that $\sum_{p\in \bbN} \sqrt{f(p)} \! < \! \infty$.  
\begin{equation*}
 \textrm{$P_{\! \rmo}$-a.s.}\quad  \sum_{n\in \bbN} f \big( \! - \! |\rmo \! \wedge \! Y_n| \big) \, < \infty \;  . 
\end{equation*}
\end{proposition}
\noi
\textbf{Proof.} Recall from (\ref{zmvkb}) in Lemma \ref{dlczidec} notations $\rmo (p)$, 
$\mathtt{Sp}$ and $\mathbf{r}_k$ and recall from Lemma \ref{dlczidec} $(i)$ the definition of the $\bbN$-valued birth-and-death process $(Z_k)_{k\in \bbN}$. 
Let $\ell \ino \bbN^*$. First observe that  
$$ W_\ell\! := \! \!\!  \sum_{0\leq n<\mathbf{r}_\ell} \!\! f\big( \!\!  -\!  |\rmo  \wedge Y_n| \big)= \sum_{0\leq k< \ell} \!\! \!\! \!\! \!\! \!\! \sum_{\quad \quad  \mathbf{r}_k \leq n< \mathbf{r}_{k+1}} \!\! \!\! \!\! \!\! \!\! \!\!  f( Z_k) \; = 
\sum_{0\leq k<\ell} (\mathbf{r}_{k +1}\! -\! \mathbf{r}_{k}) f(Z_k) . $$
By the Markov property, the sequence of random times 
$\mathbf{r}_{k +1}\! -\! \mathbf{r}_{k}$ are conditionally independent given $Z$. Thus by (\ref{gournihusj}) in 
Lemma \ref{dlczidec} $(iii)$, there is $C\ino (1, \infty)$ such that for all $s\ino [0, 1]$, 
\begin{eqnarray}
\label{hjhh}
E_{\! \rmo}  \big[e^{-sW_\ell} \big|Z  \big] \!\!\! \!  & = &\!\! \!\!   
\prod_{0\leq k< \ell} \bE \big[ e^{-s (\mathbf{r}_{k +1} - \mathbf{r}_{k}) f(Z_k)  } \big| Z \big] \geq \exp \Big( \!\! -\! C \! \! \! \sum_{0\leq k < \ell} \! \!\! \sqrt{sf(Z_k)}   \Big) \nonumber \\ 
 \!\!\! \!  & \geq  &\!\! \!\!    \exp \Big( \!\! -\! C  \sqrt{s}\sum_{p\in \bbN} N_p \sqrt{\! f(p)}   \Big), \; \, \textrm{where $N_p\! = \! \# \big\{k\ino \bbN  \!  : \!  Z_k\! = \! p \big\}$.}
\end{eqnarray}
Next observe that $E_{\rmo} [N_p] \! =\! G_\bgg (\rmo, \rmo (p))\! = \!  2\bgg/(\bgg\! -\! 1)$ by (\ref{greenoj}) in Lemma \ref{dlczidec} $(ii)$. 
Set $W\! = \! \sum_{n\in \bbN} f \big(\!\! -\! |\rmo \! \wedge \! Y_n| \big) $; by letting $\ell\! \rightarrow \! \infty$ in (\ref{hjhh}) and by Jensen's inequality, we get 
\begin{equation*}
E_{\rmo}  \big[e^{-sW} \big] \! \geq \!  \exp \Big( \!\! -\! C \sqrt{s} \sum_{p\in \bbN} \bE [N_p]   \sqrt{\! f(p)} \Big) = \exp \Big( \!\! -\! C^\prime \sqrt{s} \sum_{p\in \bbN}\sqrt{\! f(p)} \Big) 
\; ,
\end{equation*}
where $C^\prime\! = \! 2\bgg C/(\bgg\! -\! 1)$. Since $\sum_{p\in \bbN}\sqrt{\! f(p)} 
 \! < \! \infty$, we get $\lim_{s\rightarrow 0+}\bE \big[e^{-sW}  \big] \! = \! 1$ which entails that $P_{\! \rmo}$-a.s.~$W\! < \!\infty$. This completes the proof of the proposition.  \cqfd 

\begin{lemma}
\label{momtun} Let $\mu$ be a probability measure on $\bbN$ such that $\sum_{k\in \bbN} k\mu (k) \! = \! 1$. 
Let $\ftau\! = \! (\tau, \varnothing)$ be a random pointed tree such that a.s.~$|\tau|_{\! -}\! = \! 0$, $k_\varnothing (\tau)\! = \! 1$ and such that $\theta_{(1)}\tau$ is a GW($\mu$)-tree. Let $x\ino \Wid$ and let 
$\fTheta\! = \! (\tau, \varnothing \, ; (X_{\! v})_{v\in \tau})$ be a random $\Wid$-valued branching random walk that has law $Q^{_\bgg}_{^{x, \ftau}}$ conditionally given $\ftau$, as defined in Definition \ref{defimodel} $(iii)$. 
For all $y \ino \Wid$ we set
\begin{equation}
\label{defixixi}
\xi (x, y) = \bP \big( y \notin \big\{  X_{\! v}\, ; \, v\ino \tau \backslash \{ \varnothing \}  \big\} \big)  \; , 
\end{equation}
which turns out to be strictly positive. Then, we get $  1\! -\! \xi (x, y)\!  \leq\!  G_\bgg (x,y) \!  = \! \tfrac{2\bgg }{\bgg - 1} \bgg^{\, |x\wedge y|-|y|}$.
\end{lemma}
\noi
\textbf{Proof.} By a simple union bound we first get the following. 
\begin{equation*}
  \bP \big( y \ino \big\{  X_{\! v}\, ;  v\ino \tau\backslash \{ \varnothing \}\!  \big\} \;  \big| \, \ftau  \big) 
   \leq  \bE \Big[ \!\!\!  \!\!  \sum_{\quad v\in \tau \backslash \! \{ \varnothing \! \}} \!\! \!\! \!\! \un_{\{ X_v=y \}} \; \Big| \, \ftau  \Big]   \leq   \!\!\!  \!\!  \sum_{\quad v\in \tau \backslash \! \{ \varnothing \! \}}  \!\!\!  \!\!  Q^{\bgg}_{x, \ftau} \big( X_v=y \big) . 
\end{equation*}
By definition of branching random walks, $Q^{_\bgg}_{^{x, \ftau}}  ( X_v \!= \! y)\! = \! P_{\!x} (Y_{|v|}\! = \! y)$, where $P_{\! x}$ stands for the canonical law of the random walk that starts at $x$ in $\Wid$ and whose transition kernel is $p_\bgg$ as defined in (\ref{brwTinfd}) in Definition \ref{defimodel} $(iii)$.  
Thus, $ \bP \big( y \ino \big\{  X_{\! v}\, ;  v\ino \tau\backslash \{ \varnothing \}\!   \big\}  \big)  \leq \sum_{n\geq 1} \bE [N_n] \, P_{\! x} (Y_{n}\! = \! y) $,
where $N_n\! = \! \# \{ v\ino \tau\! : \! |v|\! = \! n\}$. Since the offspring distribution $\mu$ is critical, $\bE [N_n] \! =\!  1$.
Therefore, recalling Lemma \eqref{dlczidec} {\it (ii)}, we have that the right member is smaller or equal to $G_\bgg (x, y)$, which concludes the proof. \cqfd 

\begin{lemma}
\label{abldnez} 
Let $\mu$ be a probability measure on $\bbN$ such that $\sum_{k\in \bbN} k\mu (k) \! = \! 1$. 
Let us suppose that there exists $\beta \ino (0, \infty)$ such that 
$\sum_{k\in \bbN^*} \mu (k) \,  k (\log k)^{1+ \beta}\! <\! \infty$. 
For all $r\ino [0, 1]$, we set 
$\psi (r)\! = \! \sum_{k\in \bbN} r^k  \overline{\mu} (k+1)$, where we recall that for all $k\ino \bbN \backslash \{ 0\}$, 
$\overline{\mu} (k)$ stands for $\sum_{\ell \geq k} \mu(\ell) $.    
Then, there exists $C_{\beta, \mu}\ino (0, \infty)$ such that 
\begin{equation}
\label{flablel}
\forall r\ino (0, 1/2) , \qquad    0\leq 1-\psi (1-r) \leq \frac{C_{\beta, \mu}}{(\log 1/r)^{1+ \beta} } \; . 
\end{equation}
\end{lemma}
\noi
\textbf{Proof.} Let $Z$ be {an} $\bbN\backslash \{ 0\}$-valued random variable distributed according to $\overline{\mu}$. Observe that $C\! :=\! \bE [(\log Z)^{1+\beta} ]\!  = \!\sum_{l \geq k\geq 1 } \mu(l) (\log k)^{1+\beta} \!  \leq \! \sum_{k\geq 1} \mu (k)   k (\log k)^{1+ \beta}\!< \! \infty$. 
Then, for all $t \! >\! 0$ and for all $r\ino (0, 1/2)$, we get 
\begin{eqnarray*} 
1\! -\! \psi (1\! -\! r) \!\!\! & = & \!\!\!  \bE \big[ \big( 1\! -\! (1\! -\! r)^{Z-1}\big) \un_{\{ Z\leq t+1 \}} \big] + 
\bE \big[\big( 1\! -\! (1\! -\! r)^{Z-1}\big) \un_{\{ Z> t+1 \} }  \big] \\
\!\!\!  &\leq & \!\!\!  1\! -\! (1\! -\! r)^t+ \bP (Z\! >\! t+1) \;  \leq \;  2rt+ \frac{\bE [(\log Z)^{1+ \beta}]}{(\log (1+t))^{1+ \beta}}= 2rt + C(\log (1+t))^{-1- \beta}. 
\end{eqnarray*}
We choose $1+t\! = \!  r^{-1} (\log 1/r)^{-1-\beta}$ and we easily get (\ref{flablel}). \cqfd

\begin{proposition}
\label{ntgngbck} Let $\mu$ be a probability measure on $\bbN$ such that 
\begin{equation}
\label{hypomunu}
\sum_{k\in \bbN}   k \mu (k)  = 1   \quad \textrm{and} \quad \exists \beta \ino (1, \infty) \; \textrm{such that} \; \sum_{k\in \bbN} k (\log k)^{1+ \beta}  \mu(k)  <\!  \infty.  
\end{equation}

\noi
Let $\rmo \ino \Wid$ be such that $|\rmo|\! = \! 0$ and let $\fTheta\!\! := \! (\tau^*\! , \, \varrho\, ; \, (X_v)_{v\in \tau^*})$ be a $\Wbd$-valued branching random walk whose distribution is the following: $\ftau^*\! = \! (\tau^*\! , \boo)$ is an IPGW$^+(\mu)$-tree as in Definition \ref{GWptdef} $(ii)$ and conditionally given $\ftau^*$, $\fTheta$ has law $Q^{_\bgg}_{^{\rmo, \ftau^*}}$ as in Definition \ref{defimodel} $(iii)$. 
Then, with positive probability $(X_v)_{v\in \tau}$ visits its initial position only once. Namely, 
\begin{equation}
\label{lpomlpo}
\kappa_{\mu, \bgg} := \bP \big( \rmo \!  \notin \! \big\{ X_v\, ; v\ino \tau^* \backslash  \{\varrho \} \big\} \, \big) >0 \; .
\end{equation}
\end{proposition}
\noi
\textbf{Proof.} Let us recall the following notation
$$ \forall p\ino \bbN, \; \, \varrho (p) = \varrho_{| \, ]-\infty , -p]}, \; \, Y_p\! = \! X_{\varrho (p)}, 
 \quad \mathtt{Sp} \! = \! \{  \varrho (p) \, ; p\ino \bbN\} \; \, \textrm{and} \; \,
\partial \mathtt{Sp}\! = \! \big\{ u\ino \tau^* \backslash \mathtt{Sp}: \overleftarrow{u}\ino \mathtt{Sp} \big\} .$$  
By Definition \ref{GWptdef} $(ii)$ of the right part of an infinite pointed GW($\mu$)-tree, 
the subtrees $(\theta_u \tau^*)_{u\in \partial \mathtt{Sp}}$ are i.i.d.~GW($\mu$)-trees and the r.v.~$(k_{\varrho (p)} (\tau^*))_{p\in \bbN}$ are independent, $k_\boo (\tau^*)$ has law $\mu$ and for all $p\! \geq \! 1$, $k_{\boo (p)} (\tau^*) $ has law $\overline{\mu} $ that is defined by $\overline{\mu} (k)\! = \! \sum_{l\geq k} \mu (l)$, $k\ino \bbN^*$ and $\overline{\mu} (0)\! = \! 0$. 
For all $u\ino \partial \mathtt{Sp}$, we define $\Theta_u\! := \! (\tau_u; \varnothing ; (X^{u}_{\! {v}})_{v\in \tau_u} )$ as follows:  

\smallskip

\noi
$-$ $\tau_u$ is the unique tree such that $k_\varnothing (\tau_u)\! = \! 1$ and  $\theta_{(1)} \tau_u\! = \! \theta_u \tau^*$ (see (\ref{curshitr}) for the definition of $\theta_u \tau^*$). 

\smallskip

\noi
$-$ Since $\overleftarrow{u} \ino \mathtt{Sp}$, there exists $p\ino \bbN$ such that $\overleftarrow{u}\! = \! \booo (p)$ and we set $X^u_{\! \varnothing}\! = \! Y_p$. 

\smallskip

\noi
$-$ For all $v\ino \theta_u\tau$, $X^u_{\! (1)\ast v}\! = \! X_{u\ast v}$ (recall from (\ref{concbila}) 
the definition of the concatenation $u\ast v$). 

\smallskip

\noi
Namely, $\Theta_u$ is the restriction of the branching walk $\fTheta$ to the subtree stemming from $u$, including the spatial position of $\overleftarrow{u}$. Conditionally given $Y\! = \! (Y_p)_{p\in \bbN}$ and 
$ (k_{\varrho (p)} (\tau^*))_{p\in \bbN}$. 
the branching random walks $\Theta_u$, $u\ino \partial \mathtt{Sp}$, are independent; moreover, conditionally given $(Y _p, \ftau_u)$, $\Theta_u$ has law  $Q^{_d}_{^{Y_p, \ftau_u}}$; here $\ftau_u$ stands for $(\tau_u, \varnothing)$ and $p$ is such that $\overleftarrow{u}\! = \! \booo (p)$. 
Then, for all $q\ino \bbN^*$, we set 
$$A_q\! = \! \big\{ \rmo \!  \notin \! \big\{ X_v\, ; v\ino \tau^* \backslash  \{\varrho \}\! : \! |v\! \wedge \! \varrho| \! \geq \! - q \big\} \big\} \quad \textrm{and} \quad \ttH^\circ_{\rmo} \! = \! \inf \{p\! \geq\!  1 \! : \! Y_p \! = \! \rmo \} .$$ 
The event $A_q$ decreases as $q\! \rightarrow \! \infty$ to the event that the 
branching walk $(X_v)_{v\in \tau^*}$ only visits $\rmo$ once. The previous independence properties then imply the following. 
\begin{eqnarray*}
 \bP (A_q \, | \, Y, \tau^*) & =&  \un_{\{ \ttH^\circ_{\rmo}> q\} }  \!\! \!\! \!\! \!\! \!\!   \!\!   \!    \prod_{\qquad u\in \partial \mathtt{Sp}: |u| >-q } \!\! \!\! \!\!  \!\!   \!\!   \!\!  
Q^{\bgg}_{Y_{|\overleftarrow{u}|}, \ftau_u}\!  \big( \rmo \notin \{ X^u_v; v\ino \tau_u \backslash \{ \varnothing\} \} \big) \\
\textrm{and thus} \quad  \bP (A_q \, | \, Y, \mathtt{Sp}\cup\partial\mathtt{Sp}) & =&   \un_{\{ \ttH^\circ_{\rmo} > q\} }\!\! \!\! \!\! \!\! \!\!   \!\!   \!     \prod_{\qquad u\in \partial \mathtt{Sp}: |u| >-q } \!\! \!\! \!\! \!\! \!\!      \xi (Y_{|\overleftarrow{u}|}, \rmo) , 
\end{eqnarray*}
where we recall the notation $\xi (x,y)$ from (\ref{defixixi}) in Lemma \ref{momtun}. 
Since $Y$ and 
the r.v.~$k_{\boo (p)} (\tau^*)$ are independent, we get 
\begin{eqnarray*}
 \bP (A_q \, | \, Y) & = &  \un_{\{ \ttH^\circ_{\rmo} > q\} } \bE \big[ \xi (\rmo, \rmo)^{k_\boo (\tau^*) } \big] \!\!\!  \prod_{1\leq p \leq q } \!\! \bE \big[ \xi (Y_p, \rmo)^{k_{\varrho (p)}(\tau^*) -1} \big| Y \big]  \\
 & =& \un_{\{ \ttH^\circ_{\rmo} > q\} } \phi\big(  \xi (\rmo, \rmo)\big) \!\!\!  \prod_{1\leq p \leq q } \!\! \psi \big( \xi (Y_p, \rmo) \big) ,
\end{eqnarray*}
where we have set 
\begin{equation*}
\forall r \ino [0,1], \quad  \phi (r)= \sum_{k\in \bbN} r^k \mu (k) \quad \textrm{and} \quad \psi (r)\! = \! \sum_{k\in \bbN} r^k  \overline{\mu} (k+1)= \frac{1\! -\! \phi (r)}{1-r}. 
\end{equation*}
Notice that as $q\! \rightarrow \! \infty$, the event $\{ \ttH^o_{\rmo} > q\} $ decreases to the event $\{ \ttH^o_{\rmo} =\infty\}$, which has strictly positive probability by Lemma \ref{dlczidec} $(i)$. Then, we easily see that (\ref{lpomlpo}) holds true if 
\begin{equation}
\label{kravnann}
\sum_{p\in \bbN} \big( 1\! -\! 
 \psi  \big( \xi (Y_p, \rmo) \big)\big)  <\!  \infty. 
\end{equation}
But by Lemma \ref{momtun}, $1\! - \! \xi (Y_p, \rmo) \! \leq \! \tfrac{2\bgg }{\bgg - 1} \bgg^{\, |Y_p\wedge \rmo|}$ since 
$|\rmo| \! =\!  0$. By (\ref{flablel}) in Lemma \ref{abldnez}, there are two constants $C, C^\prime\ino (0, \infty)$ 
such that if $-|Y_p \wedge \rmo | \! \geq \! C$, we get 
$$ 1\! -\! 
 \psi  \big( \xi (Y_p, \rmo) \big) \leq   \frac{C_{\beta, \mu}}{\big( \log \big( \frac{\bgg-1}{2\bgg} \bgg^{ -|Y_p\wedge \rmo|} \big)\big)^{1+\beta}}\,  \leq \, \frac{\!\! C^\prime}{\big(1 \! - \! |Y_p \wedge \rmo | \big)^{1+ \beta}} $$
By Proposition \ref{quantrans} with $f(x)\! =\! (1+x)^{-1-\beta}$, if $\beta \! >\! 1$, then a.s.~$\sum_{p\geq 0} (1 \! -\!  |Y_p \wedge \rmo |)^{-1-\beta} \! < \! \infty$, which implies (\ref{kravnann}) since 
$|Y_p \wedge \rmo| \! \rightarrow \! -\infty$, as $p\! \rightarrow \! \infty$ as stated in Lemma \ref{dlczidec} $(i)$. This completes the proof of the proposition. \cqfd 

\medskip

We conclude this section with a general estimate for recurrent biased random walks on a (deterministic) rooted ordered tree $T$ that is infinite. More precisely, we fix 
$\lambda \ino (1, \infty )$ and we denote by $(Y_n)_{n\in \bbN}$ the $\lambda$-biased RW whose transition probabilities are given for all $x, y\ino T$ by 
 \begin{equation} 
\label{homsicdef}
\bP (Y_{n+1}\! = \! y \, | \, Y_n\! = \! x)= \left\{ 
\begin{array}{ll}
\frac{\lambda}{\lambda + k_x(T)}&  \textrm{if $y\! = \! \overleftarrow{x} $ and $x\! \neq\!  \varnothing$,} \\
\frac{1}{\lambda + k_x(T)}&  \textrm{if $x\! = \! \overleftarrow{y}$ and $x\! \neq\!  \varnothing$, } \\
 \frac{1}{ k_\varnothing(T)}&  \textrm{if $x\! = \! \overleftarrow{y} \!= \! \varnothing$}
\end{array} \right.
\end{equation}
and $\bP (Y_{n+1}\! = \! y \, | \, Y_n\! = \! x)\! = \! 0$ otherwise. Here, recall that $k_x(T)$ is the number of children of $x$ in $T$. We also recall that $|x|$ is the height of $x$ in $T$, that $x\! \wedge \! y$ is the most recent common ancestor of $x$ and $y$ and that $d_{\mathtt{gr}}$ stands for the graph-distance on $T$: $d_{\mathtt{gr}} (x,y)\! = \! |x|+ |y| \! -\! 2 |x\wedge y|$. 
\begin{lemma}
\label{distcontr} We keep the above notations. We assume that $Y$ is recurrent. Then, for all $\ell, n_1, n_2\ino \bbN$ such that $n_1 \! < \! n_2$, we get 
 \begin{equation} 
\label{distesti}
\bP \Big( |Y_{n_1}|+|Y_{n_2}| \! -\! 2\!\!\!\!\!  \min_{\; n_1 \leq n \leq n_2} \!\!\! \!\! |Y_n| \, \geq \, 2\ell + d_{\mathtt{gr}} (Y_{n_1},Y_{n_2}) \Big) \leq (n_2\! -\! n_1)\frac{ \lambda \! -\! 1 }{\lambda^{\ell} \! -\! 1} \; .
\end{equation}
\end{lemma}
\noi
\textbf{Proof.} For all $x\ino T$, we set $\ttH^\circ_x \! = \! \inf \{ n\ino \bbN^*\! : \! Y_n \! = \! x\}$ and $\pi (x)\! = \! (\lambda + k_{x} (T))\lambda^{-|x|}$ if $x\! \neq \! \varnothing$ and $\pi (\varnothing)\! = \! k_\varnothing (T)$. Note that $\pi$ is an invariant measure and  
standard arguments on electrical networks and random walks on graphs imply that for all distinct $x, y \ino T$, 
$\pi (x)\bP (\ttH^\circ_x \! < \! \ttH^\circ_y \, | \, Y_0 \! = \! y)$ is the effective conductance between 
$x$ and $y$ (see R.~Lyons and Y.~Peres \cite{LyoPerbook} p.~25). If {$y\! \in\,  \rgeo \varnothing , x\rgeo $}, then 
it easily implies $\bP (\ttH^\circ_x \! < \! \ttH^\circ_{^{\overleftarrow{y}}} \, | \, Y_0 \! = \! y)\! = \! (\lambda \! -\! 1) / (\lambda^{|x|-|y|+1} \! -\! 1) $. Similarly, $\bP (\ttH^\circ_x \! < \! \ttH^{\circ}_{^{\varnothing}} \, | \, Y_0 \! = \! \varnothing) \! = \! k_\varnothing (T)^{-1} (\lambda \! -\! 1) / (\lambda^{|x|} \! -\! 1) $. Thus, for all distinct $x, y \ino T$ such that 
$y\ino \lgeo \varnothing , x\rgeo $, 
\[
\bP (\ttH^\circ_x \! < \! \ttH^\circ_{^{\overleftarrow{y}}} \, |\,  Y_0 \! = \! y) \leq 
(\lambda \! -\! 1) / (\lambda^{|x|-|y|} \! -\! 1).
\]
We next define the following sequence of stopping times $(\sigma_k)_{k\in \bbN}$ by setting 
$\sigma_0\! = \! 0$ and $\sigma_{k+1}\! = \! \inf \{ n\! >\! \sigma_k : |Y_n| \! \leq \! (|Y_0| \! -\! k \! -\! 1)_+ \}$ that are a.s.~finite since $Y$ is recurrent. We then consider the (possibly empty) event $A_k \! = \! \big\{ \exists n\ino \{ \sigma_k, \ldots , \sigma_{k+1} \! -\! 1\}: Y_n \ino \lgeo \varnothing , Y_0\rgeo \; \textrm{and} \; |Y_n| \! \geq \! \ell + (|Y_0| \! -\! k \! -\! 1)_+ \big\} $. We fix $z\ino T$ such that $|z|  \! \geq \! \ell +(|z| \! -\! k \! -\! 1)_+$, and let $x, y\ino \lgeo \varnothing , z \rgeo$ be such that $|x|\! = \! \ell +(|z| \! -\! k \! -\! 1)_+$ and $|y| \! = \!  (|z| \! -\! k \! -\! 1)_+$. 
By the strong Markov property at time $\sigma_k$ and by the previous inequality for hitting times we get 
$\bP \big( A_k \cap \{ Y_0 \! = \! z \} \big) \leq (\lambda \! -\! 1) / (\lambda^{\ell} \! -\! 1)\bP (Y_0\! = \! z) $. 
Since on $A_k$, we have $|Y_0| \! \geq  \!  \ell +(|Y_0| \! -\! k \! -\! 1)_+$, the previous inequality {implies} 
$\bP ( A_k ) \! \leq \! (\lambda \! -\! 1) / (\lambda^{\ell} \! -\! 1)$. 

  Next observe that $|Y_0| + |Y_n| \! -\! 2\min_{0\leq m\leq n} |Y_m| \! - \! d_{\mathtt{gr}} (Y_0, Y_n)\! = \! 2( |Y_0\! \wedge \! Y_n| \! -\! \min_{0\leq m\leq n}\!  |Y_m| )$ and that 
$$ \big\{\,  |Y_0\! \wedge \! Y_n| \! -\!  \!\!\min_{0\leq m\leq n}  \!\! |Y_m| \geq \ell\,  \big\} \subset A_0 \cup A_1 \cup \ldots \cup A_{n-1} \; .$$
Thus, $\bP (|Y_0| + |Y_n| \! -\! 2\min_{0\leq m\leq n} |Y_m| \! \geq  \! 2\ell +d_{\mathtt{gr}} (Y_0, Y_n)) \leq n(\lambda \! -\! 1) / (\lambda^{\ell} \! -\! 1)$ and we get (\ref{distesti})
 by the Markov property at time $n_1$. \cqfd

\subsection{Law invariance.}
\label{lawinvsec}
In this section, we first define a successor map for $\overline{\bbW}_{[0, 1]}$-labelled trees that centers the spatial positions \textit{and} the genealogical tree at the individual coming next in the lexicographical order, generalizing Definition \ref{ggloupin} for the pointed $\overline{\bbW}_{[0, 1]}$-labelled trees. 
We then show that free branching random walks are invariant in law under this successor map. 
 \begin{definition}
\label{gloupin} Let $\Theta\! = \! (t,\booo \, ; \bx\! = \! (x_v)_{v\in t})$ be a labelled tree in $\overline{\bbT}^\bullet \! (\overline{\bbW}_{[0, 1]})$ or in {$\overline{\bbT}^\bullet \! (\overline{\bbW}_{\! \bgg})$}. 
\begin{compactenum}

\smallskip

\item[$(a)$] We set $\mathtt{cent} (\Theta)\! = \! (\varphi_{|\booo|}(t), \varphi_{|\booo|} (\booo) ; (\varphi_{|x_{\booo}|} (x_{\varphi_{-|\booo|} (v)}) )_{v\in \varphi_{|\booo|}(t)} )$.  
We keep calling $\mathtt{cent} (\cdot)$ the 
\textit{centering} map. As already mentioned, the spatial marks are also shifted. 

\smallskip

\item[$(b)$] We next set $\mathtt{scc} (\Theta ) \! = \! \mathtt{cent} (t, \mathtt{scc} (\booo); \bx)$ where we recall that 
$ \mathtt{scc} (\booo)$ stands for the vertex of $t$ coming next in the lexicographical order as defined by (\ref{lextredef}). 
We keep calling $\mathtt{scc} (\cdot)$ the \textit{successor map}.

\smallskip

\item[$(c)$] Recall from Definition \ref{ggloupin} $(c)$ the right-part $[(t, \boo)]^+$ of $(t,\boo)$. 
To simplify the notation, we set $(t^\prime, \booo^\prime) \! = \! [(t, \booo)]^+$ and we  
recall from Definition \ref{gloupin} $(c)$ that there exists a unique one-to-one map 
$\psi \! : \! t^\prime \! \rightarrow \! \{\booo_{| p } \, ; \,  p\! \leq \! |\booo| \} \cup  \{ v\ino t\! : \!  \booo \! <_t \! v \}$ such that $\psi (\varrho^\prime)\! = \! \varrho$ 
that preserves the relative height and 
that is increasing with respect to the lexicographical order; 
then, we set $[\Theta]^+\! = \! (t^\prime, \booo^\prime, \bx^\prime\! = \! (x^\prime_v)_{v\in t^\prime})$ where $x^\prime_v\! = \! x_{\psi (v)}$, for all $v\ino t^\prime$. We keep calling $[\Theta]^+$
the \textit{right-part} of $\Theta $.

\smallskip

\item[$(d)$] We set $\mathtt{scc}^+ (\Theta) \! = \! [\mathtt{scc} (\Theta) ]^+\! 
$; we keep calling $\mathtt{scc}^+ (\cdot)$ the \textit{right-successor} map. \cq

\end{compactenum}
\end{definition}

We next explain a way to generate a free branching random walk from an i.i.d.~field that is suited to the successor map. To that end we fix $\bt\! = \! (t,u) \ino \overline{\bbT}^\bullet$ such that 
$|t|_{\! -}\! = \! -\infty$. Let $\mathbf{o}$ be a $\bbZ_-$-indexed sequence of mutually independent uniform r.v.~on $[0, 1]$. Let $\fTheta$ be a $\Wbf$-valued branching random walk that has conditional law $Q^{_{\!}}_{\mathbf{o}, \bt}$ given $\mathbf{o}$ as in Definition \ref{defimodel} $(iv)$. To simplify we denote by $\bQ_\bt$ the (unconditional) law of $\fTheta$.      
Let $\epp_v$, $U_v$, $v\ino t$ be independent r.v.~such that $U_v$ is uniformly distributed on $[0, 1]$ and $\bP (\epp_v\! = \! 1) \! =\!  \bP (\epp_v \! = \! -1)  \! = \! 1/2$. 
From the field of i.i.d.~r.v.~$\fTheta_0\! := \! (t,u; (\epp_v, U_v)_{v\in t})$, we now explain how to 
construct a branching random walk that has law $\bQ_\bt$. 

Recall that $\lgeo v, w \rgeo$ stands for the shortest path (with respect to the graph-distance) that joins $v$ and $w$ in $t$ (along with the previous notation $\lgeo v, w \lgeo$, $\rgeo v , w \rgeo$ and $\rgeo v, w \lgeo$, introduced at the beginning of Section \ref{memesec}). Recall also that $\rgeo\!  -\! \infty, v\rgeo$ stands for the lineage of $v$: namely, 
$\rgeo \! -\! \infty, v\rgeo \! =\!  \{ v_{| l}\, ; l \! \leq \! |v| \}$. 
Recall from (\ref{comancdef}) the definition of the most recent common ancestor $v \wedge w$ in $t$. 
We first set
\begin{equation*}
\forall v\in t , \quad H_v =\!\!\!\!  \!\!\!\! \!  \sum_{\quad w\in \lgeo u, u\wedge v\lgeo }\!\!\!\!\!\! \!\!\! (-\epp_w) \, + \!\!\!\!  \!\!\!\! \! \sum_{\quad w\in \lgeo u\wedge v,  v \lgeo} \!\!\!\!  \!\!\!\! \! \epp_w \; , 
\end{equation*}
with the convention that $H_u \!= \! 0$.  
Note that $H$ takes arbitrary negative values on the lineage of $v$. Namely $H$ satisfies (\ref{ZBrRW}). It is easy to see that $(t, u; (H_v)_{v\in t})$ is a $\bbZ$-valued branching random walk whose spatial motion is that of a simple symmetric random walk on $\bbZ$  and whose initial position is a.s.~$H_u\! = \! 0$. Recall from (\ref{sdfljbcv}) the notation $v(k)$ for all integers $k\! \leq \! H_v$: namely,   
$v(k) \! \in \,  \rgeo \! -\! \infty, v\rgeo$ is such that $H_{v(k)}\! =\!  \min_{ w\in \lgeo v(k), v\rgeo} H_w \! =\! 
 k \! >\!  H_{\overleftarrow{v(k)}}$.

 As in (\ref{fgorth}), for all $v\ino t$ we set $X_v \! = \! (U_{v(k)})_{k\leq H_v}$. 
We denote $X_u$ by $\mathbf{o}$. It shows that there exists a deterministic map $F$ such that
\begin{equation} 
\label{lawalaw}
F( \fTheta_0 ) := \big( t, u; (X_v)_{v\in t}\big) \;  \textrm{, which has law $\bQ_\bt$.}
\end{equation}
Let $u^\prime \ino t$ and set $F\big( t, u^\prime ; (\epp_v, U_v)_{t\in t} \big)\! = :\! 
(t, u^\prime, (X^\prime_v)_{v\in t})$. Then, set $H^\prime_v \! = \! |X^\prime_v|$, $v\ino t$. 
We deterministically check that 
\begin{equation} 
\label{shiftospa}
\forall v\ino t , \quad H^\prime_v -H_v= -H_{u^\prime} \quad \textrm{and} \quad X^\prime_v = \varphi_{|X_{u^\prime}|} (X_v) \; .
\end{equation}

We next define $[\fTheta_0]^+$ as in Definition \ref{gloupin} $(c)$: to simplify notation, we set 
$(t^\prime, u^\prime) \! = \! [(t, u)]^+$ and we  
recall from Definition \ref{ggloupin} $(c)$ that there exists a unique one-to-one map
$\psi \! : \! t^\prime \! \rightarrow \! \{u_{| p } \, ; \,  p\! \leq \! |u| \} \cup  \{ v\ino t\! : \!  u \! <_t \! v \}$ such that $\psi (u^\prime)\! = \! u$ 
that preserves the relative height and that is increasing with respect to the lexicographical order. Then, we simply define 
$[\fTheta_0]^+ \!\! = \! \big(  t^\prime, u^\prime  , (\epp^\prime_{v}, U^\prime_v)_{v\in t^\prime} \big) $ where $(\epp^\prime_v, U^\prime_v)\! = \! (\epp_{\psi (v)}, U_{\psi (v)})$ for all $v\ino t^\prime$ and also set 
\begin{equation} 
\label{FTheplus}
\widetilde{\mathtt{scc}}^+ (\fTheta_0) = \big[ \mathtt{scc} (t,u); (\epp_v, U_v)_{{v\in t}} \big]^+ . 
\end{equation}  
By (\ref{shiftospa}) and since $X_v$ is a deterministic function of $((\epp_w, U_w), w\! 
 \in  \rgeo \! -\! \infty, u\rgeo\, \cup \, \rgeo \! -\! \infty, v\rgeo)$, it is easy to check deterministically 
\begin{equation}
\label{succplus}
F\big( \widetilde{\mathtt{scc}}^+ (\fTheta_0)\big)= \mathtt{scc}^+ (F(\fTheta_0)) \; .
\end{equation}
Then, (\ref{atschift}), (\ref{lawalaw}), (\ref{shiftospa}) and (\ref{succplus}) 
immediately entail the following lemma.   
\begin{lemma}
\label{erracina} Let $ \bt = (t,u) \ino \overline{\bbT}^\bullet$. 
We assume that $|t|_{\! -}\! = \! -\infty$. Let $\fTheta\! = \! (t, u; (X_v)_{v\in t})$ have law $\bQ_\bt$.    
Then $\mathtt{scc} (\fTheta)$ has law $\bQ_{\mathtt{scc} (\bt)}$ and $\mathtt{scc}^+ (\fTheta)$ has law 
$\bQ_{\mathtt{scc}^+ (\bt)}$. 
\end{lemma}
The next result is the key point in the proofs of the various laws of large numbers 
for the range of branching random walks that we prove in {the next section}. 
\begin{theorem}
\label{rrwergo}
Let $\mu$ be a probability measure on $\bbN$ such that $\sum_{k\in \bbN} k\mu (k) \! = \! 1$. 
Let $\mathbf{o}\! := \! (U_k)_{k\in \bbZ_-}$ be a sequence of independent uniform r.v.~on $\{ 1, \ldots, \bgg\}$. Let $\fTheta\! = \! (\tau^*\! , \, \varrho\, ; \, (X_v)_{v\in \tau^*})$ be a $\Wbd$-valued branching random walk whose distribution is the following: $\ftau^*\! = \! (\tau^*\! , \boo)$ is an IPGW$^+\! (\mu)$-tree as in Definition \ref{GWptdef} $(ii)$ and conditionally given $\ftau^*\! $ and $\mathbf{o}$, 
$\fTheta$ has law $Q^{_\bgg}_{^{\mathbf{o}, \ftau^*}}$ as in Definition \ref{defimodel} $(iii)$. Then 
\begin{equation}
\label{statioon}
\mathtt{scc}^+ (\fTheta) \overset{\textrm{(law)}}{=} \fTheta \; ,
\end{equation}
where we recall that $\mathtt{scc}^+$ stands for the right-successor map as in Definition \ref{gloupin} $(d)$. 

Next, denote by $(v_n)_{n\in \bbN}$ the sequence of the vertices that are not direct ancestors of $\boo$ listed in the lexicographical order. Namely, $v_0\! = \! \boo$, 
$v_n\! <_{\tau^*}\! v_{n+1}$ and $\{ v_n; n\ino \bbN\}\! = \! \{ v\ino \tau \! : \! \booo \! \leq_{\tau^*} \! v \}$. 
Then, 
\begin{equation} 
\label{elvjb}
\textrm{$\bP$-a.s.} \quad \frac{_1}{^n} \# \big\{ X_{v_k} ; 1\! \leq \! k\! \leq n \big\} \underset{n\rightarrow \infty}{-\!\!\! -\!\!\! -\!\!\! \longrightarrow} c_{\mu, \bgg}:=  \bP \big( \mathbf{o} \! \notin\! \{ X_{v_n}; n\! \geq \! 1\} \big) . 
\end{equation}
Furthermore, if $\mu$ satisfies (\ref{hypomunu}), then $c_{\mu, \mathtt b}\! >\! 0$.  
\end{theorem}
\noi
\textbf{Proof.} Let $\fTheta_0\! := \! (\tau^*\! ,\boo; (\epp_v, U_v)_{v\in \tau^*} \! )$ be distributed as follows: conditionally given $\ftau^*\! $, $\epp_v$, $U_v$, $v\ino t$ are independent r.v.~such that $U_v$ is uniformly distributed on 
$[0, 1]$ and $\epp_v$ is uniformly distributed on $\{ -1, 1\}$. 
Then, we set $F(\fTheta_0)\! = \! \fTheta^\prime\! = \! (\tau^*\! , \boo ; (X^\prime_v)_{v\in \tau^*})$ 
as in (\ref{lawalaw}). 
Conditionally given $\ftau^*\!$, $\fTheta^\prime$ has law $\bQ_{\ftau^*}$ and Remarks \ref{tracage} allows to take $\fTheta \! = \! \Phi_\bgg (\fTheta^\prime)$ that has the desired law. 
Lemma \ref{erracina} combined with proposition \ref{succinv} imply that 
$\mathtt{scc}^+(\fTheta^\prime)$ has the same law as $\fTheta^\prime$ and since obviously 
$\Phi_\bgg (\mathtt{scc}^+(\fTheta^\prime))\! = \!   \mathtt{scc}^+(\fTheta)$, we get (\ref{statioon}).

  For all $n\! \geq \! m \! \geq \! 0$, we set $R_{m,n} (\fTheta )\! = \! \# \big\{ X_{v_k} ; m\! <\! k\! \leq n \big\}$ and we denote by $\mathtt{scc}^+_\ell$ the $\ell$-th iterate of $\mathtt{scc}^+$. 
Observe that 
$R_{m,n} (\mathtt{scc}^+_\ell (\fTheta))\! = \! R_{m+ \ell , n+ \ell} (\fTheta)$, with an obvious notation. 
Then, the r.v.~$(R_{m,n} (\fTheta^+)_{m\geq n\geq 1})$ 
satisfy Assumptions (1.7), (1.8) and (1.9) of Liggett's version of Kingman's 
subadditive ergodic theorem (see Theorem 1.10 in \cite{Lig85} p.~1280)
that asserts that there exists a $[0, 1]$-valued r.v. $R$ such that 
$R_{0, n} (\fTheta)/ n\! \rightarrow \! R$ almost surely (and in $L^1$). 

We next prove that $R$ is a.s.~constant. Recall from (\ref{FTheplus}) the definition of $\widetilde{\mathtt{scc}}^+ (\fTheta_0)$ and we denote by $\widetilde{\mathtt{scc}}^+_\ell $ the $\ell$-th iterate of $\widetilde{\mathtt{scc}}^+ $: namely, 
$\widetilde{\mathtt{scc}}^+_\ell  (\fTheta_0) \! =\! [ (\tau^*, v_\ell); (\epp_v, U_v)_{v\in \tau^*} ]^+ $. Denote by $\ccG_\ell$ the sigma-field generated by $\widetilde{\mathtt{scc}}^+_\ell  (\fTheta_0)$. As a consequence of (\ref{succplus}), $\mathtt{\mathtt{scc}}^+_\ell (\fTheta)\! = \! \Phi_d (F (\widetilde{\mathtt{scc}}^+_\ell  (\fTheta_0))$, the r.v.~$(R_{\ell, \ell+n} 
(\fTheta))_{n\in \bbN}$ are $\ccG_\ell$-measurable. Note that the sigma-fields $\ccG_\ell$ decrease in $\ell$. Next, we set $\ccG\! = \! \bigcap_{\ell \in \bbN} \ccG_\ell$ that can be viewed as the tail sigma-field of the subtrees grafted on the infinite line of ancestors; since additional marks $(\epp_v, U_v)$ are i.i.d., Kolmogorov's zero-one law applies and $\ccG$ is $\bP$-trivial. Furthermore, the subadditivity for the $R_{m, n}$ entail that a.s.~$R\! = \! \lim_{n\rightarrow  \infty} n^{-1}R_{\ell , \ell + n} (\fTheta)$, for all $\ell$. Thus,   
$R$ is $\ccG$-measurable which implies that it is a.s.~constant.

Consequently, $R\! = \! \lim_{n\rightarrow \infty}\bE[R_{0, n} (\fTheta)]/n$. By an elementary argument, 
$$\bE[R_{0, n} (\fTheta)]\! = \! \sum_{1\leq k\leq n} \! \bP \big(X_{v_{k}} \!  \notin \! \{  X_{v_{k+1}}, \ldots, X_{v_n} \}  \big) \; .$$ 
Since the law of $\fTheta$ is preserved by the map $\mathtt{scc}^+ (\cdot)$, we get $\bP \big(X_{v_k} \! \notin \! \{ X_{v_{k+1}}, \ldots, X_{v_n} \}  \big)\! = \! \bP \big( X_{\booo} 
\! \notin \!  \{ X_{v_{1}}, \ldots, X_{v_{n-k}} \}  \big)$ and thus 
$$ \frac{1}{n} \bE[R_{0, n} (\fTheta)]\! =   \frac{1}{n} \! \!\!\!  \sum_{\; \; 0\leq k\leq n-1}\!\! \!  \bP \big( X_{\booo} 
\! \notin \!  \{ X_{v_{1}}, \ldots, X_{v_{k}} \}  \big) \underset{n\rightarrow \infty}{-\!\!\! -\!\!\! -\!\!\! \longrightarrow} c_{\mu, \bgg} $$
by Cesàro. Finally, observe that $c_{\mu, \bgg} \! \geq \! \kappa_{\mu, \bgg}$, where $\kappa_{\mu, \bgg}$ is defined in (\ref{lpomlpo}) in Proposition \ref{ntgngbck} that completes the proof of the theorem. \cqfd

\subsection{Proof of Theorem \ref{main1}.}
\label{Pfmain1}
We fix $\gamma \ino (1, 2]$. Let $\tau$ be a GW($\mu$)-tree whose offspring distribution $\mu$ satisfies $(\mathbf{H})$ in (\ref{hyposta}). We set $\ftau\! = \! (\tau, \varnothing)$.  
Let $\mathbf{o}\! := \! (U_k)_{k\in \bbZ_-}$ be a sequence of independent uniform r.v.~on $\{ 1, \ldots, \bgg \}$. Let the $\Wbd$-valued branching random walk
$\fTheta_{d} \! 
= \! (\tau\! , \, \varnothing\, ; \, (Y_v)_{v\in \tau})$ have conditional law $Q^{\bgg}_{{\mathbf{o}, \ftau}}$ given $\tau$ and $\mathbf{o}$ as in Definition \ref{defimodel} $(iii)$. 
We also introduce the $\bbW_{\! \bgg}$-valued branching random walk $\fTheta^+_{\bgg} \! 
= \! (\tau\! , \, \varnothing\, ; \, (Y^+_v)_{v\in \tau})$ that has conditional law $Q^{+\bgg}_{{\ftau}}$ given $\tau$ as in Definition \ref{defimodel} $(i)$. 
For all integers $0 \! \leq \! k \! < \! \# \tau$, we denote 
\begin{equation}
\label{rangegeg}
R_k (\fTheta_{\bgg}) \! = \! \# \{ Y_v \, ; v\ino \tau\! : v\! \leq_\tau \! v_k \} \quad \textrm{and} \quad R_k (\fTheta^+_{\bgg}) \! = \! \# \{ Y^+_v \, ; v\ino \tau\! : v\! \leq_\tau \! v_k \}, 
\end{equation}
 where $v_k$ stands for the $k$-th smallest vertex of $\tau$ with respect to the lexicographical order $\leq_\tau$. Arguing as in Le Gall and L.~\cite{LGLi16}, we derive the following result from Theorem \ref{rrwergo}. 
\begin{proposition}[adapted from Theorem~7 \cite{LGLi16}]
\label{JFLGSLin} 
Let $\gamma \ino (1, 2]$ and let $\mu$ satisfy $(\mathbf{H})$ as in (\ref{hyposta}). 
We keep the above notation and we recall from (\ref{elvjb}) 
the definition of $c_{\mu, \bgg}$ and from (\ref{rangegeg}) the definition of $R_n (\fTheta_{\bgg})$. Then,  
\begin{equation}
\label{uniprofr}
\forall \epp \ino (0, \infty), \quad \bP \Big( \! \! \!  \sup_{\; \; \;  t\in [0, 1]}\!   \tfrac{1}{n}\big| R_{\lfloor nt \rfloor} (\fTheta_{\bgg}) -c_{ \mu, \bgg} \, nt  \big| >\epp \, \Big| \, \# \tau \! = \! n  \Big) \underset{n\rightarrow \infty}{ -\!\!\! -\!\!\! -\!\!\! \longrightarrow} 0 \; .
\end{equation}
\end{proposition}
\noi
\textbf{Proof.} The proof can be adapted verbatim from the way 
Proposition~6 and Theorem~7 pp.~284-289 in \cite{LGLi16} are deduced from Proposition~3 and Theorem~4 pp.~280-284 {in} \cite{LGLi16}. Here our Theorem \ref{rrwergo} plays the role of Proposition~3 {and} Theorem~4 
in \cite{LGLi16} and we only give a brief sketch of the proof. Following the same arguments as in Proposition 6 in \cite{LGLi16},
we first prove that for any fixed $s\ino [0,1]$, 
\begin{equation*}
\forall \epp \ino (0, \infty), \quad \bP \Big(   \tfrac{1}{n}\big| R_{\lfloor ns \rfloor} (\fTheta_{\bgg}) -c_{ \mu, \bgg} \, ns  \big| >\epp \, \Big| \, \# \tau \! > \! n  \Big) \underset{n\rightarrow \infty}{ -\!\!\! -\!\!\! -\!\!\! \longrightarrow} 0 \; .
\end{equation*}
From this limit and the same absolute continuity argument as in the proof of Theorem 7 of \cite{LGLi16},
we get for any fixed $s\ino [0, 1]$,  
\begin{equation*}
\forall \epp \ino (0, \infty), \quad \bP \Big(   \tfrac{1}{n}\big| R_{\lfloor ns \rfloor} (\fTheta_{\bgg}) -c_{ \mu, \bgg} \, ns  \big| >\epp \, \Big| \, \# \tau \! = \! n  \Big) \underset{n\rightarrow \infty}{ -\!\!\! -\!\!\! -\!\!\! \longrightarrow} 0 \; .
\end{equation*}
\noi
We then get (\ref{uniprofr}) thanks to the following variant of the second Dini's theorem. \cqfd 
\begin{lemma}
\label{probDini2} For all $n\ino \bbN$, let $(\Omega_n, \ccF_{\! n}, \bP_{\! n})$ be a probability space on which a nondecreasing right-continuous process $(X^{_{(n)}}_{^t})_{t\in [0, \infty)}$ is defined. Let $(x(t))_{t\in [0, \infty)}$ be a (deterministic) continuous function. Suppose that for all $t\ino [0, \infty)$, the real valued r.v.~$X^{_{(n)}}_{^t}\! $ under $\bP_{\! n}$  tends to $x(t)$ in law. Then, for all 
$\epp \ino (0, \infty)$ and all $p\ino \bbN$, $\lim_{n\rightarrow \infty} 
\bP_{\! n} \big( \sup_{t\in [0, p]} |X^{(n)}_t \! -\! x(t)| \! >\! \epp \big) \! = \! 0$. 
\end{lemma}
\noi
\textbf{Proof of Lemma \ref{probDini2}.} We fix $p\ino \bbN$. For all $n, q\ino \bbN$, set $w_{n,q}\! = \! \max_{0\leq k\leq p2^q} |X^{(n)}_{k2^{-q}} \! -\! x(k2^{-q})|$. Since the convergence in law to a constant 
implies the convergence in probability, we easily get for all $q$ that $\lim_{n\rightarrow \infty} \bP_{\! n}(w_{n,q}\! >\! \epp)  = \! 0$. Next set $v_q\! = \! \max \{ |x(t)\! -\! x(s)|; s, t\ino [0, p] \! : \! |s\! -\! t| \! \leq \! 2^{-q} \}$ that tends to $0$ {as $q\to \infty$} since $x$ is continuous. {By monotonicity of $X^{(n)}_t$,} it is then easy to see that for all $t\ino [0, p]$, that $|X^{(n)}_t \! -\! x(t)| \! \leq\! 2v_q + 3 w_{n,q}$, which easily implies the desired result. \cqfd

\medskip

\noi
\textbf{End of the proof of Theorem \ref{main1}.} We keep the notation from above and recall from (\ref{elvjb}) 
the definition of $c_{\mu, \bgg}$ and from (\ref{rangegeg}) the definition of $R_n (\fTheta^+_{\bgg})$. We now derive Theorem \ref{main1} 
from (\ref{uniprofr}) and  Proposition \ref{coouver} that allows to change the state space from $\overline{\bbW}_{\! \bgg}$ to $\bbW_{\! \bgg}$.

For all positive integers $n $, we set $c_n\! = \! 4 \log_\bgg n$, 
$L_n\! = \! \# \{ v\ino \tau: |Y^+_v| \! \leq \! c_n \}$, $N_n\! = \! \# \{ v\ino \tau\! : |v|\! = \! n\}$ and $C_n\! = \! N_0+ N_1+ \ldots + N_n$. Let $(S_n)_{n\in \bbN}$ be a simple symmetric random walk on $\bbZ$ with initial position $S_0\! =\! 0$. By definition of $\fTheta^+_{\bgg}$, 
for all $v\ino \tau$, we get $\bE  [ f(|Y^+_v|) \, | \, \tau ]\! = \! \bE [f( |S_{|v|}|)\,  |\,  \tau]$. Then, for all $b_n\ino \bbN^*$ we get 
\begin{equation}
\label{floupii}
 \bE \big[ L_n \, \big| \, \tau\big]  = \sum_{k\in \bbN} N_k \bP \big( |S_k| \! \leq  \! c_n \big) \leq  C_{b_n} +  \sigma (b_n, \tfrac{c_n}{\sqrt{b_n}}) \# \tau \; , 
 \end{equation}
where we have set $\sigma (p,s)\! = \! \sup_{k\geq p} \bP \big( |S_k| k^{-\frac{1}{2}} \! \! \! \leq \! s \big)$ for all $p\ino \bbN^*$ and all $s\ino (0, \infty)$.
Observe that any fixed $p$, $\lim_{s\rightarrow 0+}\sigma (p, s) \! =\! 0$.   
We now take $b_n\! = \! c_n^{4} $.

We next claim that for all $\eta \ino (0, \infty)$, 
\begin{equation}
\label{seebald}
\lim_{n\rightarrow \infty} \bP\big( \tfrac{_1}{^n}C_{b_n} \! \geq \! \eta \, \big| \, \# \tau \! = \! n \big) = 0\; .
\end{equation}
\textit{Proof of (\ref{seebald}).} For all $\epp\ino (0, \infty)$ and all $y\ino \bbR$, 
we set $f_\epp(y)\! =\!  (1\! -\! \epp^{-1}(|y|\! -\! \epp)_+)_+$ that is continuous and such that $\un_{[-\epp, \epp]} \! \leq \! f_\epp \! \leq \!  \un_{[-2\epp, 2 \epp]}$. Recall from the introduction that $(H_k)_{0\leq k < \# \tau}$ stands for the height process of the tree $\tau$ and recall the convergence (\ref{Loos}). Then note that 
$$\tfrac{_1}{^n} C_{a_n \epp}\leq Q_\epp^n:=\int_0^1 f_\epp \big( \tfrac{_1}{^{a_n}} H_{\lfloor ns \rfloor}\big) \, ds \; .$$
By (\ref{Loos}), $Q_\epp^n$ under $\bP( \,\cdot \,  | \# \tau \! = \! n)$ tends in law to $Q_\epp\! := \! \int_0^1 f_\epp(H_s) \, ds$. We recall here that a.s.~for all $s\in (0, 1)$, $H_s\! >\! 0$ (see D.~\cite{Du03} or Section \ref{statreesc}). Thus $\lim_{\epp \rightarrow 0+} Q_\epp\! = \! 0$ and we get for all $\eta \ino (0, \infty)$
\begin{eqnarray}
\label{mueller}
 \limsup_{\epp \rightarrow 0+} \limsup_{n\rightarrow \infty}  \bP\big( \tfrac{_1}{^n}C_{a_n\epp } \! \geq \! \eta \, \big| \, \# \tau \! = \! n \big) & \leq & \limsup_{\epp \rightarrow 0+} \limsup_{n\rightarrow \infty}  \bP\big( Q_\epp^n \! \geq \! \eta \, \big| \, \# \tau \! = \! n \big) \nonumber\\ 
 &\leq&   \limsup_{\epp \rightarrow 0+} \bP\big( Q_\epp^n \! \geq \! \eta \, \big)= 0 . 
\end{eqnarray} 
Since $(a_n)_{n\in \bbN}$ is a $\tfrac{_{\gamma-1}}{^\gamma}$-regularly varying 
sequence, $b_n/a_n\! \rightarrow \! 0$ and (\ref{mueller}) entails (\ref{seebald}). \cq 

\medskip

Since $\tfrac{_1}{^n}C_{b_n} \ino (0, 1] $, it also implies that 
$\bE [\frac{_1}{^n} C_{b_n} | \# \tau \! = \! n ] \! \rightarrow \! 0$ as $n \! \rightarrow \! \infty$. By 
(\ref{floupii}) we get $\lim_{n\rightarrow \infty}\bE [\frac{_1}{^n}L_{n} | \# \tau \! = \! n ] \! = \! 0$. 
We next apply Proposition \ref{coouver} to $c \! =\!  \lfloor 4\log_\bgg n \rfloor \! +\! 1$: 
(\ref{vouasx}) implies  the following. 
\begin{eqnarray*}
 &  &\!\!\! \!\!\! \!\!\! \!\!\! \!\!\! \!\!\!   \bP \Big(  \sup_{s\in [0, 1]} \tfrac{1}{n} \big| R_{\lfloor ns \rfloor} (\fTheta^+_{\bgg}) -c_{\mu, \bgg} \, ns  \big| >2\epp \, \Big| \, \# \tau \! = \! n  \Big) \\
 &\leq & 2\bgg n^{-1} + \bP \big( \tfrac{1}{n}L_n \! >\! \epp  \,  \big| \, \# \tau \! = \! n  
) + \bP \Big(  \sup_{s\in [0, 1]}  \tfrac{1}{n}\big| R_{\lfloor ns \rfloor} (\fTheta_{\bgg}) -c_{\mu, \bgg} \, ns  \big| > \epp \, \Big| \, \# \tau \! = \! n  \Big)\; , 
\end{eqnarray*}
which implies (\ref{main1eff}) in Theorem \ref{main1} by (\ref{uniprofr}) since $\bE [\frac{_1}{^n}L_{n} | \# \tau \! = \! n ] \! \rightarrow \! 0$ as $n \! \rightarrow \! \infty$.  \cqfd 

\section{{Snake metrics and the Brownian cactus.}}

\subsection{Pseudo-metrics on a closed interval.}
\label{psdmsec}
\begin{definition}
\label{psddef} Let $\zeta\ino (0, \infty)$. We introduce the following spaces. 
\begin{compactenum}

\smallskip

\item[$(a)$] We denote by $\bC ([0, \zeta]^2, \bbR)$ the space of continuous 
functions from $[0, \zeta]^2$ to $\bbR$ that is a Banach space when equipped with 
uniform norm $\lVert \cdot \rVert$. 

\smallskip

\item[$(b)$] We denote by $\MMM([0, \zeta])$ the set of continuous pseudo-metrics on $[0, \zeta]$. Namely, it is the set of 
$d\ino  \bC ([0, \zeta]^2, \bbR)$ such that for all $s_1, s_2, s_3\ino [0, \zeta]$, 
$$d(s_1, s_2) \! \geq \! 0, \; \,  d(s_1, s_1)\! = \! 0, \; \, d(s_1, s_2)\! = \! d(s_2, s_1) \; \, \mathrm{and} \;\,  d(s_1, s_3)\! \leq \! d (s_1, s_2) + d(s_2, s_3). $$
\item[$(c)$] We denote by $\MMT([0, \zeta])$ the space {of continuous pseudo-metrics} $d\ino \MMM([0, \zeta])$ such that 
for all $s_1, s_2, s_3, s_4 \ino [0, \zeta]$, 
\end{compactenum}
\begin{equation}
\label{fourpoints}
 d(s_1, s_2) + d(s_3, s_4) \leq \max \big( d(s_1, s_3) + d(s_2, s_4) \, ; \,d(s_1, s_4) + d(s_2, s_3) \big) \; .
\end{equation}
\begin{compactenum}
\item[] We call the latter inequality the \textit{four points inequality}. \cq 
\end{compactenum}
\end{definition}
We easily check that $\MMT ([0, \zeta])$ and $\MMM ([0, \zeta])$ are closed subsets of 
$(\bC ([0, \zeta]^2, \bbR), \lVert \cdot \rVert)$. We shall need the following 
compactness criterion in $\MMM([0, \zeta])${.}
\begin{lemma}
\label{ascobis} For all $\eta\ino (0, \infty)$ and all $d\ino \MMM([0, \zeta])$, we set \begin{equation*}
q_\eta (d)= \max \big\{ d(s,s^\prime)\, ;\;  s,s^\prime \ino [0, \zeta] : |s\! -\! s^\prime | \leq \eta \big\} .  
\end{equation*}
Let $\cD$ be a subset of $\MMM([0, \zeta])$. Then, the closure of $\cD$ is compact {if and only if}  $\,\,  \sup_{d\in \cD} q_\eta (d) \longrightarrow 0$ as $\eta \! \rightarrow 0+$. 
\end{lemma}
\noi
\textbf{Proof.} For all $f\ino \bC ([0, \zeta]^2, \bbR)$ and all $\eta\ino (0, \infty)$, 
set 
$$\omega_\eta (f)\! = \!  \max \big\{ |f(s_1,s_1^\prime)\! -\! f(s_2,s_2^\prime) | 
\, ;\;  s_1,s_2, s_1^\prime, s_2^\prime \ino [0, \zeta] : |s_1\! -\! s_2 | \! \vee \! |s_1^\prime \! -\! s_2^\prime | 
 \leq \eta \big\} $$
that is the $\eta$-modulus of uniform continuity of $f$. The Arzel\`a-Ascoli Theorem asserts that the closure of $\cA\! \subset \! \bC ([0, \zeta]^2, \bbR)$ is $\rVert \cdot \rVert$-compact {if and only if}
$ \sup_{f\in \cA} \omega_\eta (f) \! \rightarrow \! 0$ as $\eta \! \rightarrow \!  0+$ and $\sup_{f\in \cA} |f(0, 0)| \! < \! \infty$. The desired result follows from the easy observation that 
$q_\eta (d) \! \leq \! \omega_\eta (d) \! \leq \! 2q_\eta (d)$ (here, the second inequality is a consequence of the triangle inequality) and from the obvious fact that $d(0, 0)\! = \! 0$ for any $d\ino \cD$. \cqfd

\medskip

\noi
This lemma immediately implies the following tightness criterion. 
\begin{proposition}
\label{tightness} 
Let $(\bbd_n)_{n\in \bbN}$ be a sequence of $\MMM([0, \zeta])$-valued random variables. Their laws are tight on $\big( \bC([0, \zeta]^2, \bbR), \lVert \cdot \rVert \big)$ {if and only if} for all $\epp\ino (0, \infty)$, $\sup_{n\in \bbN} \bP \big( q_\eta (\bbd_n) \! >\! \epp \big) \! \longrightarrow \! 0$ as $\eta \! \rightarrow \! 0+$.
\end{proposition}  
We now prove the following specific proposition that is used later in the proof of Theorem \ref{main2}. 
\begin{proposition}
\label{submetric} Let $\bbd_n, \bbd_n^*$, $n\ino \bbN$, be $\MMM([0, \zeta])$-valued r.v. such that 
\smallskip

\begin{compactenum}

\item[$(i)$] $\bbd_n \! \longrightarrow \! \bbd$ weakly on $\big(\bC([0, \zeta]^2, \bbR), \lVert \cdot \rVert\big)$;

\smallskip

\item[$(ii)$] For all $ n\ino \bbN$ and for all 
$ s, s^\prime\ino [0, \zeta]$, a.s.~$\bbd_n^*(s, s^\prime) \! \leq \! \bbd_n(s, s^\prime)$; 

\smallskip

\item[$(iii)$] For all $s, s^\prime\ino [0, \zeta]$, $|\bbd_n^*(s, s^\prime) \! - \! \bbd_n(s, s^\prime)|\! \longrightarrow \! 0 $ in probability.

\smallskip

\end{compactenum}

\noi
Then, $(\bbd_n^*, \bbd_n) \! \longrightarrow \! (\bbd, \bbd) $ weakly on $\big(\bC([0, \zeta]^2, \bbR), \lVert \cdot \rVert\big)^2$. 
\end{proposition} 
\noi
\textbf{Proof.} By $(i)$ and Proposition \ref{tightness}, 
for all $\epp \ino (0, \infty)$, $\lim_{\eta \rightarrow 0} \sup_{n\in \bbN} \bP ( q_\eta (\bbd_n) \! >\! \epp)\! =\! 0$. Since $\bbd_n$ and $\bbd_n^*$ are continuous, $(ii)$ actually implies that a.s.~for all $n\ino \bbN$ and for all $s, s^\prime\ino [0, \zeta]$, 
$\bbd_n^*(s, s^\prime) \! \leq \! \bbd_n(s, s^\prime)$, which immediately entails 
$q_\eta (\bbd_n^*) \! \leq \! q_\eta (\bbd_n)$. Thus, $\lim_{\eta \rightarrow 0} \sup_{n\in \bbN} \bP ( q_\eta (\bbd^*_n) \! >\! \epp)\! =\! 0$, 
for all $\epp \ino (0, \infty)$. By Proposition \ref{tightness}, the laws of the r.v.~$(\bbd_n^*, \bbd_n)$ are tight in $\bC([0, \zeta]^2, \bbR)^2$ and we get the desired result because $(i)$ and $(iii)$ easily {entail} the weak convergence of the finite dimensional marginals of $(\bbd_n^*, \bbd_n)$ to those of $(\bbd, \bbd)$. \cqfd

\medskip

\noi
\textbf{Induced metric spaces.} Let $d\ino \MMM ([0, \zeta])$. We define the relation $\sim_d$ on $[0, \zeta]$ as follows: 
for all $s_1  , s_2\ino [0, \zeta]$, $s_1 \sim_d s_2$ {if and only if} $d(s_1, s_2)\! = \! 0$. Clearly $\sim_d$ is an equivalence relation and we define the quotient space:  
\begin{equation}
\label{quotient}
E_d = [0, \zeta] / \!\! \sim_d , \quad \mathtt{proj}_d: [0, \zeta] \! \rightarrow \! E_d,\; \textrm{the canonical projection} , \; r_d = \mathtt{proj}_d (0){.}
\end{equation}
We keep denoting $d$ the metric induced by $d$ on $E_d$. Since $d$ is continuous on $[0, \zeta]^2$, $\mathtt{proj}_d$ is continuous and $(E_d, d, r_d)$ is a pointed compact and connected metric space. We also equip $E_d$ with the pushforward measure $\mu_d$ of the Lebesgue measure on $[0, \zeta]$ via the canonical projection: namely, for all nonnegative measurable functions $f$ on $E_d$, 
\begin{equation}
\label{induleb}
\int_{E_d}\!\! \! f(x) \, \mu_d (dx)\! = \! \int_0^\zeta \!\!  f(\mathtt{proj}_d (s))\,  ds \; .
\end{equation}
Note that $\mu_d $ is a finite measure with total mass $\zeta$. 

\begin{remark}
\label{Passau} Let $(E_d, d, r_d, \mu_d)$ be the compact metric space corresponding to the pseudo-metric $d\ino \MMM ([0, \zeta])$. Let $a, b\ino (0, \infty)$. Set $d^\prime (s_1, s_2) \! = \! a d (s_1/b, s_2/b)$, $s_1, s_2\ino [0, b\zeta]$. Then, we easily check that $d^\prime \ino \MMM([0, b\zeta ])$ and $(E_{d^\prime}, d^\prime, r_{d^\prime}, \mu_{d^\prime})$ is isometric to 
$(E_d, ad, r_d, b\mu_d)$. \cq 
\end{remark}

\smallskip

\noi
\textbf{Real trees.} When $d\ino \MMT([0, \zeta])$, the resulting space $E_d$ is a real tree. More precisely, real trees are 
metric spaces that extend the definition of graph-trees; they are defined as follows.  
\begin{definition}
\label{errtredef} 
let $(T, d)$ be a metric space; it is a {\it real tree} {if and only if} the following holds true.

\smallskip

\begin{compactenum}
\item[(a)] For any $\sigma_1, \sigma_2 \! \in\!  T$, there is a unique isometry 
$f:[0,d(\sigma_1,\sigma_2)] \! \rightarrow \! T$ such
that $f(0)\!=\! \sigma_1$ and $f(d(\sigma_1,\sigma_2))\! =\! \sigma_2$. Then, we set 
$\lgeo \sigma_1,\sigma_2\rgeo \! :=\! f([0,d (\sigma_1,\sigma_2)])$. 

\smallskip

\item[(b)] For any continuous injective function 
$g: [0, 1] \! \rightarrow \! T$, such that $g(0)=\sigma_1$ and $g(1)=\sigma_2$,  $g([0,1]) \! = \! \lgeo \sigma_1,\sigma_2\rgeo$. \cq 
\end{compactenum}
\end{definition}
 It turns out that the four points inequality is a metric characterisation of real trees. More precisely, if $(T, d)$ is a connected metric space, then $(T, d)$ is a real tree {if and only if} for any $\sigma_1,  \sigma_2,  \sigma_3,  \sigma_4  \in T$,
$$d(\sigma_1, \sigma_2) + d(\sigma_3, \sigma_4) \leq \big(d(\sigma_1, \sigma_3) + d(\sigma_2, \sigma_4)\big) \vee  \big( d(\sigma_1, \sigma_4) + d(\sigma_2, \sigma_3)  \big) . $$
We refer to Evans \cite{Ev08} or to Dress, Moulton and Terhalle \cite{DrMoTe96} for a detailed account on this property.

Let us introduce some notation about real trees. Let $(T,d)$ be a compact pointed real tree. We distinguish a point $r\ino T$ that is viewed as a  \textit{root}. 
We then define the \textit{length} measure $\mathtt{Length}(\cdot) $ on $T$ as the one-dimensional Hausdorff measure: namely, it is the unique Borel measure such that $\mathtt{Length} (\lgeo \sigma, \sigma^\prime \rgeo) \! = \! d(\sigma, \sigma^\prime)$, for all $\sigma, \sigma^\prime\ino T$. Let us next introduce {branch points}: 
let $\sigma_1, \sigma_2, \sigma_3 \ino T$; then the geodesic paths $\lgeo \sigma_1 , \sigma_2 \rgeo$, $\lgeo \sigma_1 , \sigma_3 \rgeo$ and $\lgeo \sigma_2 , \sigma_3 \rgeo$ have exactly one point in common that is called the {branch point of $\sigma_1, \sigma_2, \sigma_3$} and that is denoted by $\mathtt{br} (\sigma_1, \sigma_2, \sigma_3)$; namely 
\begin{equation}  
\label{braachpt}
\{ \mathtt{br} (\sigma_1, \sigma_2, \sigma_3)\}= \lgeo \sigma_1 , \sigma_2 \rgeo \cap \lgeo \sigma_1 , \sigma_3 \rgeo \cap \lgeo \sigma_2 , \sigma_3 \rgeo . 
\end{equation}
If we view $T$ as a family tree whose ancestor is $r$, then $\mathtt{br} (\sigma_1, \sigma_2, r)$ is the most recent common ancestor of $\sigma_1$ and $\sigma_2$ and we use the following notation 
\begin{equation*}  
\sigma_1 \wedge \sigma_2 = \mathtt{br} (\sigma_1, \sigma_2, r).
\end{equation*}
We next introduce the (extended) \textit{degree} of any point 
$\sigma \! \in \! T$ as follows.
\begin{equation}
\label{Traun}
\textrm{$\mathtt{deg} (\sigma)$ is the (possibly infinite) number of connected components of 
the open set $T\backslash \{ \sigma\}$.}
\end{equation} 
We say that $\sigma$ is a \textit{branch point} if 
$\mathtt{deg}(\sigma) \! \geq\!  3$; we say that $\sigma$ is a \textit{leaf} if $\mathtt{deg} (\sigma) \! =\!  1$ and we say that $\sigma$ is \textit{simple} if $\mathtt{deg} (\sigma) \! = \! 2$. 
We shall use the following notation for the set of leaves of $T$
\begin{equation*}
\mathtt{Lf} (T):= \big\{ \sigma \! \in \!T \backslash \{ r \} : \mathtt{deg} (\sigma) \! =\! 1 \big\} \; . 
\end{equation*}
\begin{definition} Let $\mu$ be a finite Borel measure on $T$; then $(T,d,r, \mu)$ is a \textit{continuum real tree} in the sense of Aldous \cite{Al93}
if $T$ is compact and if  
\begin{equation} 
\label{CRTdef}
\textrm{$\mu$ is diffuse, the topological support of $\mu$ is $T$ and $\mu \big(T\backslash \mathtt{Lf} (T) \big)\! = \! 0$ .} 
\end{equation}
\end{definition}

\vspace{-3mm}

\begin{lemma}
\label{critCRT} Let $\zeta\ino (0, \infty)$ and $d\ino \mathbf{MT} ([0, \zeta])$. To simplify we denote by $(T,d , r, \mu)$ the pointed and measured compact real tree induced by the pseudo-metric $d$ and we denote by $p\! : \! [0, \zeta] \rightarrow T$ the canonical projection. We set 
\begin{equation}
\label{Udefff}
U\! =\! \big \{ s\ino (0, \zeta) : \; \forall s^\prime \ino [0, \zeta] \backslash \{ s \} , \; d(0, s)+ d(s, s^\prime) \! > \! 
d(0, s^\prime) \big\}\; .
\end{equation}
We assume that $[0, \zeta] \backslash U$ is Lebesgue-negligible. Then, $T$ is a continuum real tree as defined in (\ref{CRTdef}).  
\end{lemma}
\noi
\textbf{Proof.} By construction and since $p$ is continuous, the topological support of $\mu$ is $T$. 
We next show that $p(U)\! \subset \! \mathtt{Lf} (T)$. Let $s\ino (0, \zeta)$ be such that $p(s)\! \neq \! r$ and such that $\mathtt{deg} (p(s) ) \! \geq \! 2$. Let $s^\prime$ be such that $p(s^\prime)$ belongs to a connected component of $T\backslash \{ p(s) \}$ that does not contain $r $. Then, 
$p(s)\!  \in \, \rgeo r , p(s^\prime) \lgeo \, $, which implies that 
$d(p(s^\prime), r)\! = \! d( r, p(s))+ d(p(s), p(s^\prime))$, namely, $d(0, s^\prime) \!  = \! d(0, s)+ d(s, s^\prime)$. Thus, $s\! \notin \! U$. This proves that $p(U)\! \subset \! \mathtt{Lf} (T)$. 
Since $[0, \zeta] \backslash U$ is Lebesgue negligible, $T\backslash \mathtt{Lf} (T)$ is $\mu$-negligible.

It remains to prove that $\mu$ is diffuse. Let $s\ino [0, \zeta]$. If $s\! \notin \! U$, then $p(s)\ino T\backslash \mathtt{Lf} (T)$ that is $\mu$-negligible; thus $\mu (\{ p(s) \})\! = \! 0$. 
Assume now that $s\ino U$ and suppose that $s^\prime $ is such that $d(s,s^\prime)\! = \! 0$. The triangle inequality for $d$ implies that $d(0, s)\! = \! d(0, s^\prime)$. Since $s\ino U$, we get $s\! = \! s^\prime$. Thus $p^{-1} (\{ p(s)\} ) \! =\!  \{ s \}$ and $\mu (\{ p(s) \})\! = \! 0$, which completes the proof of the lemma. \cqfd 

\begin{example}
\label{rempli} (\textit{Real trees spanned by graph-trees}) Let $\mathtt{T}$ be a discrete graph-tree with a special vertex $\rho$ viewed as a root. We denote by $d_{\mathtt{gr}}$ the graph-distance on $\mathtt{T}$. The real tree spanned by $\mathtt{T}$ is obtained by joining neighbouring vertices of $\mathtt{T}$ by a unit-length segment of the real line with its own metric. Formally, we can take $\widetilde{\mathtt{T}}\! = \! \{ (\rho, 1) \} \cup \bigcup_{v\in \mathtt{T} \backslash \{ \rho\}} \{ v\} \! \times \! (0, 1]$ and for all $\sigma\! = \! (v, s)$ and $\sigma^\prime \! = \! (v^\prime , s^\prime)$ in $\widetilde{\mathtt{T}}$, we set:
 \begin{equation*} 
\widetilde{d}_{\mathtt{gr}}^{\, \mathtt{T}} (\sigma, \sigma^\prime)= \left\{ 
\begin{array}{ll}
d_{\mathtt{gr}} (v, v^\prime) +s+s^\prime -2 &  \textrm{if $v\! \wedge \!v^\prime \!  \notin \! \{ v, v^\prime\} $,} \\
d_{\mathtt{gr}} (v, v^\prime) +s-s^\prime  &  \textrm{if $v^\prime \! = \! v\! \wedge \!v^\prime$ and $v\! \neq \! v^\prime$, } \\
|s-s^\prime|   &  \textrm{if $v\! = \!v^\prime$, }
\end{array} \right.
\end{equation*}
where $v\! \wedge \! v^\prime$ is the most recent common ancestor of $v$ and $v^\prime$ when we view $\mathtt{T}$ as a family tree whose ancestor is $\rho$. Clearly, $(\widetilde{\mathtt{T}}, \widetilde{d}_{\mathtt{gr}}^{\, \mathtt{T}})$ is a real tree and if we identify $\mathtt{T}$ with 
$\mathtt{T} \! \times \! \{ 1\} \! \subset \! \widetilde{\mathtt{T}}$, we easily check that $\widetilde{d}_{\mathtt{gr}}^{\, \mathtt{T}}$ extends $d_{\mathtt{gr}}$. 
Note that $(\widetilde{\mathtt{T}}, \widetilde{d}_{\mathtt{gr}}^{\, \mathtt{T}})$ is compact {if and only if} $\mathtt{T}$ is finite. \cq 
\end{example}
\begin{example}
\label{rltrcdng} (\textit{Real trees coded by continuous functions}) Let $\zeta \ino (0, \infty)$ and 
let $h\! : \! [0, \zeta] \! \rightarrow \! \bbR$ be continuous. For all $s_1, s_2\ino [0, \zeta]$, we set 
\begin{equation} 
\label{djhczkk}
m_h (s_1, s_2) = \!\! \!\!\!\!\!\! \!\!\!\! \!\!\!\!\inf_{\,  \qquad s\in [s_1 \wedge s_2, s_1 \vee s_2] }\!\!\!\!\!\!\!\!\!\!\!\! \!\!\!\! h (s) \qquad  \textrm{and} \quad d_h (s_1, s_2)= h(s_1)+ h(s_2) - 2 m_h (s_1, s_2) \; .
\end{equation}
We easily check that $d_h \ino \MMT ([0, \zeta])$ and to simplify we denote by $(T_h, d_h, r_h, \mu_h)$ 
the induced  metric space {$(E_{d_h}, d_h, r_{d_h}, \mu_{d_h})$} as defined in (\ref{quotient}) 
and (\ref{induleb}); 
we also denote by $p_h \! : \! [0, \zeta] \! \rightarrow \! T_h$ the canonical projection (instead of $\mathtt{proj}_{d_h}$); $(T_h, d_h, r_h, \mu_h)$ is a pointed measured compact real tree that is the \textit{tree coded by the function $h$}.  

Let $h^\prime \! : \! [0, \zeta] \! \rightarrow \! \bbR$ be continuous. Then observe that 
\begin{equation}
\label{contitreed}
 \forall s_1, s_2, \ino [0, \zeta], \quad \big| d_h (s_1, s_2)\! -\! d_{h^\prime} (s_1, s_2)\big| \leq  4 \sup_{s\in [0, \zeta] } |h(s)\! -\! h^\prime (s)| \; , 
\end{equation}
which shows the continuity in $ \MMT ([0, \zeta])$ of the application $h\mapsto d_h$. \cq 
\end{example}

\begin{remark}
\label{pastt} Let $(T_h, d_h, r_h)$ be the real tree coded by the continuous function 
$h\! : \! [0, \zeta] \! \rightarrow \! \bbR$ and the canonical projection $p_h$ as explained above. Suppose that $h$ is nonnegative and that $h(0)\! = \! 0$. Then $d_h(0, s)\! = \! h(s)$ and for all $s_1, s_2 \ino [0, \zeta]$, we easily get $m_h (s_1, s_2)\! = \! d_h \big( r_h  ,  p_h (s_1) \! \wedge \! p_h (s_2) \big)$.  

  However, let $d\ino \mathbf{MT} ([0, \zeta])$ and set $h(s)\! = \! d(0, s)$, $s\ino [0, \zeta]$ 
 that is a continuous function. The real tree $T\! = \! [0, \zeta] /\!  \sim_d $ is in general different from the the real tree $T_h$ coded by $h$: namely $d\! \neq \! d_h$. But let us mention that all compact real trees can be coded by (many) continuous functions (see D.~\cite{Duq} for more details on the coding of real trees). \cq        
\end{remark}
\noi
\textbf{Gromov-Hausdorff-Prokhorov metric.} Let $(E_1, d_1, r_1, \mu_1)$ and $(E_2, d_2, r_2, \mu_2)$ be two pointed measured compact metric spaces: here $\mu_1$ and $\mu_2$ are finite measures on the respective Borel sigma-fields of $E_1$ and $E_2$, and $r_1\ino E_1$ and $r_2 \ino E_2$ are distinguished points. The pointed \textit{Gromov-Hausdorff-Prokhorov distance} (the \textit{GHP-distance} for short) between $E_1$ and $E_2$ is defined by
\begin{multline}
\label{defGHP}
\bdelta_{\mathtt{GHP}} (E_1, E_2)\! = \! \inf \!  \Big\{ d_{E}^{\textrm{Haus}} \! \big(\phi_1 (E_1), \phi_2 (E_2)\big)  \\ 
+ 
d_E (\phi_1 (r_1), \phi_2 (r_2))  + d_{E}^{\textrm{Prok}} \big(\mu_1 \! \circ \! \phi^{-1}_1\! \! , \mu_2 \! \circ \! \phi^{-1}_2\big)  \Big\}. 
\end{multline}
Here, the infimum is taken over all Polish spaces $(E, d_E)$ and all isometric embeddings $\phi_i: E_i
\hookrightarrow E$, $i\ino \{ 1, 2\}$; 
$d_{E}^{\textrm{Haus}}$ stands for the Hausdorff distance on the space of compact subsets of $E$ (namely, 
$d_{E}^{\textrm{Haus}}(K_1, K_2)\! = \! \inf \{ \epp \ino (0, \infty): K_1 \! \subset \!  K_2^{(\epp)} \; \textrm{and} \;   
K_2 \! \subset \!  K_1^{(\epp)} \}$, where $A^{(\epp)}\! = \! \{ y\ino E: d_E(y, A) \! \leq \! \epp \}$ for all non-empty 
$A\! \subset \! E$); 
$d_{E}^{\textrm{Prok}}$ stands for the Prokhorov distance on the space of finite Borel measures on $E$ (namely, 
$d_{E}^{\textrm{Prok}}(\mu, \nu)\! = \! \inf \{ \epp \ino (0, \infty) \! : \! \mu (K)\! \leq \! \nu(K^{(\epp)}) + \epp  \; \textrm{and} \;   
\nu (K)\! \leq \! \mu(K^{(\epp)}) + \epp , \; \forall K \! \subset \! E \; \textrm{compact} \}$); for all 
$i\ino \{ 1, 2\}$, $\mu_i \! \circ \! \phi^{-1}_i$ stands for the pushforward measure of $\mu_i$ via $\phi_i$. 

\begin{remark}
\label{Kremsmunster} Let $(E,d)$ be a Polish space and let $a,b, c\ino (0, \infty)$. 
We denote by $d_{\mathtt{Prok}}$ (resp.~$d^a_{\mathtt{Prok}}$) 
the Prokhorov distance on the space of finite Borel measures on $(E,d)$ (resp.~$(E, ad)$) and we denote by $d_{\mathtt{Haus}}$ (resp.~$d^a_{\mathtt{Haus}}$) 
the Hausdorff distance on the space of the compact subsets of $(E,d)$ (resp.~$(E, ad)$). 
First note that $d^a_{\mathtt{Haus}}\! = \! ad_{\mathtt{Haus}}$. 
Then, observe that 
$d_{\mathtt{Prok}} (b\mu, c\mu) \! = \! |b\! -\! c| \mu (E)$ and that 
$ (a\wedge b) d_{\mathtt{Prok}} (\mu, \nu)  \! \leq \! d^a_{\mathtt{Prok}} (b\mu, b\nu)  \! \leq \! (a\vee b)
d_{\mathtt{Prok}} (\mu, \nu) $. 

Let $(E_1, d_1, r_1, \mu_1)$ and $(E_2, d_2, r_2, \mu_2)$ be two pointed measured compact metric spaces. We set $(E^\prime_i, d^\prime_i, r^\prime_i, \mu^\prime_i)\! = \! (E_i, ad_i, r_i, b\mu_i)$, $i\ino \{ 1, 2\}$. Then, it is easy to check that $(a\wedge b) \bdelta_{\mathtt{GHP}} (E_1, E_2)\! \leq \!  \bdelta_{\mathtt{GHP}} (E^\prime_1, E^\prime_2)\! \leq \!  (a\vee b) \bdelta_{\mathtt{GHP}} (E_1, E_2)$. \cq 
\end{remark}

\begin{example}
\label{rempli2} Let $\mathtt{T}$ be a finite graph-tree that is equipped with its graph distance and with its counting measure: $\mathtt{m} \! = \! \sum_{v\in \mathtt{T}} \delta_v$. Let $(\widetilde{\mathtt{T}}, \widetilde{d}_{\mathtt{gr}}^{\, \mathtt{T}})$ be the compact real tree spanned by $\mathtt{T}$ (see Example \ref{rempli}). We equip $\widetilde{\mathtt{T}}$ with its length measure $\mathtt{Length}$. 
Up to obvious identifications, we can assume that $\mathtt{T}\! \subset \! \widetilde{\mathtt{T}}$. 
Then, we easily get $d_{\mathtt{Haus}} (\mathtt{T}, \widetilde{\mathtt{T}}) \! \leq \! 1$, $d_{\mathtt{Prok}} (\mathtt{m}, \mathtt{Length})\!  \leq \! 2$ and thus, $ \bdelta_{\mathtt{GHP}} (\mathtt{T}, \widetilde{\mathtt{T}}) 
\! \leq \! 3$. \cq 
\end{example}  

We next recall from Theorem 2.5 in Abraham, Delmas and Hoscheit \cite{AbDeHo13}
the following assertions: $\bdelta_{\mathtt{GHP}} $ is a pseudo-metric (i.e.~it is symmetric and it satisfies the triangle inequality) and $\bdelta_{\mathtt{GHP}} (E_1, E_2)\! = \! 0$ {if and only if} $E_1$ and $E_2$ are \textit{isometric}, namely {if and only if} there exists a bijective isometry $\phi \! : \! E_1 \! \rightarrow \! E_2$ such that $\phi(r_1)\! = \! r_2$ and such that $\mu_2 \! = \! \mu_1 \circ  \phi^{-1}$. Denote by $\bbM$ the isometry classes of pointed measured compact metric spaces. Then, Theorem 2.5 in Abraham, Delmas and Hoscheit \cite{AbDeHo13} asserts {that}
\begin{equation}
\label{rappGHP} \textrm{$(\bbM, \bdelta_{\mathtt{GHP}})$ is a complete and separable metric space. }
\end{equation}
\begin{proposition} 
\label{groweak} Let $d, d^\prime\ino \MMM([0, \zeta])$ and let $E_d$ and $E_{d^\prime}$ be the induced pointed measured compact metric spaces as defined by (\ref{quotient}) and (\ref{induleb}). Then
\begin{equation}
\label{ineqw}
\bdelta_{\mathtt{GHP}} (E_{d},E_{d^\prime}) \leq \tfrac{3}{2} \lVert d \! -\! d^\prime \rVert \; .
\end{equation}
\end{proposition}
\noi
\textbf{Proof.} We use the notation from (\ref{quotient}) and we set 
$\cR \! = \! \{ (\mathtt{proj}_d (s), \mathtt{proj}_{d^\prime} (s)); s\ino [0, \zeta] \}$ that is a relation on $E_{d}\! \times E_{d^\prime}$ since for all $x\ino E_d$ and all $x^\prime \ino E_{d^\prime}$, $\cR \cap (\{ x\} \times E_{d^\prime})$ and $\cR \cap ( E_{d} \times  \{ x^\prime \} )$ are not empty. By the triangle inequality for pseudo-metrics, we get 
\begin{equation}
\label{disineq}
\mathtt{dis} (\cR):= \sup \big\{  \big| d(x,y)\! -\! d^\prime (x^\prime, y^\prime)  \big| \, ; \; (x,x^\prime), (y, y^\prime) \ino \cR  \big\} \leq  \lVert d \! -\! d^\prime \rVert \; .
\end{equation}
We next set $G\! = \! E_d \sqcup E_{d^\prime}$ that is a disjoint union, and we define a metric $\delta$ on $G$ as follows: $\delta$ coincides with $d$  on $E_{d} \! \times \! E_{d}$ and with $d^\prime$ on $E_{d^\prime} \! \times \! E_{d^\prime}$ and for all $x\ino E_d$ and all $x^\prime \ino E_{d^\prime}$, we set 
$$\delta  (x,x^\prime)  = \inf \big\{ d(x,y) + \tfrac{1}{2} \mathtt{dis} (\cR) + d^\prime \! (y^\prime, x^\prime)\, ; \; (y,y^\prime) \ino \cR \big\} \; .$$ 
We easily check that $\delta$ is a separable compact metric on $G$. Thus, $(G, \delta)$ is Polish. The inclusions $E_d\hookrightarrow G$ and  $E_{d^\prime}\hookrightarrow G$ are isometries and since 
$(r_d, r_{d^\prime})\! = \! (\mathtt{proj}_d (0), \mathtt{proj}_{d^\prime} (0))\ino \cR$, we get 
\begin{equation}
\label{hausdfdf}
\delta_{\mathtt{Haus}} (E_d, E_{d^\prime}) \vee \delta (r_d , r_{d^\prime}) \leq \tfrac{1}{2} \mathtt{dis} (\cR) \; ,
\end{equation}
where $\delta_{\mathtt{Haus}} $ stands for the Hausdorff distance of the space of compacts subsets of $G$. 

  Let $K$ be a closed (and thus compact) subset of $G$. We set $Q\! = \! E_{d}\cap K$ and $Q^\prime\! = \! E_{d^\prime}\cap K$ that are also compact subsets of $G$. 
 Note that $\mu_d (K)= \mu_d (Q)$ and $\mu_{d^\prime} (K)\! = \! \mu_{d^\prime} (Q^\prime)$. We set 
$ C \! =\!  \mathtt{proj}_d^{-1} (Q)$, $Q^\prime_0 \!= \! \mathtt{proj}_{d^\prime} (C)$ and  $C^\prime \! =\!  \mathtt{proj}_{d^\prime}^{-1} (Q^\prime_0)$. Since $\mathtt{proj}_d$ and 
$\mathtt{proj}_{d^\prime}$ are continuous, $C$ and $C^\prime$ are compact subsets of $[0, \zeta]$ and $Q^\prime_0 $ is a compact  subset of $E_{d^\prime}$. Denote by $\ell$ the Lebesgue measure on $[0, \infty)$. Observe that $C \! \subset \! C^\prime$, which implies that $\mu_d (K) \! = \! \mu_d (Q)\! = \! \ell (C) \! \leq  \! \ell (C^\prime)\! = \! \mu_{d^\prime} (Q^\prime_0 )$. 
Let $x^\prime \ino Q_0^\prime$. There exists $s\ino C$ such that $\mathtt{proj}_{d^\prime} (s)\! = \! x^\prime$. But note that $x\! :=\! \mathtt{proj}_d(s) \ino Q$ and thus $(x,x^\prime)\ino \cR$. Set $\eta\! := \!  \frac{_1}{^2} \mathtt{dis} (\cR)$; then we have proved that $\mu_d (K)\! \leq \! \mu_{d^\prime} (K^\eta)$ where 
$K^\eta\! = \! \{ z\ino  G\! : \!  \delta (z, K) \! \leq \! \eta \}$. Similarly, we prove that $ \mu_{d^\prime} (K) \! \leq \! \mu_{d} (K^\eta)$ which implies that $\delta_{\mathtt{Prok}} (\mu_{d}, \mu_{d^\prime})\! \leq \! \eta$, where $\delta_{\mathtt{Prok}}$ stands for the Prohorov distance on the space of finite Borel measures on $G$. This inequality 
combined with (\ref{hausdfdf}) implies that $\bdelta_{\mathtt{GHP}} (E_{d},E_{d^\prime})\!  \leq \! \frac{_3}{^2} \mathtt{dis} (\cR)$, which entails (\ref{ineqw}) by (\ref{disineq}).  \cqfd 
{
\subsection{Scaling limit of the range of biased RWs {on} trees.}
\label{singRW}
Let $T$ be an infinite rooted ordered tree. We fix $\lambda \ino (1, \infty)$ and we denote by $(Y_n)_{n\in \bbN}$ the $\lambda$-biased RW {on} $T$ whose transition probabilities are given by (\ref{homsicdef}). We make $T$ a real tree by joining neighbouring vertices by a line isometric to $[0, 1]$ as explained in Example \ref{rempli} and 
we keep denoting {by} $(T, d_{\mathtt{gr}}, \varnothing)$ the resulting rooted real tree. We also denote simply by $(Y_s)_{s\in [0, \infty)}$ the continuous interpolation 
of $Y$: namely, for all $n\ino \bbN$ and all $s\ino [n, n+1]$, $Y_s$ is the unique point of the line $\lgeo Y_n , Y_{n+1} \rgeo$ in the (spanned) real tree $T$ such that $d_{\mathtt{gr}} (Y_n, Y_s)\! = \! s\! -\! n$. 

 Let $\zeta\ino (0, \infty)$. For all $n\ino \bbN$, we then set 
$$\cR_n \! = \! \{ Y_s \, ; \, s\ino [0, \zeta n] \}$$ 
that is a random compact real tree. We equip $\cR_n$ with the occupation measure $m_{\mathtt{occ}}^{_{(n)}}$ induced by the RW, namely $\int_{\cR_n}\!  f(\sigma)  m_{\mathtt{occ}}^{_{(n)}}  (d\sigma) \! = \! \int_0^\zeta f(Y_{sn}) \, ds $.

\begin{proposition}
\label{cvrgrwT} We keep the above notation. Let $(a_n)_{n\in \bbN}$ be a sequence of positive real numbers. We assume that $(Y_n)_{n\in \bbN}$ is recurrent, that $\lim_{n\rightarrow \infty} a_{^n}^{_{-1}}\log n   \! = \! 0$ and that there exists a continuous random process $(H_s)_{s\in [0,\zeta]}$ such that 
\begin{equation}
\label{hypdistoro}
\Big(\tfrac{1}{a_n} d_{\mathtt{gr}} (\varnothing , Y_{ns}) \Big)_{\! s\in [0, \zeta]} \overset{\textrm{(law)}}{\underset{n\rightarrow \infty}{-\!\!\! -\!\!\!-\!\!\! -\!\!\!  -\!\! \! \longrightarrow}} (H_s)_{s\in [0, \zeta]}
\end{equation}
weakly on $\bC ([0, \zeta], \bbR)$. Then, jointly with (\ref{hypdistoro}), the following convergence 
\begin{equation*}
\Big( \cR_n \, , \, \tfrac{1}{a_n} d_{\mathtt{gr}}\,  , \, \varnothing \, , \,   m_{\mathtt{occ}}^{_{(n)}} \Big) \overset{\textrm{(law)}}{\underset{n\rightarrow \infty}{-\!\!\! -\!\!\!-\!\!\! -\!\!\!  -\!\! \! \longrightarrow}} \; 
 \big( T_H, d_H, r_H, \mu_H \big)
\end{equation*}
holds weakly on the space of rooted measured compact metric spaces $\bbM$ equipped with the Gromov-Hausdorff-Prokhorov distance. Here, $ \big( T_H, d_H, r_H, \mu_H \big)$ stands for the real tree coded by $H$ as in Example \ref{rltrcdng}. 
\end{proposition}
\noi
\textbf{Proof.} For all $s\ino [0, \zeta]$, we set 
$H^{_{(n)}}_s\! = \! \tfrac{1}{a_n} d_{\mathtt{gr}} (\varnothing , Y_{ns})$ and for all $s, s^\prime \ino [0, \zeta]$, $d_n (s,s^\prime)\! = \!  H^{_{(n)}}_s+H^{_{(n)}}_{s^\prime}\! -\! 2\min_{r\in [s\wedge s^\prime, s\vee s^\prime]} H^{_{(n)}}_r$, that is the tree 
pseudo-distance coded by $H^{_{(n)}}_{\!^{\!}}$ as in Example \ref{rltrcdng}. Clearly, by (\ref{contitreed}), (\ref{hypdistoro}) implies that $d_n \! \rightarrow \!  d_H$ weakly on $\bC ([0, \zeta]^2, \bbR)$. Then we set $d^*_n (s, s^\prime) \! = \! a_n^{-1} d_{\mathtt{gr}} (Y_{ns}, Y_{n s^\prime})$. Thus $d^*_n (s, s^\prime) \! \leq \! d_n (s, s^\prime)$ and  
by (\ref{distesti}) in Lemma \ref{distcontr}, we get for all $\epp \! >\! 0$ {that}
\begin{eqnarray*}
\bP \! \big( d_n (s, s^\prime) \! -\! d_n^* (s, s^\prime) \! >\! \epp \big) &\!\!\!\!\! \leq  &\!\!\!\! \bP \Big( |Y_{\! \lfloor ns \rfloor }|+|Y_{\! \lfloor ns^\prime \rfloor }| \! -\! 2\!\!\!\!\! \!\!  \!\!\!\!\!  
\min_{\quad \, \lfloor ns\rfloor \leq k \leq \lfloor ns^\prime \rfloor }\!  \!\! \!\!\!\!\!   \!\!\! \!\! |Y_{\! k}| \, \geq   \epp a_n \! -\! 2  + d_{\mathtt{gr}} \big(Y_{\! \lfloor ns \rfloor},Y_{\! \lfloor ns^\prime \rfloor} \big) \Big) \\ 
&\!\!\!\!\! \leq  &\!\!\!\!  (\lfloor ns^\prime \rfloor\! -\! \lfloor ns \rfloor)\frac{ \lambda \! -\! 1 }{\lambda^{\frac{_1}{^2}\epp a_n-1} \! -\! 1}\; \,  \underset{n\rightarrow \infty}{ -\!\!\!-\!\!\! -\!\!\!  -\!\! \! \longrightarrow} 0. 
\end{eqnarray*}
By Proposition \ref{submetric}, $(d^*_n, d_n) \! \rightarrow \! (d_H, d_H)$ weakly on $\bC([0, \zeta]^2, \bbR)^2$. This implies the desired result by Proposition \ref{groweak} and because the compact measured rooted tree induced by the pseudo-metric $d_n^*$ is isometric to $( \cR_n  ,  \frac{_1}{^{a_n}} d_{\mathtt{gr}} ,  \varnothing  ,    m_{\mathtt{occ}}^{_{(n)}} )$. \cqfd

\medskip

We then derive from the previous proposition (the easy part) and from a result due to A.~Dembo and N.~Sun \cite{DeSu12} (the difficult part) 
that {the scaling limit of} the range of critical biased RWs on a supercritical multi-type GW-tree is the Brownian tree. More precisely, we consider a $N$-type GW-tree; for all $k\ino \{ 1, \ldots, N\}$ we denote by $\bmu (k, \cdot)$ a probability {measure} on $\bbW_{\! N} \! = \! \bigcup_{n\in \bbN} \{ 1, \ldots, N\}^n$, the set of finite words written in the alphabet $\{ 1, \ldots, N \}$. For all $p\ino (0, \infty)$, the $p$-th moment matrix $M_p \! = \! (m_p (k, \ell))_{1\leq k, \ell \leq N}$ is given by  
$$ m_p (k, \ell)= \sum_{n\in \bbN} \!\!  \!\!  \!\!  \!\!  \!\!  \!\!  \sum_{\qquad \; 1\leq w_1, \ldots, w_n \leq N} \!\!  \!\!  \!\!  \!\!  \!\!  \!\!  \!\!  \bmu \big(k, (w_1, \ldots, w_n)\big)  \Big( \# \big\{j\ino \{ 1, \ldots, n \} : w_j\! = \! \ell \big\}\Big)^p . 
$$
We assume the following. 

\smallskip

\begin{compactenum}
\item[$(a)$] There exists $p\ino (4, \infty)$ such that for all $k, \ell \ino \{ 1, \ldots, N\}$, $m_p (k, \ell) \! < \! \infty $.

\smallskip

\item[$(b)$] There exists an integer $n_0\! \geq \! 1$, such that $M^{n_0}_1$ has only strictly positive entries.  

\smallskip

\item[$(c)$] Let $\lambda$ be the Perron-Frobenius eigenvalue of $M_1$. We assume that $\lambda \! >\! 1$.  
\end{compactenum}

\smallskip
 
\noi 
Let $\bT\!= \! (T, \varnothing, (a_x)_{x\in T})$ be a $N$-type GW$(\bmu)$-tree. Namely, $T$ is a random rooted ordered tree, $a_x\ino \{ 1, \ldots, N\}$ is the type of vertex $x\ino T$ and $\bT$ satisfies the following property: recall that $k_\varnothing (T)$ is the number of children of the root and recall that 
{for all $j\ino \{ 1, \ldots, k_\varnothing(T)\}$}, $\theta_{(j)} \bT $ stands for the tree stemming form the $j$-th child of {$\varnothing$} (equipped with the types of corresponding vertices). Then, conditionally given $a_\varnothing$, 
the types $\big( a_{(1)}, \ldots, a_{(k_\varnothing (T))} \big)$ of the children of the root have conditional law $\bmu (a_\varnothing , \cdot)$ and conditionally given 
$\big( a_{(1)}, \ldots, a_{(k_\varnothing (T))} \big)$, the subtrees 
$\theta_{(1)} \bT , \ldots , \theta_{(k_\varnothing (T))} \bT$ are independent GW$(\bmu)$-trees. 

Then, conditionally given $T$, let $(Y_n)_{n\in \bbN}$ be the $\lambda$-biased RW on $T$ started at $\varnothing$ and denote by $(Y_s)_{s\in [0, \infty)}$ its continuous interpolation as explained above. 
Under {Assumptions} $(a)$, $(b)$ and $(c)$ on $\bmu$, Theorem 1.1 in A.~Dembo and N.~Sun \cite{DeSu12} p.~3 asserts that there exists 
$\sigma \ino (0, \infty)$, a constant that depends only on $\bmu$, such that a.s.~on the {non-extinction event 
$\{|T|\! = \! \infty \}$}, 
\begin{equation}
\label{DemboSun}
\Big(\tfrac{1}{\sigma \sqrt{n}} d_{\mathtt{gr}} (\varnothing , Y_{ns}) \Big)_{\! s\in [0, \infty)} \overset{\textrm{(law)}}{\underset{n\rightarrow \infty}{-\!\!\! -\!\!\!-\!\!\! -\!\!\!  -\!\! \! \longrightarrow}} \big( |B_s| \big)_{s\in [0, \infty)}
\end{equation}
where $(B_s)_{s\in [0, \infty)}$ is a standard linear Brownian motion such that $B_0\! = \! 0$ a.s. By {Proposition \ref{cvrgrwT}}, we immediately get the following {new result}. 
\begin{corollary}
\label{corodiff} We keep the notation as above. Let $\zeta \ino (0, \infty)$. Let $\cR_n\! = \! \{ Y_{s} ; s\ino [0, n\zeta] \}$ be the range of $Y$ up to time $n\zeta$ and let $m_{\mathtt{occ}}^{_{(n)}}$ be the occupation measure of $Y$ on $\cR_n$ as introduced before. 
Under Assumptions $(a)$, $(b)$ and $(c)$ on $\bmu$, a.s.~on the event 
{$\{|T|\! = \! \infty\}$}, we get 
\begin{equation}
\label{Ariach}
\Big( \cR_n \, , \, \tfrac{1}{\sigma \sqrt{n}} d_{\mathtt{gr}}\,  , \, \varnothing \, , \,   m_{\mathtt{occ}}^{_{(n)}} \Big) \overset{\textrm{(law)}}{\underset{n\rightarrow \infty}{-\!\!\! -\!\!\!-\!\!\! -\!\!\!  -\!\! \! \longrightarrow}} \; 
 \big( T_{|B|}, d_{|B|}, r_{|B|}, \mu_{|B|} \big)
\end{equation}
weakly on the space of rooted measured compact metric spaces $\bbM$ equipped with the Gromov-Hausdorff-Prokhorov distance. Here, $ \big(  T_{|B|}, d_{|B|}, r_{|B|}, \mu_{|B|} \big)$ stands for the real tree coded by the reflected Brownian motion $( |B_s| )_{s\in [0, \zeta]}$. 
\end{corollary}

In literature,
the first scaling limit for the range of tree-valued RWs appears in D.~\cite{Du05}: in this paper the tree is {$\mathtt{b}$-ary} and the RW is slightly super-critical (see Theorem 2.1 \cite{Du05} p.~2224; {Lemma 3.7 \cite{Du05} p.~2241} also contains a local law of large numbers for the range). 
When $T$ is a supercritical {single-type} GW-tree, Y.~Peres and {O.~Zeitouni} \cite{PerZei08} have first proved 
 (\ref{DemboSun}) when the offspring distribution has 
exponential moments (see Theorem 1 \cite{PerZei08}, p.~596). 
Then, E.~Aïd\'ekon and L.~de Raph\'elis in \cite{AidRap17} have proved (\ref{DemboSun}) 
for supercritical {single-type} GW-tree under a second moment assumption and they also proved (\ref{Ariach}) in these cases (Theorem 1.1 \cite{AidRap17}, p.~645). In the same article, they extend (\ref{DemboSun}) and (\ref{Ariach}) to RWs in random environment on GW-trees (see Theorem 6.1 \cite{AidRap17}, p.~660).

\subsection{Snake metrics.}
\label{snkmtrsec}
Snakes are path-valued processes that provide a nice parametrization of the spatial positions of a population whose genealogy is a continuum tree and that are scaling limits of branching random walks. Snake processes, and in particular the Brownian snake, has been introduced by J-F.~Le Gall in \cite{LG93bis} to study fine properties of super-Brownian motion. 
In this section, we first 
recall basic definitions on snakes, in a deterministic setting and in dimension $1$. 
Then, we introduce a pseudo-metric derived from a snake, we study its continuity properties and we show that snake metrics actually yield real-trees. Finally, we prove elementary geometric properties of such real trees that are deterministic version of the cactus introduced by Curien, Le Gall and Miermont in \cite{CuLGMi13}.  
\begin{definition}
\label{rePopol} We fix $\zeta \ino (0, \infty)$ and we denote by $\bC ([0, \zeta], \bbR)$ the space of the continuous functions from $[0, \zeta]$ to $\bbR$; it is a Banach space when equipped with the uniform norm $\lVert \cdot \rVert_\infty$. We also denote by $\bC_0 $ the space of continuous functions from $[0, \infty)$ to $\bbR$ that is a Polish space when equipped with the following metric: 
\begin{equation}
\label{uttcpct}
\forall w, w^\prime \ino \bC_0, \quad \delta_{\mathtt{u}} (w,w^\prime)= \sum_{n\in \bbN} 2^{-n-1} \min\!  \big( 1 \, , \!\!\!\! \sup_{\quad r\in [0, n]}  \!\!\!\! |w(r)\! -\! w^\prime (r)| \, \big) \; .
\end{equation}
\begin{compactenum}
\item[$(a)$] We denote by $\bC ([0, \zeta],  \bC_0)$ the space of the $\delta_{\mathtt{u}}$-continuous functions from $[0, \zeta]$ to $\bC_0$ equipped with the distance 
$\Delta^{\!*} (w, w^\prime)\! = \! \sup_{s\in [0, \zeta] } \delta_{\mathtt{u}} (w_s (\cdot), w_s^\prime (\cdot) )$, for all $w, w^\prime \ino \bC ([0, \zeta], \bC_0)$.   
We next equip the product space $\bC ([0, \zeta], \bbR) \times \bC ([0, \zeta],  \bC_0 )$ with the following distance: for all  $(h,w), (h^\prime, w^\prime) \ino \bC ([0, \zeta], \bbR) \! \times \! \bC ([0, \zeta], \bC_0)$, 
\begin{equation*}
\Delta \big(  (h,w), (h^\prime, w^\prime)\big)= \lVert h\! -\! h^\prime\rVert_\infty + \Delta^{\!*}  (w, w^\prime)\; .
\end{equation*}

\smallskip

\item[$(b)$] We denote by $\fSigma ([0, \zeta])$ the space 
of the $\bbR$-valued \textit{snakes}; namely, the space of $(h,w) \ino \bC ([0, \zeta], \bbR) \! \times \! \bC ([0, \zeta], \bC_0)$ that satisfy $h\! \geq \! 0$, $h(0)\! = \! h(\zeta) \! = \! 0$ and the following. 
\begin{compactenum}

\smallskip

\item[$(b1)$] For all $s\ino [0, \zeta]$ and for all $r\ino [h(s), \infty)$, $w_s (r)\! = \! w_s (h(s))\! =: \! \widehat{w}_s $. 

\smallskip

\item[$(b2)$] For all $s_1, s_2 \ino [0, \zeta]$ and for all $r\ino [0, m_h (s_1, s_2)]$, $w_{s_1} (r)\! = \! w_{s_2} (r)$, where we recall from (\ref{djhczkk}) the definition of $m_h (s_1, s_2)$. 
\end{compactenum}
We refer to $(b2)$ as to the \textit{snake property}. The function $h$ is called the \textit{lifetime process} and the function $(\widehat{w}_s)_{s\in [0, \zeta]}$ is called the \textit{endpoint process} of the snake. \cq 
\end{compactenum}
\end{definition}
We easily check that $\big( \bC ([0, \zeta], \bC_0), \Delta^{\!*} \big)$ and $\big(  \bC ([0, \zeta], \bbR) \! \times \! \bC ([0, \zeta], \bC_0), \Delta \big)$ are Polish spaces and that $\fSigma ([0, \zeta])$ is a $\Delta$-closed subset.  

The following lemma is used in the proof of Theorem \ref{main2}. 
\begin{lemma}
\label{Innsbruck} Let $(h,w)\ino \fSigma ([0, \zeta])$. For all $\eta \ino (0, \infty)$, we set $\omega_\eta (w)\!  = \! 
\sup \big\{ \delta_{\mathtt{u}} (w_s, w_{s^\prime}) \, ; s,s^\prime\ino [0, \zeta] \! : \! |s\! -\! s^\prime| \! \leq \! \eta \big\} $ and $\omega_\eta (\widehat{w})\!  = \! \sup \{ |\widehat{w}_s \! -\! \widehat{w}_{s^\prime}| \, ; s,s^\prime\ino [0, \zeta ] \! : \! |s\! -\! s^\prime| \! \leq \! \eta \}$. Then, $\omega_\eta (w) \! \leq \! 2 \omega_\eta (\widehat{w})$. 
\end{lemma}
\noi
\textbf{Proof.} For all $s\ino [0, \zeta]$ and for all $r\ino [0, h(s)]$ set $\alpha_{s,r}\! = \! \sup \{ s^\prime\ino [ 0, s] \! : \! h(s^\prime) \! \leq \! r \}$ and $\beta_{s, r}\! = \!\inf \{ s^\prime\ino [ s, \zeta] \! : \! h(s^\prime) \! \leq \! r \}$. Then, Definition \ref{rePopol} $(b2)$ implies that $w_s(r)\! = \! w_{s^\prime} (r)$ for all $s^\prime \ino [\alpha_{s, r}, \beta_{s, r}]$ and in particular, $ \widehat{w}_{\alpha_{s, r}}\! = \!  \widehat{w}_{\beta_{s, r}}\! = \! w_s(r)$ since $r\! = \! h( \alpha_{s, r})\! = \!  
h( \beta_{s, r})$.

Next, fix $s, s^\prime \ino  [0, \zeta]$ such that $|s\! -\! s^\prime| \! \leq \! \eta $. 
To simplify we set $m\! = \! m_h (s,s^\prime)$.  
Let $r_* \ino [m, h(s)]$ 
be such that $|w_s (r_*) \! -\!  w_s (m) |\! = \! \max_{r\in [m, h(s) ]} |w_s(r)\! -\!  w_s (m) |$. Suppose that 
$s\! \leq \! s^\prime$ (resp.~that $s^\prime \! \leq \! s$), 
then $\beta_{s, r_*}, \beta_{s, m}\ino [s, s^\prime]$ (resp.~$\alpha_{s, r_*}, \alpha_{s, m}\ino [s^\prime, s]$) 
and $|w_s (r_*) \! -\!  w_s (m) |\! = \! |\widehat{w}_{\beta_{s, r_*}}- \widehat{w}_{\beta_{s, m}}|$ (resp.$= \! |\widehat{w}_{\alpha_{s, r_*}}- \widehat{w}_{\alpha_{s, m}}|$). Thus,  
$ \max_{r\in [m, h(s) ]} |w_s(r)-  w_s (m) | \! \leq \! \omega_\eta (\widehat{w}) $.  Similarly, $ \max_{r\in [m, h(s^\prime) ]} |w_{s^\prime} (r)\! -\!  w_{s^\prime} (m) | \! \leq \! \omega_\eta (\widehat{w}) $. By 
Definition \ref{rePopol} $(b2)$, we get $w_s(m)\! = \! w_{s^\prime}(m)$ and $\sup_{r\in [0, \infty) } |w_s (r)\! -\! w_{s^\prime} (r)| \! \leq  \!  \max_{r\in [m, h(s) ]} |w_s(r)\! -\!  w_s (m) | +\max_{r\in [m, h(s^\prime) ]} |w_{s^\prime} (r)\! -\!  w_{s^\prime} (m) |$, which easily implies the desired result. \cqfd

\begin{definition}
\label{defsnkme} Let $\zeta \ino (0, \infty)$ and let $(h,w)\ino \fSigma ([0, \zeta])$. Recall from (\ref{djhczkk}) the definition of $m_h (\cdot, \cdot)$ and recall from Definition \ref{rePopol} $(b1)$ 
the definition of $\widehat{w}$. For all $s_1, s_2\ino [0, \zeta]$, we set 
\begin{eqnarray}
\label{snkmeff}
M_{h,w} (s_1, s_2) \!\! \! \! \!& =&\!\!  \! \!\!  \min \! \Big( \!  \min \! \big\{  w_{s_1} (r) ; r\ino [m_h (s_1, s_2) , h(s_1)] \big\}  ,  \min \! \big\{ w_{s_2} (r) ; r\ino [m_h (s_1, s_2) , h(s_2)] \big\} \! \Big)  \nonumber \\
 &   &\!\!\!\! \!\!\!\! \!\!\!\! \!\!\!\!  \textrm{and} \quad d_{h,w} (s_1, s_2)\! = \! \widehat{w}_{s_1}+  \widehat{w}_{s_2} -2M_{h,w} (s_1, s_2) \; .
\end{eqnarray}
We call $d_{h,w}$ the \textit{snake metric associated with $(h,w)$} (see Lemma \ref{snkmtrcL} below). \cq
\end{definition}
\begin{lemma} 
\label{Hermagor} Let $\zeta \ino [0, \infty)$ and for all $n\ino \bbN$, let $(h_n , w^{_{(n)}})\ino \fSigma ([0,\zeta])$ that converges to $(h, w)$ in $\bC([0, \zeta], \bbR) \! \times \! \bC([0, \zeta ], \bC_0)$ equipped with $\Delta$ as in Definition \ref{rePopol} $(a)$. Recall that $\lVert \cdot \rVert$ stands for the uniform norm on 
$\bC([0, \zeta]^2, \bbR)$. Then, 
$ \lim_{n\rightarrow \infty}  \lVert d_{h_n} \! -\!  d_{h} \rVert\! = \!  \lim_{n\rightarrow \infty}  \lVert d_{h_n, w^{(n)}} \! -\!  d_{h,w} \rVert = 0 $ (see \eqref{djhczkk} and \eqref{snkmeff}). 
\end{lemma}
\noi
\textbf{Proof.} The first limit follows from (\ref{contitreed}). To prove the second one, we fix $\epp \ino (0, 1)$ and we set $a\! = \! 1+ \sup_{n\in \bbN} \max_{s\in [0, \zeta]} |h_n (s)| $. 
Let $n_1\ino \bbN$ such that $\Delta ((h_n,w^{_{(n)}}), (h,w)) \! \leq \! 2^{-a-2} \epp$, for all $n\! \geq \! n_1$. Thus, for all $n\! \geq \! n_1$ and for all $s\ino [0, \zeta]$,  
\begin{equation}
\label{Melk}
 \lVert w^{{(n)}}_{s} \! -\! w_s \rVert_\infty\! : = \! \sup_{r\in [0, \infty)} 
| w^{{(n)}}_{s} (r) \! -\! w_s (r)| \! = \!  \sup_{r\in [0, a]} 
| w^{{(n)}}_{s} (r) \! -\! w_s (r)| \leq \epp . 
\end{equation}
Fix $s_1,s_2 \ino [0, \zeta]$ and set $m_n \! = \! m_{h_n} (s_1,s_2)$ and $m \! = \! m_{h} (s_1,s_2)$. 
By (\ref{Melk}), we get 
\begin{equation}
\label{Spitz}
\forall n\! \geq \! n_1, \forall j\ino \{ 1, 2\} , \qquad  \Big| \!\! \!\!\!\!\! \!\! \!\! \! \min_{\quad \quad r\in [m_n, h_n(s_j) ]} \!\!\!\!\!\! \!\!\!\!\!\! \! \!  w^{_{(n)}}_{^{s_j}} (r) \, -\! \!  \!\! \!\!\!\!\! \!\! \!\! \!  \min_{\quad \quad r\in [m_n, h_n(s_j) ]} \!\! \!\!\!\!\! \!\! \!\! \! \!\!  w^{_{\! }}_{^{s_j}} (r) \;  \Big| \leq  \lVert w^{{(n)}}_{s_j} \! -\! w_{s_j} \rVert_\infty \leq \epp . 
\end{equation}

 Next, for all uniformly continuous $g\! : \! [0, \infty) \! \rightarrow \! \bbR$, and for all $\eta \ino (0, \infty)$, we use the notation 
$\omega (g, \eta)\! = \! \sup \{ |g(s) \! -\! g(s^\prime)|; s, s^\prime\ino [0, \infty) \! : \! |s\! -\! s^\prime| \! \leq \! \eta  \}$ for the $\eta$-uniform modulus of continuity of $g$. 
We recall that $\lVert h\! -\! h_n\rVert_\infty\! \leq  \! \Delta ((h_n,w^{_{(n)}}), (h,w)) $. 
Observe that for all $j\ino \{ 1, 2\}$ and all $n\ino \bbN $, 
\begin{equation}
\label{Baden}
\Big| \!\! \!\!\!\!\! \!\! \!\! \! \min_{\quad \quad r\in [m_n, h_n(s_j) ]} \!\!\!\!\!\! \!\!\!\!\!\! \! \!  w^{_{\! }}_{^{s_j}} (r) -  \!\! \!\!\!\!\! \!\! \!\! \!  \min_{\quad \quad r\in [m, h(s_j) ]} \!\! \!\!\!\!\! \!\! \!\!   w^{_{\! }}_{^{s_j}} (r)  \Big| \leq  \omega \big( w_{s_j} , \lVert h\! -\! h_n 
  \rVert_\infty \big) \; .
\end{equation}
By the definition (\ref{snkmeff}) of $d_{h_n, w^{(n)}}$, by (\ref{Melk}), (\ref{Spitz}) and (\ref{Baden}), 
for all $n\! \geq \! n_1$, we get the following:  
\begin{equation}
\label{Lavamund}
\big| \, d_{h_n, w^{(n)} } (s_1, s_2) \! -\! d_{h, w} (s_1, s_2)\,  \big| \leq 4\epp + 3 \omega \big( w_{s_1} , \lVert h\! -\! h_n 
  \rVert_\infty \big)  +3 \omega \big( w_{s_2} , \lVert h\! -\! h_n  \rVert_\infty \big)  . 
\end{equation}
Since $w$ is uniformly $\delta_{\mathtt{u}}$-continuous on $[0, \zeta]$, there exist $\sigma_1, \ldots, \sigma_p\in [0, \zeta]$ such that for all $s\ino [0, \zeta]$, there exists $k\ino \{ 1, \ldots, p\}$, such that $\lVert w_s \! -\! w_{\sigma_k} \rVert_\infty \! 
\leq \! \epp$, which implies that $ |\omega \big( w_{s} ,\eta \big)\! -\! \omega \big( w_{\sigma_k} ,\eta \big)| \! \leq \! 2\epp $, for all $\eta\ino (0, \infty)$. By (\ref{Lavamund}), it implies that 
$$\lVert d_{h_n, w^{(n)}} \! -\!  d_{h,w} \rVert \leq 16 \epp + 6 \max_{1\leq k\leq p} \omega \big( w_{\sigma_k} , \lVert h\! -\! h_n \rVert_\infty \big), $$
which implies $\limsup_{n\rightarrow \infty}\lVert d_{h_n, w^{(n)}} \! -\!  d_{h,w} \rVert  \l \leq \! 16 \epp$, for all $\epp \ino (0, \infty)$. It completes the proof of the lemma. \cqfd 
\begin{remark}
 \label{quotsnk} Let $(h,w)\ino \fSigma ([0, \zeta])$. Definition \ref{rePopol} $(b)$ means 
 that $w$ is actually defined on the real tree $(T_h, d_h, r_h, \mu_h)$ coded by 
 $h$ (as defined in Remark \ref{rltrcdng}). 
 Indeed, let $s_1, s_2\ino [0, \zeta]$ 
 be such that $d_h(s_1, s_2)\! = \! 0$; then $h(s_1)\! = \! h(s_2)\! = \! m_h (s_1, s_2)$ and Definition \ref{rePopol} $(b)$ implies that $w_{s_1}\! = \! w_{s_2}$. Up to a slight abuse of notation, 
it therefore makes sense to define $w$ on $T_h$ as follows: for all $\sigma \ino T_h$ and for all $s\ino [0, \zeta]$ such that $\sigma\! = \!  p_h \! (s)$, then 
 \begin{equation}
 \label{snktree}
w_\sigma (\cdot ):= w_s (\cdot)  \quad  \textrm{and} \quad  \widehat{w}_\sigma := \widehat{w}_s   \; , 
\end{equation}
where we recall that $p_h\! : \! [0, \zeta] \! \rightarrow\!  T_h$ stands for the canonical projection. It is easy to check that $ w\! :\!  T_h \rightarrow \! \bC_0$ is continuous. Moreover,  
Definition \ref{rePopol} $(b)$ combined with the argument of the proof of Lemma \ref{Innsbruck} 
{entails} the following.
\begin{equation*}
\forall \sigma \ino T_h, \; \forall \gamma \ino \lgeo r_h, \sigma \rgeo, \quad  w_\sigma \big( d_h (r_h, \gamma)\big)= \widehat{w}_\gamma \; .
\end{equation*}
We also get the following: let $\sigma_1, \sigma_2 \ino T_h$ and $s_1, s_2 \ino [0, \zeta]$ such that $p_h (s_i)\! = \! \sigma_i$, $i\ino \{ 1, 2\}$; then 
\begin{equation}
\label{frounsyy}
 M_{h,w} (s_1, s_2) =\!\!  \!\!\!\!\min_{\quad \gamma \in  \lgeo   \sigma_1 ,  \sigma_2 \rgeo }\!\!\!\!\!  \widehat{w}_\gamma \quad \textrm{and thus} \quad d_{h,w} (s_1, s_2)= \widehat{w}_{\sigma_1} + \widehat{w}_{\sigma_2} -2 \!\!  \!\!\!\!\min_{\quad \gamma \in  \lgeo   \sigma_1 ,  \sigma_2 \rgeo }\!\!\!\! \! \widehat{w}_\gamma \; .
\end{equation}
Up to a slight abuse of notation, it makes sense to view  
$M_{h, w}$ and $d_{h, w}$ as continuous functions from $T_h\! \times \! T_h$ to $\bbR$. \cq 
\end{remark}

\begin{lemma}
\label{snkmtrcL} Let $\zeta \ino (0, \infty)$ and let $(h,w)\ino \fSigma ([0, \zeta])$. Let $d_{h,w}$ be the associated snake metric as in Definition \ref{defsnkme}. 
Then, $d_{h,w} \ino \MMT([0, \zeta])$. Namely, it is a continuous pseudo-metric on $[0, \zeta]$ that satisfies the four points inequality (\ref{fourpoints}). 
\end{lemma}
\noi
\textbf{Proof.} We first prove the continuity of $d_{h,w}$. To that end, for all $a, b\ino [0, \infty)$ and all $s\ino [0, \zeta]$ we set $\phi_{a, b} (s)\! = \! \min \{ w_s (r)\, ; r\ino [a\wedge b, a\vee b ] \}$. First note that $|\phi_{a,b}(s)\! -\! \phi_{a, b} (s^\prime) | \! \leq \! \sup_{r\in [0, a\vee b ] }|w_s (r) \! -\! w_{s^\prime} (r) |$. 
Then, observe that for a fixed $s$, $(a, b) \! \mapsto \phi_{a, b} (s)$ is continuous. This easily implies that $(a,b,s)\! \mapsto \! \phi_{a, b} (s)$ is continuous. Therefore, $d_{h,w}$ is continuous on $[0, \zeta]^2$.

Let $s_1, s_2, s_3, s_4 \ino [0, \zeta]$. We set $X_1\! = \! M_{h,w} (s_1, s_2)+ M_{h,w} (s_3, s_4)$, $X_2\! = \! M_{h,w} (s_1, s_3)+ M_{h,w} (s_2, s_4)$ and 
$X_3\! = \! M_{h,w} (s_1, s_4)+ M_{h,w} (s_2, s_3)$, so that $d_{h,w}$ satisfies the four points inequality (\ref{fourpoints}) if $X_1 \ge \min \big\{X_2, X_3\big\}$, which is a consequence of 
\begin{equation}
\label{freerz}
\#\big\{ i \ino \{ 1, 2, 3 \} : X_i \! =\!  \min (X_1, X_2, X_3) \big\} \geq 2 \; .
\end{equation}
To prove (\ref{freerz}) it is convenient to work on $(T_h, d_h, r_h)$ that is the pointed compact real tree coded by $h$ as explained in {Example} \ref{rltrcdng}: 
recall that $p_h \! : \! [0, \infty] \! \rightarrow \! T_h$ stands for the canonical projection and recall from (\ref{snktree}) the definition of $w_\sigma$ and $\widehat{w}_\sigma$ for all $\sigma \ino T_h$.  For all $i\ino \{ 1, \ldots , 4 \}$, we set $\sigma_i \! = \! p_h (s_i)$ and we recall from (\ref{frounsyy}) that $M_{h,w} (\sigma_i, \sigma_j) \! = \! M_{h, w}(s_i, s_j)$. Recall from (\ref{braachpt}) the definition of branch points in $T_h$.  Since (\ref{freerz}) does not depend on a specific indexation of the $s_i$, without loss of generality we can assume that 
$\gamma\! := \! \mathtt{br} (\sigma_1, \sigma_2, \sigma_3)\!= \!  \mathtt{br} (\sigma_1, \sigma_2, \sigma_4)$, that $\gamma^\prime\! = \! \mathtt{br} (\sigma_3, \sigma_4, \sigma_1)\! = \! 
\mathtt{br} (\sigma_3, \sigma_4, \sigma_2)$ and that 
$$ a:=\!\!\!\!\!\!\!\!  \min_{\quad \sigma \in \lgeo  \gamma ,\,  \sigma_1\rgeo} \!\!\!\!\!\!  \widehat{w}_\sigma \, \leq b:= \!\!\!\!\!\!\!\! \min_{\quad \sigma \in \lgeo  \gamma , \, \sigma_2\rgeo} \!\!\!\!\!\!  \widehat{w}_\sigma\, , \quad d:=\!\!\!\!\!\!\!\! \min_{\quad \sigma \in \lgeo  \gamma^\prime  \!  , \, \sigma_3\rgeo} \!\!\!\!\!\!  \widehat{w}_\sigma  \, \leq e:= \!\!\!\!\!\!\!\! \min_{\quad \sigma \in \lgeo  \gamma^\prime \! , \,\sigma_4\rgeo} \!\!\!\!\!\!  \widehat{w}_\sigma \quad  \textrm{and} \quad a \leq d . $$ 
We also set $c\! : =\! \min_{\sigma \in \lgeo \gamma , \gamma^\prime \rgeo}  \widehat{w}_\sigma $. Then, $X_1\! = \! a+d$, 
$X_2\! = \! (a\wedge c)+ (b\wedge c \wedge e)$ and 
$X_3\! = \! (a\wedge c) + (b\wedge c\wedge d)$. We have four cases to consider;  $(i)$: if $c\! \leq \! a$, then $X_1\! = \! a+d \! \geq \!  X_2\! = \! X_3 \! = \! 2c$; $(ii)$: if 
$a \! \leq \! c \! \leq \! d$, then $X_1\! = \! a+d \! \geq \!  X_2\! = \! X_3 \! = \! a+ (b\wedge c) $; $(iii)$: if $d\! \leq \! c$ and $b \! \leq \! d$, then $X_1\! = \! a+d \! \geq \!  X_2\! = \! X_3 \! = \! a+b$; $(iv)$: if $d\! \leq \! c$ and $d \! \leq \! b$, then $X_1 \! = \! X_3\! = \! a+d \! \leq \! X_2 \! = \! a+ (b \wedge c \wedge e) $. This proves (\ref{freerz}) and it completes the proof of the lemma.\cqfd 

\begin{definition}
\label{notasnkm} Let $\zeta \ino (0, \infty)$ and let $(h,w)\ino \fSigma ([0, \zeta])$. 
Let $d_{h,w}$ be the associated snake metric (as in Definition \ref{defsnkme}). Since $d_{h,w} \ino \MMT([0, \zeta])$ by Lemma \ref{snkmtrcL}, we denote by $T_{h,w}\! = \! [0, \zeta]/\! \! \sim_{d_{h,w}}$ the corresponding real tree and to simplify we denote by $p_{h,w}: [0, \zeta] \! \rightarrow \! T_{h,w}$ the canonical projection, by $r_{h,w}\! = \! p_{h,w} (0)$ the root 
of $T_{h,w}$, and by $\mu_{h,w}$ the measure on $T_{h,w}$ induced by the Lebesgue measure on $[0, \zeta]$ via $p_{h,w}$: namely, 
$$\int_{T_{h,w}} \! f(x)  \, \mu_{h,w} (dx) \! = \!\! \int_0^\zeta \!  f(p_{h,w} (s)) \, ds \; ,$$ 
for all continuous $f\! : \! T_{h,w} \! \rightarrow \! \bbR$.
Since the pseudo-metric $d_{h,w}:[0, \zeta]^2 \! \rightarrow \! [0, \infty)$ is continuous, so is $p_{h,w}$ and 
$(T_{h,w}, d_{h,w}, r_{h,w}, \mu_{h,w})$ is a pointed measured 
compact real tree that we call the \textit{snake tree} associated with $(h,w)$. By Lemma \ref{Hermagor} and Proposition \ref{groweak}, $(h,w)\ino\fSigma ([0, \zeta])\! \mapsto \! T_{h,w}\ino \bbM$ is $(\Delta, \bdelta_{\mathtt{GHP}})$-continuous. \cq  
\end{definition}

\begin{remark}
\label{Villach}
Let $(h, w)\ino \fSigma ([0, \zeta])$ and $ \alpha , a, b\ino (0, \infty)$. 
For all $s\ino [0, b \zeta]$, we set $h^\prime (s)\! = \! \alpha h(s/b)$ and $w^\prime_s (r)\! =\!  aw_{s/b} (r/\alpha )$. Then 
$(h^\prime, w^\prime) \ino \fSigma ([0, b\zeta])$ and thanks to Remark \ref{Passau}, we easily check 
that $(T_{h^\prime}, d_{h^\prime}, r_{h^\prime}, \mu_{h^\prime})$ is isometric to $(T_h, \alpha d_h, r_h, b\mu_h)$ and that  
$(T_{h^\prime, w^\prime }, d_{h^\prime, w^\prime}, r_{h^\prime, w^\prime}, \mu_{h^\prime, w^\prime})$ is isometric to 
$(T_{h, w}, ad_{h, w}, r_{h, w}, b\mu_{h, w})$. \cq  
\end{remark} 

\begin{remark}
\label{SktPaul}
Let $(h, w) \ino \fSigma ([0, \zeta])$; let $(T_{h}, d_h, r_h, \mu_h)$ and 
$(T_{h,w}, d_{h,w},r_{h,w},\mu_{h,w})$ be the compact real trees coded by resp.~$h$ and $(h,w)$ and recall that $p_h$ and $p_{h,w}$ stand for the canonical projections from $[0, \zeta]$ to resp.~$T_h$ and $T_{h,w}$. Observe that it actually makes sense to define a function $\mathtt{y}: T_h \! \rightarrow \!  T_{h,w}$ by setting 
\begin{equation}
\label{Neunkirchen}
\forall s\ino [0, \zeta] , \quad \mathtt{y}(p_{h} (s))= p_{h,w} (s) \; . 
\end{equation}
It is easy to deduce from (\ref{frounsyy}) that 
\begin{equation}
\label{Judenburg}
d_{h,w} \big( \mathtt{y}(\sigma), \mathtt{y}(\sigma^\prime) \big)\! = \! \widehat{w}_{\sigma}+  \widehat{w}_{\sigma^\prime}\! -\! 2 \min_{\gamma \in \lgeo \sigma , \sigma^\prime \rgeo} \widehat{w}_\gamma, 
\end{equation}
It implies that $\mathtt{y}$ is continuous and surjective. Note that $\mu_{h,w}$ is the image of $\mu_h$ via $\mathtt{y}$. \cq
\end{remark}

We next prove two results that deal with basic geometric properties of {snake trees}. The first one provides conditions for a snake tree to be a continuum real tree.  
\begin{lemma}
\label{snkCRT} Let $(h,w)\ino \fSigma([0, \zeta])$. We set 
\begin{equation}
\label{Uhwdef}
U_{h,w}\! = \! \Big\{ s\ino [0, \zeta]: \, h(s) \! >\! 0 \quad \textrm{and} \quad  \forall \epp \ino (0, h(s)), \!\!\!\!\!\! \!\!\!\!\inf_{\qquad r \in [h(s) -\epp , h(s) ]} \!\!\!\!\!\!\!\!\!\!\!\!\!\! \!\!  w_s(r)  \;  < \! \widehat{w}_s \,   \Big\} 
\end{equation}
and we assume that $[0, \zeta ]\backslash U_{h,w}$ is Lebesgue-negligible. We also assume that the tree $(T_h, d_h, r_h, \mu_h)$ coded by $h$ is a continuum real tree as defined in (\ref{CRTdef}). 
Then, the snake tree $(T_{h,w}, d_{h,w}, r_{h,w}, \mu_{h,w})$ is also a continuum real tree. 
\end{lemma}
\noi
\textbf{Proof.} We set $V\! = \! U_{h,w}\cap p_h^{-1} (\mathtt{Lf} (T_h))$. Since $[0, \zeta ]\backslash U_{h,w}$ is Lebesgue-negligible and since $T_h$ is a continuum real tree, $[0, \zeta]\backslash V$ is Lebesgue-negligible. To conclude the proof we are going to show that $V \! \subset \! U$, where $U$ is defined in (\ref{Udefff}) with $d$ replaced by $d_{h,w}$, so that Lemma \ref{critCRT} implies the desired result. To that end, we fix $s\ino V$ and $s^\prime \ino [0, \zeta]$ be distinct from $s$. Suppose that $h(s)\! = \! m_h (s, s^\prime)$; if there is $s^{\prime \prime} \ino (s\wedge s^\prime, s\vee s^\prime)$ such that $h(s^{\prime \prime}) \! >\! m_h(s, s^\prime)$, then $p_h(s) \! \in \! \lgeo r, p_h(s^{\prime \prime}) \lgeo \, $ which is impossible since $p_h(s)$ has to be a leaf of $T_h$. Therefore, if $h(s)\! = \! m_h (s, s^\prime)$, then 
$h(s^{\prime \prime}) \! =\! m_h(s, s^\prime)$ for all $s^{\prime \prime} \ino [s\wedge s^\prime, s\vee s^\prime]$; thus $ [s\wedge s^\prime, s\vee s^\prime] \! \subset \! p_h^{-1} (\{ p_h(s)\})$ and $\mu_h$ would have an atom at $p_h(s)$, which is impossible since $T_h$ is a continuum real tree. Thus, we have proved that $h(s) \! >\! m_h (s, s^\prime)$. Therefore $\min_{r\in [m_h(s,s^\prime) , h(s) ]} w_s(r) \! < \! \widehat{w}_s$ since $s\ino U_{h,w}$. Next, we set 
$b\! = \! \min \{ w_s(r)\, ; r\ino [0, m_h(s,s^\prime) ] \}$; by Definition \ref{rePopol} $(b)$, 
$b$ is also equal to $ \min \{w_{s^\prime}(r) \, ;  r\ino [0, m_h(s,s^\prime) ] \}$. We also set $a\! = \!\min \{ w_s(r) \, ;  r\ino [m_h(s,s^\prime) , h(s)] \}$ and $a^\prime \! = \! \min \{w_{s^\prime} (r) \, ;  
r\ino [m_h(s,s^\prime) , h(s^\prime)] \}$. Then, 
\begin{eqnarray*}
d_{h,w} (0, s)+ d_{h,w} (s, s^\prime) \! -\! d_{h,w} (0, s^\prime) & =& 2 \widehat{w}_{s} +2M_{h,w} (0, s^\prime) -2M_{h,w} (0, s)-2M_{h,w} (s, s^\prime) \\
& =& 2 \widehat{w}_{s}-2a + 2 \big( a + (b\! \wedge \! a^\prime) \! -\! (b \! \wedge\!  a) \!  -\! (a \! \wedge \! a^\prime) \big) . 
\end{eqnarray*}
We next check that $ q\! :=\! a + (b\! \wedge \! a^\prime) \! -\! (b \! \wedge\!  a) \!  -\! (a \! \wedge \! a^\prime) \! \geq \! 0$. Indeed, if $a\! \leq \! a^\prime$, $q\! = \! (b\! \wedge \! a^\prime) \! -\! (b \! \wedge\!  a)  \! \geq \! 0$; if $b\! \leq \!  a^\prime \! \leq \! a$, $q\! = \! a\! -\! a^\prime \! \geq \! 0$; if $a^\prime \! \leq \!  b \! \leq \! a$, $q\! = \! a\! -\! b \! \geq \! 0$; if 
$a^\prime \! \leq \!  a \! \leq \! b$, $q\! = \! 0$. Since we have proved that $a\! < \! \widehat{w}_s$, we get $d_{h,w} (0, s)+ d_{h,w} (s, s^\prime) \! >\! d_{h,w} (0, s^\prime)$. This implies that $V \! \subset \! U$, which completes the proof. \cqfd

\medskip

Recall from (\ref{Traun}) the definition of the degree of a point in a real tree and recall that a point is a branch point if its degree is $\geq \!3$. 
The following lemma provides conditions ensuring that {snake} trees have only binary branch points. 
\begin{lemma}
\label{BruckanderMur} Let $(h, w) \ino \fSigma ([0, \zeta])$. Let $T_{h}$ and $T_{h,w}$ be the compact real trees coded by resp.~$h$ and $(h,w)$. Recall from Remark \ref{quotsnk} that $w$ can 
actually be defined on $T_h$. Let $\cD$ be a countable dense set of points of $T_h$. 
We consider two cases. In the first case, we assume the following. 

\begin{compactenum}

\smallskip

\item[$(a)$] $\widehat{w}_0\! = \! 0$ and $\widehat{w}_s \ino [0,\infty)$ for all $s\ino [0, \zeta]$.

\smallskip

\item[$(b)$] For all distinct $\sigma , \sigma^\prime \ino \cD$ such that 
$\min_{\gamma \in \lgeo \sigma , \sigma^\prime \rgeo} \widehat{w}_\gamma \! > \! 0 $, there is at most one $\sigma_0\! \in \,  \rgeo \sigma , \sigma^\prime \lgeo$ such that 
$\widehat{w}_{\sigma_0}\! = \!  \min_{\gamma \in \lgeo \sigma , \sigma^\prime \rgeo} \widehat{w}_\gamma$ and when there is one, it is never a branch point of $T_h$.  
\end{compactenum}

\smallskip

\noi
Then, for all $x\ino T_{h,w}\backslash \{ r_{h,w}\}$, $\mathtt{deg} (x) \ino \{ 1, 2, 3\}$. In the second case, we assume the following 

\begin{compactenum}

\smallskip

\item[$(a^\prime)$] $\widehat{w}_0\! = \! 0$.

\smallskip

\item[$(b^\prime)$] 
For all distinct $\sigma , \sigma^\prime \ino \cD$, 
there is at most one $\sigma_0\! \in \,  \rgeo \sigma , \sigma^\prime \lgeo$ such that 
$\widehat{w}_{\sigma_0}\! = \!  \min_{\gamma \in \lgeo \sigma , \sigma^\prime \rgeo} \widehat{w}_\gamma$ and when there is one, it is never a branch point of $T_h$.  
\end{compactenum}

\smallskip

\noi
Then, {for all} $x\ino T_{h,w}$, 
$\mathtt{deg} (x) \ino \{ 1, 2, 3\}$. 
\end{lemma}
\noi
\textbf{Proof.} Let us consider the first case. We assume $(a)$ and $(b)$. We shall 
argue by contradiction. To that end, suppose that $x_0\ino T_{h,w}\backslash \{ r_{h,w}\}$ is such that 
$\mathtt{deg} (x_0)\! \geq \! 4$. Recall from (\ref{Neunkirchen}) in Remark \ref{SktPaul} the definition of the continuous surjective function $\mathtt{y}\! : \! T_h \! \rightarrow \! T_{h,w}$. Then, choose $x_1, x_2, x_3\in \mathtt y (\mathcal{D})$ 
such that $r_{h,w}, x_1, x_2, x_3$ are in distinct connected components of $T_{h,w}\backslash \{ x_0\} $. Let $\sigma_0\ino T_h$ and $\sigma_1, \sigma_2, \sigma_3 \ino \mathcal{D}$ be such that $\mathtt{y} (\sigma _i)\! = \! x_i$, $i\ino \{ 0, 1, 2,3\}$. We first claim that for all 
$i, j\ino \{ 1, 2, 3\}$ distinct, 
\begin{equation}
\label{zerotoi}
0< d_{h,w}(r_{h,w}, x_0) = \widehat{w}_{\sigma_0}   =  \min_{ \gamma \in \lgeo \sigma_0, \sigma_i \rgeo}  \widehat{w}_\gamma \; \leq  \inf_{\gamma \in \lgeo \sigma_i, \sigma_j \rgeo} \widehat{w}_\gamma \; .
\end{equation}
\textit{Proof of (\ref{zerotoi}).} Note that $\mathtt{y} (r_h)\! = \! r_{h,w}$. 
Since $x_0\! \neq \! r_{h,w}$, by 
$(a)$ and (\ref{Judenburg}), we get $0 \! <\! d_{h,w} (x_0, r_{h,w}) \! = \! 
\widehat{w}_{\sigma_0} + \widehat{w}_{r_h}\!  -\! 
2\min_{\lgeo r_h, \sigma_0 \rgeo} \widehat{w}\! = \! \widehat{w}_{\sigma_0}$.  
Fix $i\ino \{ 1, 2, 3\}$. Since $x_0\! \in \, \rgeo r_{h,w}, x_i \lgeo $, (\ref{Judenburg}) implies 
$\widehat{w}_{\sigma_i}\! = \! d_{h,w} (r_{h,w}, x_i) 
 =   d_{h,w} (r_{h,w}, x_0)+ d_{h,w} (x_0, x_i)$ and thus 
 $ \widehat{w}_{\sigma_i}= 2\widehat{w}_{\sigma_0} + \widehat{w}_{\sigma_i}   - 2\min_{ \gamma \in \lgeo \sigma_0, \sigma_i \rgeo}  \widehat{w}_\gamma$, 
which implies the second equality in (\ref{zerotoi}). We complete the proof of (\ref{zerotoi}) by noting that $\lgeo \sigma_i, \sigma_j\rgeo \!  \subset \! \lgeo \sigma_0, \sigma_i\rgeo \cup \lgeo \sigma_0, \sigma_j\rgeo$. \cq 

\smallskip

For all $i\ino \{ 1, 2, 3\}$, let $C_i$ be the connected component of $T_{h,w}\backslash \{ x_0\}$ that contains $x_i$. By connectivity, 
there are $\gamma_1, \gamma_2 \ino \lgeo \sigma_1, \sigma_2 \rgeo$ 
such that $\mathtt y (\gamma_1)\! =\! \mathtt y (\gamma_2)\! =\!    x_0$ and 
$\lgeo \sigma_i , \gamma_i \lgeo  \subset \!   \mathtt y^{-1}(C_i)$, $i \ino \{ 1, 2\}$. Since 
$\mathtt{y} (\gamma_i) \! = \! x_0$, $\widehat{w}_{\sigma_0}+\widehat{w}_{\gamma_i}-2 \min_{\gamma \in \lgeo \sigma_0, \gamma_i \rgeo} \widehat{w}_\gamma= 0$, which implies that $\widehat{w}_{\gamma_i}\! = \! \widehat{w}_{\sigma_0}$. By (\ref{zerotoi}), we get  
$\widehat{w}_{\gamma_{1}}\! = \! \widehat{w}_{\gamma_{2}} = \min_{\gamma \in \lgeo \sigma_1, \sigma_2 \rgeo}  \widehat{w}_\gamma$. Then Assumption $(b)$ implies that 
$\gamma_1\! = \! \gamma_2 $. 

  We next introduce $\beta\! =\!  \mathtt{br} (\sigma_1, \sigma_2, \sigma_3)$, the branch point of $\sigma_1, \sigma_2, \sigma_3$ as defined in (\ref{braachpt}). Without loss of generality we can assume that $\gamma_1\ino \lgeo \sigma_1 , \beta \rgeo$ and thus $\gamma_1\ino \lgeo \sigma_1 , \sigma_3 \rgeo$. By connectivity, 
there is $\gamma_3 \ino \lgeo \sigma_1, \sigma_3 \rgeo$ 
such that $y(\gamma_3)\! =\!    x_0$ and $\lgeo \sigma_3 , \gamma_3 \lgeo  \subset \!   \mathtt y^{-1}(C_3)$. It implies that $\widehat{w}_{\gamma_1} \! = \! \widehat{w}_{\gamma_3} \! = \! \widehat{w}_{\sigma_0}\! = \! \min_{\gamma \in \lgeo \sigma_1, \sigma_3 \rgeo}  \widehat{w}_\gamma$. Then Assumption $(b)$ implies that $\gamma_3\! = \! \gamma_1\! = \! \gamma_2 $. 
By definition,  the subsets $\lgeo \sigma_i , \gamma_i \lgeo\, $, $i \ino \{ 1, 2, 3\}$ are pairwise disjoints and therefore $\gamma_3$ is a branch point of $T_h$ (it is actually $\mathtt{br} (\sigma_1, \sigma_2, \sigma_3)$) such that $\widehat{w}_{\gamma_3} \! = \! \min_{\lgeo \sigma_1, \sigma_3 \rgeo}  \widehat{w}_\gamma$ which is not possible by Assumption $(b)$. This proves the first point of the lemma by contradiction. 
The exact same arguments can be used to prove the second case of the lemma. \cqfd 

\subsection{Stable L\'evy trees.}
\label{statreesc} 
In this section that contains no new result, we briefly recall basic definitions and properties on stable L\'evy trees. L\'evy trees are 
a class of random compact metric spaces that have been introduced by Le Gall and Le Jan in \cite{LGLJ98} (and further studied in Le Gall and D.~\cite{DuLG02}) as the genealogy of continuous state branching processes. Among stable trees, 
Aldous' continuum random tree corresponds to the Brownian case (see Aldous \cite{Al93} {and here below)}. Stable trees (and more generally L\'evy trees) are the scaling limit of Galton--Watson trees as recalled in Introduction (see also below). 

More precisely, let $\gamma \in (1, 2]$ be the index 
of a spectrally positive stable L\'evy process $X\! = \! 
(X_s)_{s\in [0, \infty)}$: namely, the  law of $X$ is characterised by its Laplace exponent that is given by  
$$ \forall s ,\lambda \ino [0, \infty), \quad \bE [\exp (-\lambda X_s)] = \exp (s\lambda^\gamma) \; .$$
Note that $X$ is a Brownian motion when $\gamma\! =\!  2$ and we shall refer to this case as the {\it Brownian case}. As shown in \cite{LGLJ98} (see also \cite{DuLG02}, Chapter 1), there exists a {\it continuous process} $H(X)= (H_s(X) )_{ s \in [0, \infty)}$ such that for any $s \in [0, \infty)$,
\begin{equation}
\label{Hlimit}
H_s(X) = \lim_{\varepsilon\to 0} \frac{1}{\epp}\int_0^t {\bf 1}_{\{I^r_s<X_r<I^r_s+\epp\}}\,dr , 
\end{equation}
where $I^r_s$ stands for $\inf_{r\leq u\leq s} X_u$ and the convergence in \eqref{Hlimit} holds in probability.
The process $H (X)$ is the $\gamma$-{\it stable height process}. Note that in the Brownian case, $H(X)$ is simply the reflected Brownian motion $s\!\mapsto \! 2(X_s \! -\! I^0_s)$. From the scaling property of $X$ and from (\ref{Hlimit}), we see that for any $r \ino (0, \infty)$, 
$(r^{(\gamma-1)/\gamma}{H_{t/r}(X)} )_{t\in [0, \infty)}$ 
has the same distribution as $H(X)$.

As in the discrete setting, the process $H(X)$ encodes a family of continuum real trees: each excursion of $H(X)$ above $0$ corresponds to an excursion of $X$ above its infimum (it is obvious in the Brownian case) 
and each excursion of $H(X)$ above $0$ corresponds to a single continuum real tree of the family. 
The scaling property of $H(X)$ yields the following definition for the normalised excursion of $H(X)$ that is provided 
in D.~\cite{Du03} (see p.1005): 
set $\ell_1 \! = \! \max \{ s\ino [0, 1]\! : \! X_s\! \! = \! I_s \}$, 
$r_1 \! = \! \min \{ s\ino [1, \infty)\! : \! X_s\! \! = \! I_s \}$ and $\zeta_1\! = \! r_1 \! -\! \ell_1$. Then, 
we set:
\begin{equation*}
\forall s\ino [0, 1], \quad H_s := \zeta_1^{-(1-\frac{_1}{^\gamma}) } H_{\ell_1 + s \zeta_1} (X)\; .
\end{equation*}
The process $H$ 
is taken as the definition of the normalised excursion of the $\gamma$-stable height process.

As shown in Theorem 1.4.4 in Le Gall and D.~\cite{DuLG02}, for all {$\alpha \ino (0, 1 \! -\! \tfrac{_1}{^\gamma})$}, 
$\bP$-a.s.~$H(X)$ is $\alpha$-locally Hölder-continuous. 
Thus, the same holds true for $H$: namely, 
\begin{equation}
\label{steyr}
\textrm{{for all $\alpha \ino (0, 1 \! -\! \tfrac{_1}{^\gamma})$}, a.s.~$H$ is $\alpha$-Hölder continuous.}
\end{equation}

\noi 
\textbf{Normalised stable Lévy trees.} We call \textit{normalised $\gamma$-stable L\'evy tree} the real tree 
$$(T_H, d_H, r_H, \mu_H)$$ 
coded by the function $H$, as defined in (\ref{djhczkk}). Recall from 
(\ref{contitreed}) that $h\! \mapsto \! d_h$ is a $4$-Lipschitz function from 
$\bC([0, 1], \bbR)$ to $\bC([0, 1]^2, \bbR)$, both equipped with uniform norms. 
Thus, $d_H$ is a measurable random element of $\mathbf{MT} ([0, 1])$ 
and by Proposition \ref{groweak}, $T_H$ is a measurable random element of the Polish space 
$(\mathbb{M}, \bdelta_{\mathtt{GHP}})$ as {mentioned} in (\ref{rappGHP}).

Let us briefly recall some geometric properties of $T_H$. 
One can prove that a.s.~$T_H$ is a {continuum} real tree (as defined in (\ref{CRTdef})) 
and that the set of {its} branch points is a countable dense set; moreover 
we recall from Theorem 4.6 {in} Le Gall and D.~\cite{DuLG05} (p.~583) the following result.  
 \begin{equation} 
\label{Ennsdorf} 
\textrm{A.s.~for all $\sigma \ino T_H$, $\mathtt{deg} (\sigma) \ino A_\gamma$, where $A_\gamma$} 
= \left\{ 
\begin{array}{ll}
\{ 1,2, \infty\} &  \textrm{if $\gamma\ino (1, 2)$,} \\
\{1, 2, 3\} &  \textrm{if $\gamma \! =\! 2$.}
\end{array} \right.
\end{equation}

 \medskip
 
\noi \textbf{The contour of a discrete trees.} We briefly recall how to code discrete trees by various functions. Let $t\ino \bbT$ be a \textit{finite} (rooted ordered) {tree} as in Definition \ref{ort}. To simplify, we set $\# t\! =\! n$;      
since $t$ is finite, we can list the vertices of $t$ in the lexicographical order: $u_0\! = \! \varnothing\! <_t \! u_1 \! <_t \! \ldots \! <_t u_{n-1} $; for all $j\ino \{ 0, \ldots, n\! -\! 1\}$, we set $H_j (t)\! = \! |u_j|$ that is the height of the $j$-th vertex of $t$. By convenience we also set $H_n(t)\! = \! 0$. The function $(H_k(t))_{0\leq k\leq n}$ is the \textit{discrete height function of $t$}. Note that $H(t)$ entirely encodes $t$. 
 
 We also introduce another coding function known as the \textit{contour} (or the \textit{depth-first exploration}) function of $t$. Informally, we embed $t$ into the clockwise oriented half plane so that order on siblings corresponds to 
 orientation; we think of a particle that visits the tree at unit speed, that starts at the root and that goes from the left to right, backtracking as less as possible; we denote by $v(k)$ the vertex visited at time $k$ during this depth-first exploration. 
 The particle crosses exactly twice each edge of $t$: once upwards and once downwards; so, the total time needed to go back to the root is $2n\! -\! 2$, that is twice the number of edges. Namely,  
$\{ (v(k), v(k+1))  ; 0\! \leq \! k \! \leq \! 2n\! - \! 2 \} \! = \! 
\{ (u,\overleftarrow{u}), (\overleftarrow{u}, u); u\ino t\backslash \{ \varnothing\} \}$. Then, for all 
$k\ino \{ 0, \ldots, 2n\! -\! 2\}$, we set $C_k (t)\! = \! |v(k)|$ and we also set $C_{2n-1} (t)\! = \! C_{2n} (t)\! = \! 0$. 
We call \textit{contour function} of $t$ the linear interpolation of the $C_k (t)$ that we still denote by $(C_s(t))_{s\in [0, 2n]}$.

The discrete height function and the contour function of $t$ are related as follows. 
For all $0\!\leq \! j\! \leq  \! n$, we set $b_j\! = \! 2j \! -\! H_j(t)$; we easily check that $(b_j)$ is an increasing sequence from $0$ to $2n$. Moreover, for all $j\! < \! n\! -\! 1$, observe that for all $s\ino [b_j, b_{j+1}]$, 
 \begin{equation} 
\label{Crigo}
C_s(t)
= \left\{ 
\begin{array}{ll}
H_ j(t) \!  -\! (s\! -\! b_j) &  \textrm{If $s\ino [b_j, b_{j+1} \! -\! 1] $,} \\
s\! -\! b_{j+1}\! +\! H_{j+1} (t) &  \textrm{If $s\ino [b_{j+1} \! -\! 1 , b_{j+1}]$}
\end{array} \right.
\end{equation}
and that $C_{s}(t)\! = \! H_{n-1} (t)-(s\! -\! b_{n - 1}  )$ if $s\ino [b_{n-1}, b_n]$. 

The pointed measured compact real tree coded by the contour function of $t$ is described as follows: to simplify notation, we set $C_s \! =\! C_{2s} (t)$, $s\ino [0, n]$ and we recall from Example \ref{rempli} the definition of $(\widetilde{t}, \widetilde{d}^{{\, t}}_{\mathtt{gr}})$, the compact real tree spanned by the graph-tree $t$. It is easy to check that the pointed compact measured real tree $(T_{C},d_C, r_C, \mu_C)$ coded by $C$ is isometric to 
$(\widetilde{t}, \widetilde{d}^{\, t}_{\mathtt{gr}}, \varnothing, \widetilde{\mathtt{m}})$ where 
$\widetilde{\mathtt{m}}\! = \! \delta_\varnothing + \mathtt{Length}$, where $\delta_\varnothing$ stands for the Dirac mass at the root $\varnothing$. 

\smallskip

We next recall the following limit theorem from D.~\cite{Du03} (Theorem 3.1, p.~1006). 
\begin{theorem}[Theorem 3.1 \cite{Du03}]
\label{rappcvtr} We fix $\gamma \ino (1, 2]$. Let $\tau$ be a Galton--Watson tree whose offspring distribution satisfies $(\mathbf{H})$ as in (\ref{hyposta}). 
Then, there exists a nondecreasing sequence of positive real numbers $(a_n)_{n\in \bbN}$ that tends to $\infty$ such that 
\begin{eqnarray*}
\lefteqn{ \Big( \tfrac{1}{a_n} H_{\lfloor ns \rfloor}(\tau) \, ; \, \tfrac{1}{a_n} C_{2ns}(\tau) \Big)_{\! \! s\in [0,1]} \; \textrm{under $\; \bP \big(  \cdot  \big| \, \# \tau \! = \! n \big) $ } } \hspace{40mm} 
                    \nonumber\\
 & &\overset{\textrm{(law)}}{\underset{n\rightarrow \infty}{-\!\!\! -\!\!\!-\!\!\! -\!\!\!  -\!\! \! \longrightarrow}} \; (H_s, H_s)_{s\in [0, 1]}  
\end{eqnarray*}
where $H$ is the normalised excursion of the $\gamma$-stable height process. 
\end{theorem}

\subsection{One-dimensional reflected Brownian snakes.}
\label{rfbrsnsc}
In this section we briefly introduce normalised one-dimensional (reflected) Brownian snake with $\gamma$-stable branching mechanism (i.e.~whose lifetime process is a $\gamma$-stable height process). 
For more details  we refer to the monograph of Le Gall and D.~\cite{DuLG02} (Chapter 4, pp.~107-149). 

Let $I$ be an interval of $\bbR$ whose interior is not empty. 
Let $(Z^{x}_{{r}})_{r\in [0, \infty)}$ be an $I$-valued continuous Markov process starting at $Z^x_0\! = \! x\ino I$. We shall restrict to the two following cases: 
\begin{compactenum}

\smallskip

\item[$\bullet$] either $I$ is $\bbR$ and $Z^x$ is  
a \textit{Brownian motion}, 

\smallskip

\item[$\bullet$] or $I\! = \! [0, \infty)$ and $Z^x$ is 
a \textit{Brownian motion that is reflected at $0$}. 
\end{compactenum}

\smallskip

\noi
Recall from Definition \ref{rePopol} that $\bC_0$ stands for the space of continuous functions from $[0, \infty)$ to $\bbR$ equipped with (Polish) topology of uniform convergence on all compact subsets. We let $w\ino \bC_0$ be $I$-valued and we 
fix two nonnegative real numbers $b\! \geq \! a $. 
We then denote by $R_{a,b}(w, dw^\prime)$ the law on $\bC_0$ of the process $W(\cdot)$ that is defined as follows: 
\begin{compactenum}

\smallskip

\item[$\bullet$] for all $r\ino [0, a]$, $W(r)\! = \! w(r)$;

\smallskip

\item[$\bullet$]  the process $(W(a+r))_{r\in [0, \infty)}$ has the same law as $(Z^{w(a)}_{r\wedge (b-a)})_{r\in [0, \infty)}$.
\end{compactenum}    

\smallskip

\noi
It is clear that $(a,b,w)\!  \longmapsto \! R_{a,b} (w,dw^\prime)$ is weakly continuous on $\bC_0$. 

  We next fix $\zeta \ino (0, \infty)$ and 
$h\ino \bC([0, \zeta], [0, \infty))$ such that $h(0)\! = \! h(\zeta)\! = \! 0$.
Recall from (\ref{djhczkk}) the definition of $m_h (\cdot, \cdot)$. It is easy to check that one 
defines a $\bC_0$-valued process $(W_s (\cdot))_{s\in [0, \zeta]}$ by specifying its 
finite dimensional marginals as follows: 
\begin{compactenum}

\smallskip

\item[$\bullet$] $W_0 (r)\! = \! 0$, $r\ino [0, \infty)$.  

\smallskip
\item[$\bullet$] For all $0 \! \leq \! s_1 \! \leq \! \ldots \! \leq  \! s_p \! \leq \! \zeta$, $(W_{s_k})_{1\leq k\leq p}$ has law  
\end{compactenum} 
 \begin{equation}
 \label{margisn} 
 R_{m_h(0, s_1) , h(s_1)} (W_0, dw_1) R_{m_h(s_1, s_2) , h(s_2)} (w_1, dw_2) \ldots R_{m_h(s_{p-1}, s_p) , h(s_p)} (w_{p-1}, dw_p) \; .
 \end{equation}
When the $Z^x$ are Brownian motions, the resulting collection of 
$\bC_0$-valued r.v.~$(W_s)_{s\in [0, \zeta]}$ is a 
\textit{Brownian snake with lifetime process $h$ and initial position $0$}. 
When the $Z^x$ are reflected Brownian motions, $W$ is a \textit{reflected Brownian snake with lifetime process $h$ and initial position $0$}. The following $\bbR$-valued process \begin{equation*}
\forall s\ino [0, \zeta], \quad \widehat{W}_s= W_s (h(s)), 
\end{equation*}
is the \textit{endpoint process} of the snake $W$. 
\begin{remark}
\label{Raxalpen} If $(W_s)_{s\in [0, \zeta]}$ is a Brownian snake with lifetime process $h$ and initial position $0$, then clearly $(|W_s|)_{s\in [0, \zeta]}$ is a reflected Brownian snake with lifetime process $h$ and initial position $0$. \cq 
 \end{remark}

Let us discuss the regularity of Brownian snakes. First, denote by $M^{_h}_{^{s_1, \ldots , s_p}} $ the law on $\bC_0^p$ defined by (\ref{margisn}). It is easy to check that $(s_1, \ldots, s_p; h)\ino [0, \zeta]^p\! \times \! \bC([0, \zeta], \bbR) \! \longmapsto \!   M^{_h}_{^{s_1, \ldots , s_p}}$ is weakly continuous on $\bC_0^p$. Although $h$ is continuous, the process $W$ may not be continuous in general.  
However, the arguments used in the proof of Proposition~4.4.1 in Le Gall and D.~\cite{DuLG02} (Chapter 4 p.~127) show that if $h$ is Hölder, then $W$ has a continuous modification. 
More precisely, for all $\zeta, c\ino (0, \infty) $ and all $\alpha \ino (0, 1]$, denote by $\texttt{Höl}_{c, \alpha}([0, \zeta])$ the set of functions 
$h\ino \bC([0, \zeta], [0, \infty))$ such that $h(0)\! = \! h(\zeta) \! = \! 0$ and $|h(s)\! -\! h(s^\prime)|\! \leq \! c|s\! -\! s^\prime|^\alpha $, for all {$s , s^\prime \ino [0, \zeta]$}. From the proof of Proposition~4.4.1 in Le Gall and D.~\cite{DuLG02}, we easily  
get the following lemma. 
\begin{lemma}
\label{Kolsnake} Let $\zeta \ino (0, \infty)$ and let $h\ino \mathtt{Ho}^{\!\!\!\!^{..}} \mathtt{l}_{c, \alpha}([0, \zeta])$. 
Let $(W_s)_{s\in [0, \zeta]}$ be a Brownian snake (or a reflected Brownian snake) 
with lifetime process $h$ and initial position $0$. Then, for all $p\ino (1, \infty)$, 
\begin{equation*}
\forall s, s^\prime \ino [0, \zeta], \quad \bE \Big[ \!\!\!\! \! \sup_{\quad r\in [0, \infty)} \!\!\!\! \!\! \big| W_s (r) \! -\! W_{s^\prime} (r)\big|^p  \, \Big] \leq c_p |s\! -\! s^\prime|^{\alpha p/2} \; , 
 \end{equation*} 
where $c_p\! = \! (2p \sqrt{c}/(p\! -\! 1))^p \bE [|Z^0_1|^p]$. By Kolmogorov's criterion, $W$ admits a modification that is $\alpha/2$-Hölder with respect to the distance $\delta_{\mathtt{u}}$ on $\bC_0$ (see (\ref{uttcpct}) for the definition of 
$\delta_{\mathtt{u}}$). 
\end{lemma}

We now assume that $h$ is Hölder; we keep denoting $W$ the continuous version of the (possibly reflected) Brownian snake with lifetime process $h$ and we briefly discuss 
properties of $W$ for later use.

First observe that a.s.~$(h,W)$ is a snake 
as in Definition \ref{rePopol} $(b)$. 
Then, note that for all $s\ino [0, \zeta]$, $(W_s(r))_{r\in [0, \infty)}$ 
has the same distribution as $(Z^0_{r\wedge h(s)} )_{r\in [0, \infty)}$, 
which easily implies the following: for all bounded and measurable 
$F\! : \! \bC_0 \! \rightarrow \! \bbR$, 
\begin{equation}       
\label{onept}
\bE \Big[ \int_0^\zeta \!\!\! F \big(W_s (\cdot)\big)\,  ds \Big]= \int_0^\zeta \!\!\! \bE \big[ F \big( \big( Z^0_{r\wedge h(s)} \big)_{r\in [0, \infty)}\big)\big] \, dr \; . 
\end{equation}

To define the normalised (reflected) Brownian snake with $\gamma$-stable branching mechanism, we need to shows that the law of a (reflected) Brownian snake is a Borel-measurable function of 
its lifetime process. More precisely, 
since $\texttt{Höl}_{c, \alpha}([0, \zeta])$ is a compact subset of $\bC([0, \zeta], [0, \infty))$ equipped with uniform convergence and since 
$\texttt{Höl}_{c, \alpha}([0, \zeta]) \! \subset \! \texttt{Höl}_{\zeta^{\alpha -\alpha^\prime} 
c^\prime, \alpha^\prime}([0, \zeta]) $ for all $0\! < \! \alpha^\prime \! \leq \! \alpha$ and all $c^\prime \! \geq \! c\! >\! 0$, the set of Hölder lifetime-functions  
$\texttt{Höl}\! := \! \bigcup_{c\in (0, \infty)} \bigcup_{\alpha \in (0,1] } \texttt{Höl}_{c, \alpha}([0, \zeta])$ is a countable union of compact subsets and therefore a Borel subset of $\bC ([0, \zeta], [0, \infty))$. For all 
$h \ino \texttt{Höl}$, denote by $P_h( dW)$ (resp.~$P^+_h (dW)$) 
the law on $\bC ([0, \infty), \bC_0)$ of the Brownian snake (resp.~of the reflected Brownian snake) 
$W$ with lifetime-function $h$ and initial position $0$ and extend $h\! \mapsto\! P_h$ 
(resp.~$h\! \mapsto\! P^+_h$) to $\bC ([0, \zeta], [0, \infty))$ trivially by taking $P_h$ (resp.~$P^+_h$) equal to the Dirac mass at the null function if $h$ is not a Hölder lifetime-function. 
Now recall that for all integers $p\! \geq \! 1$ and for all bounded 
continuous $F\! : \! \bC_0^p \! \rightarrow \! \bbR$, the map $(s_1, \ldots, s_p; h) \mapsto \int F(W_{s_1}, \ldots, W_{s_p} ) \, P_h (dW)$ is continuous; the same holds with $P^+_h$. 
Since the Borel sigma field of $\bC([0, \zeta], \bC_0)$ equipped with $\Delta^{\! *}$ (as in Definition \ref{rePopol} $(a)$) is generated by the finite dimensional marginals, a standard monotone class argument implies that for all Borel subsets $B$ of the Polish space $(\bC([0, \zeta], \bC_0), \Delta^{\! *})$, the map $h\ino \bC([0,\zeta], [0, \infty)) \! \longmapsto \! P_h (B)\ino [0, \infty)$ is Borel-measurable. Namely: 

\smallskip

\begin{compactenum}
\item[] \textit{$P_h(dW)$ is a Borel regular kernel from $\bC([0,\zeta], [0, \infty))$ equipped with the uniform norm to the Polish space $(\bC([0, \zeta], \bC_0), \Delta^{\! *})$. The same holds with $P^+_h$}. 
\end{compactenum}
\begin{definition}
\label{brosnadef} Let $\gamma \ino (1, 2]$.  
Let $H$ be the normalised excursion of the $\gamma$-stable height process. 
Recall from (\ref{steyr}) that $H$ is a.s.~Hölder-continuous. Thus, it makes sense to define the 
\textit{Brownian snake} (resp.~\textit{reflected Brownian snake}) \textit{with $\gamma$-stable branching mechanism} as the $\bC([0, \zeta], \bbR) \! \times \! \bC([0, \zeta], \bC_0)$-valued r.v.~$(H, W)$ whose regular conditional distribution given $H$ is $P_H (dW)$ 
(resp. $P^+_H(dW)$), as defined previously.   \cq
\end{definition}

  Let $(H, W)$ be a normalised Brownian snake with {$\gamma$-stable branching mechanism}. Recall that the real tree $(T_H, d_H, r_h, \mu_H)$ coded by $H$ is the $\gamma$-stable Lévy tree. 
Since a.s.~$(H,W)\ino \fSigma ([0, 1])$, (\ref{snktree}) in Remark \ref{quotsnk} applies to $(H, W)$ and 
it makes sense  (up to a slight abuse of notation) 
to define $W$ and $\widehat{W}$ on $T_H$. To prove geometric properties of the $\gamma$-stable (reflected) Brownian cactus, we shall need the following result that is recalled from 
Lemma 6.4 in Le Gall and D.~\cite{DuLG05} (p.~600). 
\begin{lemma}[Lemma 6.4 \cite{DuLG05}]
\label{Bregenz} 
Let $(H, W)$ be a normalised Brownian snake with {$\gamma$-stable branching mechanism}. Then, conditionally given $H$, $(\widehat{W}_\sigma)_{\sigma \in T_H}$ is a centered Gaussian process whose covariance is characterised by the following: 
\begin{equation}
\label{Mayerling} 
\forall \sigma, \sigma^\prime \ino T_H, \quad \int \!\! P_H (dW) \, \big| \widehat{W}_{\sigma^\prime} \! - \widehat{W}_{\sigma} \big|^2 = d_H (\sigma, \sigma^\prime) \; .
\end{equation}
Moreover, for all $\epp \ino (0, 1/2)$, 
conditionally given 
$H$, $\sigma \ino T_H \! \longmapsto \! \widehat{W}_\sigma$ is $\big(\frac{1}{2} \! -\! \epp \big)$-Hölder-continuous. 
\end{lemma}

\subsection{Reflected Brownian cactus with stable branching mechanism.} 
\label{Propcact} 
\begin{definition}
\label{brocacdef} Let $\gamma \ino (1, 2]$ and let $(H, W)$ be a normalised Brownian (resp.~reflected Brownian) snake with {$\gamma$-stable branching mechanism} as in Definition \ref{brosnadef}. 
Recall the notation $\fSigma ([0, 1])$ from Definition \ref{rePopol} $(b)$ 
and recall that $(H,W)\ino \fSigma ([0, 1])$. The 
\textit{Brownian} (resp.~\textit{reflected Brownian}) 
\textit{cactus with $\gamma$-stable mechanism} is the pointed measured compact real tree $(T_{H,W}, d_{H,W},  r_{H,W},  \mu_{H,W})$ as specified in Definitions \ref{defsnkme} and \ref{notasnkm}. \cq 
\end{definition}
Recall from (\ref{CRTdef}) the definition of continuum real trees.   
Recall from (\ref{Traun}) the definition of the degree of a point in a real tree. 
\begin{proposition}  
\label{Dornbirn} Let $\gamma \ino (1, 2]$. Let $(H, W)$ be a normalised Brownian (resp.~reflected Brownian) snake with $\gamma$-stable branching mechanism and let $T_{H, W}$ the corresponding cactus as in Definition \ref{brocacdef}. 
Then, the following holds true. 

\begin{compactenum}

\smallskip

\item[$(i)$] $(T_{H,W}, d_{H,W},  r_{H,W},  \mu_{H,W})$ is a.s.~a continuum real tree.

\smallskip

\item[$(ii)$] A.s.~for all $x\ino T_{H, W} \backslash \{ r_{H,W} \}$, $\mathtt{deg} (x) \ino \{ 1,2, 3\}$. Moreover, if $T_{H,W}$ is a reflected Brownian cactus, then $\mathtt{deg} ( r_{H,W})\! = \! \infty$ 
and if $T_{H,W}$ is a Brownian cactus, then $\mathtt{deg} ( r_{H,W})\! = \! 1$. 

\smallskip

\end{compactenum}
\end{proposition}
\begin{rem}
\label{Traunsee1} As recalled in (\ref{Ennsdorf}), when $\gamma\ino (1, 2)$, the $\gamma$-stable L\'evy tree $T_H$ has infinite branch points: namely, there is a countable dense set of points $\sigma\ino T_H$ such that $\mathtt{deg} (\sigma)\! = \! \infty$. However, note that the corresponding Brownian cactus $T_{H,W}$ has only binary branch points (except possibly at the root in the reflected case). \cq 
\end{rem}

\begin{rem}
\label{Traunsee2} As already mentioned, the Brownian Cactus has been introduced by N.~Curien, J-F.~Le Gall and G.~Miermont in \cite{CuLGMi13} to study planar maps: roughly speaking the Brownian cactus {corresponds} to the case of a quadratic branching mechanism $\gamma\! = \! 2$ and the spatial motion of the snake is a (unreflected) linear Brownian motion. In this case they proved much finer geometrical results: {a.s.~}the upper-local density for typical points is $4$ (Proposition 5.1 {\cite{CuLGMi13}} p.~364) and the Hausdorff dimension is $4$ (Corollary 5.3 {\cite{CuLGMi13}} p.~365). See also the article by J-F.~Le Gall \cite{LG15} where the level sets of the Brownian cactus are studied to derive results on the Brownian maps.  \cq 
\end{rem}
\noi
\textbf{Proof of Proposition \ref{Dornbirn}.} To prove $(i)$, we use Lemma \ref{snkCRT}. Indeed recall that $(T_H, d_H, r_H, \mu_H)$ is a.s.~a continuum real tree and recall from (\ref{Uhwdef}) the definition of the subset $U_{H,W} \! \subset \! [0, 1]$. Recall that $(Z^0_r)_{r\in [0, \infty)}$ stands for a one-dimensional Brownian motion (resp.~reflected Brownian motion) {starting at 0} that is independent from $H$. We denote by $\ell$ the Lebesgue measure on $[0, 1]$.
Then, by (\ref{onept}), conditionally given $H$, we get a.s.~
\begin{eqnarray*} 
\int \! P_H (dW) \ell (U_{H, W})\!\!  &  = & \!\!\! 
\int_0^1 \!\! \un_{\{ H_s   >  0 \} } \,  P_H \Big(  \, \forall  \epp \ino (0, H_s)   : \!   
\!\!\!\!\!\!\!\!\!\!  \!\! \inf_{\qquad r\in [H_s -\epp, H_s] }\!\!\!\!\!  \!\!\!\!\!  \!\! W_s (r) \! <  \widehat{W}_s \Big) \, ds  \\ 
&= & \int_0^1 \!\!  \un_{\{ H_s   >  0 \} } \,   \bP 
\Big( \,     \forall \epp \ino (0, H_s)   : \!   
\!\!\!\!\!\!\!\!\!\!  \!\! \inf_{\qquad r\in [H_s -\epp, H_s] }\!\!\!\!\!  \!\!\!\!\!  \!\! \!\! Z^0_r \! < Z^0_{H_s}  \, \Big| \, H  \Big) \, ds  = 1 . 
\end{eqnarray*}
Thus, a.s.~$\ell ([0, 1]\backslash U_{H, W})\! = \! 0$ and Lemma \ref{snkCRT} implies $(i)$. 

  We next use Lemma \ref{BruckanderMur} to prove $(ii)$. Recall that (up to a slight abuse of notation) $W$ can be defined on $T_H$. We first consider the {unreflected} Brownian case and recall from Lemma \ref{Bregenz} that conditionally given $H$, 
$(\widehat{W}_\sigma)_{\sigma\in T_H}$ is a centered Hölder-continuous Gaussian process whose covariance is {characterised} by (\ref{Mayerling}). 
According to Remark \ref{Raxalpen}, 
$|\widehat{W}|$ is the endpoint process of the $\gamma$-stable reflected Brownian snake. 
We now work conditionally given $H$. 

Let $\cD$ be a countable dense subset of $T_H$. 
Let $\sigma , \sigma^\prime \ino \cD$ be distinct; denote by $\phi\! : \! [0, d_H (\sigma, \sigma^\prime)] \! \rightarrow \! \lgeo \sigma, \sigma^\prime \rgeo_{T_H}$ the isometry such that $\phi (0)\! = \! \sigma$ and $\phi (d_H (\sigma, \sigma^\prime)) \! = \! \sigma^\prime$. The covariance of $\widehat{W}$ combined with the fact that Brownian motion is reversible entails that the process 
$Z_r\! := \! \widehat{W}_{\phi (r)} \! -\! \widehat{W}_{\sigma}$, $r\ino [0, d_H (\sigma , \sigma^\prime)]$, has the same law conditionally given $H$ as a 
Brownian motion starting at $0$: it therefore a.s.~reaches its infimum at a unique time $r_0 \! \in ]0, d_H (\sigma, \sigma^\prime)[ \,$ and the 
law of $r_0$ conditionally given $H$ is diffuse, which entails that a.s.~$\phi (r_0)$ cannot be a branch point of $T_H$ since the set of branch points of $T_H$ are countable. The same holds true for $|\widehat{W}| $ on 
$\lgeo \sigma, \sigma^\prime\rgeo_{T_H}$ if $\inf_{\gamma \in  \lgeo \sigma , \sigma^\prime \rgeo_{T_H}} |\widehat{W}_\gamma| \! >\! 0$. Since $\cD$ is countable, it implies 
that conditionally given $H$, $\widehat{W}$ satisfies Assumptions $(a^\prime)$ and $(b^\prime)$ in  Lemma \ref{BruckanderMur} and that $|\widehat{W}|$ satisfies Assumptions $(a)$ and $(b)$ in Lemma \ref{BruckanderMur}. This lemma then implies the first assertion of $(ii)$.  
   
 Let us complete the proof of $(ii)$. Recall from Remarks \ref{SktPaul} that it is makes sense to define $\mathtt{y}^+\! : \! T_H \! \rightarrow \! T_{H, |W|}$ by setting 
$\mathtt{y}^+ (p_H (s))\! = \! p_{H, |W|} (s)$, $s\ino [0, 1]$. Let $\sigma\ino T_H$ be such that $\mathtt{y}^+ (\sigma) \! \neq \! r_{H, |W|}$. Then, observe that $\{ \gamma \ino \lgeo r_H , \sigma \rgeo_{T_H} \! : \! |\widehat{W}_\gamma| \! >\! 0 \} $ has infinitely many connected components: let $\gamma$ and $\gamma^\prime$ be in two distinct such connected components: by (\ref{Judenburg}) we get 
$d_{H, |W|} (\mathtt{y}^+ (\gamma), \mathtt{y}^+ (\gamma^\prime)) \! = \! \widehat{W}_{\gamma} +  \widehat{W}_{\gamma^\prime}\! = \! d_{H, |W|} (r_{H,W}, \mathtt{y}^+ (\gamma^\prime)) +d_{H, |W|} (r_{H,W}, \mathtt{y}^+ (\gamma))$, which easily implies that $\mathtt{y}^+ (\gamma)$ and $\mathtt{y}^+ (\gamma^\prime)$ belong to distinct connected components of $T_{H, |W|}\backslash \{ r_{H, |W|} \}$. 
It implies that a.s.~$\mathtt{deg} (r_{H, |W|})\! = \! \infty$. 

 We define similarly $\mathtt{y}\! : \! T_H \! \rightarrow \! T_{H, W}$ by setting  $\mathtt{y} (p_H (s))\! = \! p_{H, W} (s)$, $s\ino [0, 1]$. Let $x, x^\prime$ and $r_{H, W}$ be distinct points of $T_{H, W}$ 
and denote by $z$ their branching point; let $\sigma $ and $\sigma^\prime \ino T_H$ be such that 
$\mathtt{y} (\sigma)\! = \! x $ and $\mathtt{y} (\sigma^\prime)\! = \! x^\prime $ and recall that $\mathtt{y} (r_H)\! = \! r_{H, W}$. Then 
\begin{eqnarray*}
2 d_{H, W} (r_{H, W} , z)&= &   d_{H, W} (r_{H, W} , x)+ d_{H, W} (r_{H, W} , x^\prime)- { d_{H, W}
(x, x^\prime) } \\
&  \overset{\textrm{by (\ref{Judenburg})}}{=}& 2 \!\! \!\! \!\!  \!\! \!\!   \!\! \min_{\qquad \gamma \in \lgeo \sigma , \sigma^\prime \rgeo_{T_H} }\!\! \!\!  \!\! \!\! \!\! \widehat{W}_\gamma \; -   2 \!\!  \!\! \!\! \!\!  \!\! \!\! \min_{\qquad \gamma \in \lgeo r_H  ,\sigma  \rgeo_{T_H} } \!\! \!\! \!\!  \!\! \!\!  \widehat{W}_\gamma \; - 2 \!\! \!\! \!\!  \!\! \!\! \!\!  \min_{\qquad \gamma \in \lgeo r_H  , \sigma^\prime \rgeo_{T_H} } \!\! \!\! \!\!  \!\! \!\!  \widehat{W}_\gamma . 
\end{eqnarray*}
Now observe that a.s.~$ \min_{\gamma \in \lgeo r_H  ,\sigma  \rgeo_{T_H} } \! \!   \widehat{W}_\gamma $ and  $ \min_{\gamma \in \lgeo r_H  ,\sigma^\prime  \rgeo_{T_H} } \!\!  \widehat{W}_\gamma $ are strictly negative quantities and that the least of the two has to be {smaller} than $\min_{\gamma \in \lgeo \sigma , \sigma^\prime \rgeo_{T_H} }\! \! \widehat{W}_\gamma$. This prove that $d_{H, W} (r_{H, W} , z)\! >\! 0$ and that $x$ and $x^\prime $ are in the same connected component of $T_{H, W}\backslash \{ r_{H, W} \}$. Thus, $\mathtt{deg} (r_{H, W} )\! = \! 1$, which 
completes the proof of $(ii)$.  \cqfd

\smallskip

\section{{Scaling limit of the range of the BRW.}}
\label{Pflimth}
\subsection{Continuous interpolation of discrete snakes.} 
\label{intpolsnk}
We briefly devise here a natural 
way to embed a tree-valued branching walk into a continuous branching motion, which takes its values in 
the real tree spanned by a discrete space-tree and is indexed by the real tree spanned 
by a genealogical tree.

To that end, let us consider a (possibly infinite) graph tree $T$ equipped with its graph-distance $d^{_T}_{^\mathtt{gr}}$. {$T$ will be the space tree.} We distinguish a special vertex $\mathtt{o} \ino T$ that plays the role of a root. 
Here, $T$ can be a rooted ordered tree as in Definition \ref{ort} or the $\mathtt{b}$-ary tree $\bbW_{\! \mathtt{b}} $ or the free tree $\bbW_{\! [0, 1]}$ as introduced in (\ref{Wdesun}). Recall from Example \ref{rempli} the definition of $(\widetilde{T}, \widetilde{d}^{_{\, T} }_{^{\mathtt{gr}}})$ of the real tree spanned by $T\! $. We make the necessary identifications to assume that $T \! \subset \! \widetilde{T}\! $.

We first explain how to obtain a continuous interpolation of a path $v(k)\ino T$, $k\ino \{ 0, \ldots, n \}$ (by a path, we mean a sequence of adjacent vertices). It is easy to check that there exists a continuous map $\widetilde{v}\! : \! [0, n] \! \rightarrow \! \widetilde{T}$ such that $\widetilde{v} (k)\! = \! v(k)$ and for all $0\! \leq \! k\! < \! n$, 
\begin{equation*}     
\forall s\ino [k, k+1], \quad \widetilde{v} (s)\ino \lgeo v(k), v(k\! +\! 1) \rgeo \;\,  \textrm{is such that 
$\, \widetilde{d}^{\, \, T }_{\mathtt{gr}} (v(k) ,\widetilde{v} (s))\! = \! s\! -\! k$. }  
\end{equation*}
The continuous path $(\widetilde{v}(s))_{s\in [0,n]}$ is the continuous interpolation of {the} path $(v(k))_{k\in\{0, \dots, n\}}$. 

 We next extend this interpolation to $T$-valued branching walks. 
Let $t\ino \bbT$ be a \textit{finite} (rooted ordered) tree as in Definition \ref{ort}. {$t$ will be the genealogical tree.} We denote by 
$(\widetilde{t}, \widetilde{d}^{\, t }_{\mathtt{gr}})$ the compact real tree spanned by $t$ and we assume that $t\! \subset \! \widetilde{t}$. Let $(Y_v)_{v\in t}$ be a $T$-valued branching walk indexed by $t$. Namely, for all $v\ino t$, $Y_v\ino T$ and it  
satisfies the following conditions: $(a)$  $Y_\varnothing\! := \! \mathtt{o}$; $(b)$ if $v, v^\prime\ino t$ are neighbours, so are $Y_v$ and $Y_{v^\prime}$ in $T$. We easily check that there exists a continuous map 
$\widetilde{Y}\! : \! \widetilde{t} \! \mapsto \!  \widetilde{T}$ such that $\widetilde{Y}_{v}\! = \! Y_v$, $v\ino t$, by specifying the following: for all $\sigma \ino \widetilde{t}$ there exist two neighbouring vertices $v, v^\prime \ino t$ such that $\sigma \ino \lgeo v, v^\prime\rgeo_{\widetilde{t}}$ and 
\begin{equation}
\label{Yexten}
\textrm{$\widetilde{Y}_\sigma$ is the only point of $\lgeo Y_v, Y_{v^\prime} \rgeo_{\widetilde{T}}$ such that 
$\widetilde{d}^{\, T}_{\mathtt{gr}} (Y_v, \widetilde{Y}_\sigma)\! = \!  \widetilde{d}^{\, t}_{\mathtt{gr}} (v, \sigma)$.  }
\end{equation}

We next introduce the \textit{spatial contour} associated with $(Y_v)_{v\in t}$. To that end, denote by $v(k) \ino t$, $0\! \leq \! k \! \leq \! 2\#t $, the sequence of vertices of $t$ visited by the contour (or the depth-first) exploration of $t$ as recalled in (\ref{Crigo}) and 
note that they form a path in $t$. We denote by $\widetilde{v}\! : \! [0, 2\#t] \!\rightarrow \! \widetilde{t}$ its continuous interpolation. Observe that 
$C_s(t)\! = \!  \widetilde{d}^{\, t }_{\mathtt{gr}}( \varnothing, \widetilde{v}(s))$, for all $s\ino [0, 2\#t ]$, where $C(t)$ stands for the contour function of $t$. Then, for all $s\ino [0, 2 \# t]$, we define a continuous map 
$V_s\! : \! [0, \infty) \! \rightarrow \! \widetilde{T}$ as follows: 
\begin{equation} 
\label{spacctr} 
\left\{ 
\begin{array}{l}
\!\!  \textrm{for all $r\ino [0, C_s(t)]$, $V_s(r)\! = \! \widetilde{Y}_\sigma$ where $\sigma\ino \lgeo \varnothing , \widetilde{v} (s) \rgeo_{\widetilde{t}} $ is such that $ r\! = \!  \widetilde{d}^{\, t }_{\mathtt{gr}}( \varnothing,\sigma)$;} \\
\!\!  \textrm{for all $r \ino [C_s(t), \infty)$, $V_s(r) \! = \! \widetilde{Y}_{\widetilde{v} (s)} $.}
\end{array} \right.
\end{equation}

We call $(V_s (\cdot))_{s\in [0, 2\# t]}$ the \textit{spatial contour} associated with the branching walk $(Y_v)_{v\in t}$. 
It is easy to check here that $s\mapsto V_s$ is continuous from $[0, 2\# t]$ to the space of $\widetilde{T}$-valued continuous functions equipped with the uniform distance. We also define its endpoint process as follows: 
\begin{equation}
\label{kotschach}
\forall s\ino [0, 2 \# t], \quad \widehat{V}_s= V_s (C_s(t))\; .
\end{equation}

We then define the \textit{continuous interpolation of the discrete snake associated with $(Y_v)_{v\in t}$} and its endpoint process by 
\begin{equation}  
\label{contsnk}
\forall s\ino [0, 2\# t], \quad  \gzW_s (r)= \widetilde{d}^{{\, T}}_{{\mathtt{gr} }} (\mathtt{o}, V_s(r)), \; r \ino [0, \infty) \quad \textrm{and} \quad \widehat{\gzW}_s\! = \! \widetilde{d}^{\, T}_{\mathtt{gr} } (\mathtt{o}, \widehat{V}_s). 
\end{equation}
By Remark \ref{adjadja} below, $(\gzW_s (\cdot) )_{s\in [0, 2\#t]}$ is the continuous interpolation of the $\bbN$-valued branching walk 
$(d_{\mathtt{gr}}^T (\mathtt{o}, Y_v))_{v\in t}$. Moreover, observe that $\gzW$ is a snake whose lifetime function is the contour $C(t)$ of $t$ as in Definition \ref{rePopol}: namely, $(C_s (t), \gzW_{\! s} )_{s\in [0, 2\#t]}\ino \fSigma ([0, 2 \# t])$.

\begin{remark}
\label{adjadja}
Let $T^\prime$ be another graph-tree and let $\phi\! : \! T\! \rightarrow \! T^\prime$ be an adjacency-preserving map so that the image by $\phi$ of a path in $T$ is a path in $T^\prime$ (for instance, $T\! = \! \bbW_{\! [0, 1]}$, $T^\prime\! = \! \bbW_{\! \bgg}$ and $\phi\! = \! \Phi_\bgg$, the $\bgg$-contraction as in Definition \ref{tracage}; or $T^\prime\! = \! \bbN$ and $\phi (v)\! = \! d_{\mathtt{gr}} (\mathtt{o}, v)$) 
Consequently, $(Y^\prime_v)_{v\in t} \! := \! (\phi (Y_v))_{v\in t}$ is a $T^\prime$-valued branching walk indexed by $t$. 
Then, $\phi$ clearly extends to a continuous map from $\widetilde{T}\! \rightarrow \! \widetilde{T}^\prime$  by requiring that 
$ \widetilde{d}^{\, T^\prime}_{\mathtt{gr}} (\phi (\sigma), \phi (u))\! = \!   \widetilde{d}^{\, T}_{\mathtt{gr}} (\sigma, u)$ for all $\sigma \ino \lgeo u,v \rgeo_{\widetilde{T}}$ with $u$ and $v$ adjacent in $T$.   
Then, observe that the continuous interpolation of $Y^\prime$ is the image by $\phi$ of the continuous interpolation of $Y$: namely, $\widetilde{Y}^\prime\! \! = \! \phi (\widetilde{Y})$. \cq 
\end{remark}

We next provide basic properties related to counting and occupation measures on the range of $Y$ and its interpolation $\widetilde{Y}$. We first set  
\begin{equation}
\label{Weissbach}
\cR \! = \! \big\{ Y_v; v\ino t\big\} \quad \textrm{and} \quad \widetilde{\cR} \! = \! \big\{ \widehat{V}_s; s\ino [0, 2 \# t]\big\}= \big\{ \widetilde{Y}_\sigma; \sigma \ino \widetilde{t}\,  \big\}\; .
\end{equation}
Observe that $\cR$ {(resp.~$\widetilde{\cR}$)} is a subtree of $T$ {(resp.~$\widetilde{T}$)} and note that $\widetilde{\cR}$ is the real tree spanned by $\cR$. We equip $\cR$ and $\widetilde{\cR}$ with the respective occupation measures $\bm_{\mathtt{cont}}$ and $\widetilde{\bm}_{\mathtt{cont}}$ of the endpoint process of the spatial contour $\widehat{V}$ (here, the subscript $\mathtt{cont}$ is for \textit{contour}). Namely, for all bounded measurable $f\! : \! \widetilde{\cR} \! \rightarrow \! \bbR$, 
\begin{equation}
\label{Gloggnitz}
\int_{\cR} f \, d \bm_{\mathtt{cont}} = \! \sum_{0\leq k< 2 \# t}\!\!  f \big( \widehat{V}_k\big) \quad \textrm{and} \quad \int_{\widetilde{\cR}} f \, d \widetilde{\bm}_{\mathtt{cont}} = \! \int_0^{2 \# t} \!\!\! \! \!  f\big( \widehat{V}_s\big) \, ds . 
\end{equation}
We next introduce the \textit{spatial contour pseudo-distance} $d_{\mathtt{cont}}$ associated with $Y$: 
\begin{equation}
\label{Globasnitz} 
\forall s,s^\prime \ino [0, 2\# t], \quad d_{\mathtt{cont}} (s,s^\prime)= \widetilde{d}^{\, T}_{\mathtt{gr}} (\widehat{V}_s, \widehat{V}_{s^\prime} ) \; .
\end{equation}
Then, we easily check that $d_{\mathtt{cont}} \ino \MMT([0, 2\# t])$ and that 
\begin{equation}
\label{Troepolach} 
\textrm{the corresponding pointed measured compact real tree is isometric to $\big(  \widetilde{\cR},  \widetilde{d}^{\, T}_{\mathtt{gr}} , \mathtt{o}, 
\widetilde{\bm}_{\mathtt{cont}} \big)$.}  
\end{equation}
\begin{lemma}
\label{Tobiach} We keep the above notation. We denote by  $\widetilde{d}_{\mathtt{Haus}}$ the Hausdorff distance on the compact subsets of 
$(\widetilde{T}, \widetilde{d}_{^\mathtt{gr}}^{_{\, T}})$ and we denote by $\widetilde{d}_{\mathtt{Prok}}$ the Prokhorov distance on the finite Borel measures on $(\widetilde{T}, \widetilde{d}_{^\mathtt{gr}}^{_{\, T}})$. 
Then, $\widetilde{d}_{\mathtt{Haus}} (\cR, \widetilde{\cR}) \! \leq\!  1$, $\widetilde{d}_{\mathtt{Prok}} (\widetilde{\bm}_{\mathtt{cont}}, \bm_{\mathtt{cont}}) \! \leq \! 1 $. 
Moreover, set $\bm_{\mathtt{occ}}\! = \! \sum_{v\in t} \delta_{Y_{v}}$, 
that is the occupation measure of $Y$. Then, we also get 
$\widetilde{d}_{\mathtt{Prok}} (\bm_{\mathtt{cont}}, 2\bm_{\mathtt{occ}}) \! \leq \! 1 $ and thus, 
\begin{equation}
\label{Voelkermarkt}
\widetilde{d}_{\mathtt{Prok}} (\widetilde{\bm}_{\mathtt{cont}}, 2\bm_{\mathtt{occ}}) \! \leq \! 2 \; .
\end{equation}
\end{lemma}
\noi
\textbf{Proof.} Clearly, $\widetilde{d}_{\mathtt{Haus}} (\cR, \widetilde{\cR}) \! \leq\!  1$. 
Next denote by $\pi$ the image of the Lebesgue measure on $[0, 2\# t]$ 
on $\widetilde{T} \! \times \!  \widetilde{T} $ via the map $s\! \mapsto \! (\widehat{V}_{ \lfloor s \rfloor} , \widehat{V}_s)$. By the definition (\ref{Gloggnitz}), the first marginal of $\pi$ is $\bm_{\mathtt{cont}}$ and the second one is $\widetilde{\bm}_{\mathtt{cont}}$. Moreover by definition of the interpolation and of the spatial contour $ \widetilde{d}_{^\mathtt{gr}}^{_{\, T}} (\widehat{V}_{ \lfloor s \rfloor} , \widehat{V}_s) \! \leq \! 1$. This easily implies that 
$\widetilde{d}_{\mathtt{Prok}} (\widetilde{\bm}_{\mathtt{cont}}, \bm_{\mathtt{cont}}) \! \leq \! 1 $.

   To prove (\ref{Voelkermarkt}), we use the following coupling: recall that $(v(k))_{0\leq k\leq 2 \# t }$ stands for the contour (or depth-first) exploration of $t$ (and recall the convention $
   v(2\# t \! -\! 1)\! = \! v(2 \# t)\! = \! \varnothing$). For all $k\ino \{ 0, \ldots, 2\#t \! - \! 3\}$, we set $\rho (k)\! = \! v(k)$ if $|v(k)|\! = \! 1+ |v(k+1)|$ and $\rho (k)\! = \! v(k+1)$ if $|v(k+1)| \! = \! 1+|v(k)|$ and we also set  $\rho (2\# t \! -\! 2)\! = \! \rho(2 \# t\! -\! 1)\! = \! \varnothing$. Note that the image measure of the counting measure on $\{0, \ldots, 2\# t-1 \}$ via $\rho$ is $2\sum_{u\in t} \delta_u$, namely twice the counting measure on $t$ and thus, the image measure of the counting measure on $\{0, \ldots, 2\# t \! -\! 1 \}$ via $Y_{\rho (\cdot)}$ is $2\bm_{\mathtt{occ}}$. 
Then, denote by $\varpi$ the image measure on $T\times T$ of the counting measure on $\{0, \ldots, 2\# t \! -\! 1\}$ via the map $k \! \mapsto \! ( Y_{\rho (k)}, \widehat{V}_k)$. We just have proved that the first marginal of $\varpi$ is $2\bm_{\mathtt{occ}}$; the second marginal of $\varpi$ is by definition $\bm_{\mathtt{cont}}$. Since $Y_{\rho (k)}$ is either $\widehat{V}_k$ or $\widehat{V}_{k+1}$, we get $d_{^\mathtt{gr}}^{_T} (Y_{\rho (k)}, \widehat{V}_k) \! \leq \! 1$. 
This easily implies that $\widetilde{d}_{\mathtt{Prok}} (\bm_{\mathtt{cont}}, 2\bm_{\mathtt{occ}}) \! \leq \! 1 $ and (\ref{Voelkermarkt}) follows immediately from the previous bound. \cqfd

\medskip

In the following lemma we provide bounds to compare the occupation measure induced by the spatial contour and the counting measure on $\cR$. 
\begin{lemma}
\label{Loibach} We keep the above notation. To simplify, we set $n\! = \! \# t$ and $n^\prime \! = \! \# \cR$. 
We denote by $u_0, u_1, \ldots , u_{n-1}$ the sequence of the vertices of $t$ listed in the lexicographical order. Let $c\ino (0, \infty)$. We set $ \alpha \! =\!  \max_{1\leq i\leq n }  | \# \{ Y_{u_j}; 0\! \leq \! j \! < \! i \} \! - \! 2c i | $ and $\beta \! = \!  \max_{s \in (0, 2n] } | \#\{ \widehat{V}_{l} \, ; \; 0\! \leq \! l \! < \! \! \lceil  s \rceil \} \! -\!  c s |$. 
Then, 
\begin{equation}
\label{SktAndrae} 
\beta \leq \alpha + 3c + c \!\! \max_{s\in [0, 2n]} \!\! C_s (t) \; .
\end{equation}
We next denote by $\bm_{\mathtt{count}} = \sum_{x\in \cR} \delta_x$, the counting measure on $\cR$ 
and for all $\eta\ino (0, \infty)$, we set $q(  d_{\mathtt{cont}}, \eta)\! = \! \max 
\big\{ d_{\mathtt{cont}} (s,s^\prime) ; s, s^\prime \ino [0, 2n] : |s\! -\! s^\prime| \! \leq \! \eta \big\}$, where $d_{\mathtt{cont}}$ is defined by (\ref{Globasnitz}). Then, 
\begin{equation}
\label{SktVeitanderGlan}
\widetilde{d}_{\mathtt{Prok}} \big(\bm_{\mathtt{occ}} ,  \tfrac{n}{n^\prime}\bm_{\mathtt{count}}   \big) \leq 1 + 2q\big(  d_{\mathtt{cont}}, \frac{_{4\beta+1}}{^c} \big) \; .
\end{equation}
\end{lemma}
\noi
\textbf{Proof.} By convenience we set $u_n\! = \! \varnothing$. Let $i\ino \{ 0, \ldots, n\! -\! 1\}$. 
Recall from the definition of the contour function of $t$ in (\ref{Crigo}) that $b_i\!  :=\! 2i-|u_i| \! = \! \inf \{ k\ino \{ 0, \ldots, 2n\} : v(k) \! = \! u_i \}$. Therefore, for all $k\ino \{ b_i +1, \ldots , b_{i+1}\}$, 
$\{ v(l); 0\! \leq \! l\! < \! k\}\! = \! \{ u_j; 0\! \leq \! j \! < \! i+1\}$ and we get 
$$\# \big\{ \widehat{V}_l ;  0\! \leq \! l\! < \! k \big \}\! = \! \# \big\{ Y_{u_j}; 0\! \leq \! j \! < \! i+1 \big\}\! = \!  \# \big\{ Y_{u_j}; 0\! \leq \! j \! < \! i +1\big\}\! -2c(i+1)+ c(2(i+1) \!-\! k)+ ck \; .$$ 
Note that $b_i \! < \! k \! \leq \! b_{i+1}$ implies that $|2(i+1)-k| \! \leq \! 2 + \max_{0\leq j \leq n} |u_j|\! = \! = 2+ \max_{s\in [0, 2n]} C_s (t)$. Thus, we get 
$$  \max_{1\leq k\leq 2n}  \Big| \# \big\{ \widehat{V}_{l} \, ; \; 0\! \leq \! l \! < \! \! k \big\} - c k \Big| \leq \alpha + 2c + {c}\max_{s\in [0, 2n]} C_s (t) \; ,$$
which immediately implies (\ref{SktAndrae}).  

We next prove (\ref{SktVeitanderGlan}). To simplify set $J(k)\! = \!   \# \big\{ \widehat{V}_{l} \, ; \; 0\! \leq \! l \! < \! \! k \big\} $, for all $k\ino \{ 1, \ldots, 2n\}$ and thus $\beta \! = \! 
\max_{s\in (0, 2n]} |J( \lceil s \rceil  )\! -\! cs |$. 
For all $i^\prime \ino \{1, \ldots, n^\prime\}$ we next set $k(i^\prime)\! =\! \inf \{ k\ino \{ 1, \ldots, 2n \}\!  : 
J(k)\! = \! i^\prime \}$ and 
$Z_{i^\prime}\! = \! \widehat{V}_{k(i^\prime)-1}$. Note that the $Z_{i^\prime}$ are distinct and that $\cR \! = \! \{ Z_{i^\prime}; 1\! \leq \! i^\prime \! \leq \! n^\prime \}$. For all $i^\prime, j^\prime\ino \{ 1, \ldots, n^\prime \}$, $|k(i^\prime) \! -\! k(j^\prime)| \! \leq \!  c^{-1} ( 2\beta + |i^\prime \! -\! j^\prime | )$ since 
$J(k(i^\prime))\! = \! i^\prime $ and $J(k(j^\prime))\! = \! j^\prime $ and by definition of $\beta$. 
Then we get the following. 
\begin{equation}
\label{Grafenstein}
d_{\mathtt{gr}}^T \big( Z_{i^\prime} , Z_{j^\prime} \big) \! = \! 
d_{\mathtt{cont}} \big(   k(i^\prime) \! -\! 1,  k(j^\prime) \! -\! 1 \big) \! \leq \! 
q \big( d_{\mathtt{cont}}, |k(j^\prime) \! -\! k(i^\prime)  | \big)  \! 
 \leq \!  q \big( d_{\mathtt{cont}},   \frac{_{2\beta+|i^\prime  - j^\prime |}}{^c} \big).
\end{equation}
Next observe that for all $k(i^\prime) \! \leq  \! k \! < \! k (i^\prime+1)$, we get $J(k)\! = \! i^\prime$ and 
$$ d_{\mathtt{gr}}^T \big( \widehat{V}_{k-1} , Z_{J(k)} \big) =d_{\mathtt{cont}} ( k\! -\! 1, k(i^\prime) \! -\! 1)\leq q \big( d_{\mathtt{cont}}, k(i^\prime \! +1) \! -\! k(i^\prime)  \big)\!\! \!  \overset{\textrm{by (\ref{Grafenstein})}}{\leq } \!\! \! 
q \big( d_{\mathtt{cont}},   \frac{_{2\beta+1}}{^c} \big) \; $$
Thus, we have proved that 
\begin{equation}
\label{hgghgh}
\max_{s\in (0, 2n]} d_{\mathtt{gr}}^{T} 
\big(\widehat{V}_{\lfloor s \rfloor } , Z_{J(\lceil s \rceil)}  \big) \! \leq \! q \big( d_{\mathtt{cont}},   \frac{_{2\beta+1}}{^c} \big) \; .
\end{equation}
Next observe that $|n^\prime \! -\! 2cn| \! \leq \! \beta$ and that 
$$ J_{\lceil s \rceil} -\lceil \frac{_{n^\prime}}{^{2n}}s \rceil  =  J_{\lceil s \rceil} -cs - \frac{_s}{^{2n}} (n^\prime \! -\! 2cn) +  \frac{_{n^\prime}}{^{2n}}s -\lceil \frac{_{n^\prime}}{^{2n}}s \rceil . $$
Thus, $\max_{s\in (0, 2n]} |J( \lceil s \rceil  )- \lceil n^\prime s/2n \rceil | \! \leq \!  1+ 2\beta$ and by (\ref{Grafenstein}), $d^{_T}_{^\mathtt{gr}}  \big(  Z_{J(\lceil s \rceil)} , 
Z_{\lceil n^\prime s/2n \rceil} \big) \! \leq \!    q \big( d_{\mathtt{cont}},   \frac{_{4\beta+1}}{^c} \big)$. This inequality combined with (\ref{hgghgh}) implies that 
\begin{equation}
\label{SpittalanderDrau}
\max_{s\in (0, 2n]} d_{\mathtt{gr}}^{T} 
\big(\widehat{V}_{\lfloor s \rfloor } , Z_{\lceil n^\prime s/2n \rceil} \big)  \leq     2q \big( d_{\mathtt{cont}},   \frac{_{4\beta+1}}{^c} \big)\; .
\end{equation}
We then denote by $\pi $ the image on $T\! \times \! T$ 
of the Lebesgue measure on $[0, 2n]$ via the map $s\! \mapsto \!  \big(\widehat{V}_{\lfloor s \rfloor } , Z_{\lceil n^\prime s/2n \rceil} \big)$. The first marginal of $\pi$ is by definition 
$\bm_{\mathtt{cont}}$ and we easily see 
that the second marginal of $\pi$ is $\frac{2n}{n^\prime} \bm_{\mathtt{count}}$. 
By (\ref{SpittalanderDrau}), it implies that 
$\widetilde{d}_{\mathtt{Prok}} \big( \bm_{\mathtt{cont}} , \tfrac{2n}{n^\prime} 
\bm_{\mathtt{count}} \big)\!  \leq \! 2
q\big(  d_{\mathtt{cont}}, \frac{_{4\beta+1}}{^c} \big) $. 
Recall from Lemma \ref{Tobiach} that $\widetilde{d}_{\mathtt{Prok}} (\bm_{\mathtt{cont}}, 2\bm_{\mathtt{occ}}) \! \leq \! 1 $. Thus, we get 
$\widetilde{d}_{\mathtt{Prok}} \big( 2\bm_{\mathtt{occ}} , \tfrac{2n}{n^\prime} 
\bm_{\mathtt{count}} \big)\!  \leq \! 1+ 2
q\big(  d_{\mathtt{cont}}, \frac{_{4\beta+1}}{^c} \big) $, which implies (\ref{SktVeitanderGlan}) by Remark \ref{Kremsmunster}.  \cqfd

\noi
\textbf{The range of the branching random walk.}
We next discuss the connection between the graph-metric of the range of branching random walks and the associated snake metric. More precisely, let  
$t\ino \bbT$ be a \textit{finite} (rooted ordered) tree as in Definition \ref{ort} and let $(Y_v)_{v\in t}$ be a $\bbW_{\! [0, 1]}$-valued branching random walk with law $Q^+_{\mathbf{t}}$ as in Definition \ref{defimodel} $(ii)$. Recall from (\ref{Weissbach}) that $\cR \! = \! \{ Y_v; v\ino t\}$ 
that is a subtree of $\bbW_{\! [0, 1]}$ as in Definition \ref{Wsubtree}. 
We denote by $d_{\mathtt{gr}}$ the graph distance in $\bbW_{\! [0, 1]}$ 
and we recall from Corollary \ref{freeran} the fundamental property that makes the range $\cR$ tractable: 
\begin{equation}
\label{refreeran}
\textrm{$\bP$-a.s.~for all $u,v\ino t$,}  \quad d_{\mathtt{gr}} (Y_u, Y_v)\! = \! |Y_u|+|Y_v|\! -\! 2 
\min_{w\in \lgeo u, v \rgeo} |Y_w| \; . 
\end{equation}
Next, let {$(\widetilde{t}, \widetilde{d}^{{\, t}}_{\mathtt{gr}})$} and $(\widetilde{\bbW}_{\! [0, 1]}, \widetilde{d})$ be the real trees spanned by respectively $t$ and $\bbW_{\! [0, 1]}$; with a slight abuse of notation, we suppose that $t\! \subset \! \widetilde{t}$ and that $\bbW_{\! [0, 1]}\! \subset \! \widetilde{\bbW}_{\! [0, 1]}$. Let 
$\widetilde{Y}$ be the continuous interpolation of $Y$ in $\widetilde{\bbW}_{\! [0, 1]}$ as defined by (\ref{Yexten}) (with $(T, \mathtt{o})\! = \! (\bbW_{\! [0, 1]} , \varnothing)$) 
and recall from (\ref{Weissbach}) that $\widetilde{\cR}\! = \! \{ \widetilde{Y}_\sigma ; \sigma \ino \widetilde{t} \}$ is the tree spanned by $\cR$. 
We easily check that (\ref{refreeran}) extends to $\widetilde{\cR}$. Namely, 
\begin{equation}
\label{spfreeran}
\textrm{$\bP$-a.s.~for all $\sigma, \sigma^\prime\ino \widetilde{t}$,}  \qquad \widetilde{d}  (\widetilde{Y}_{\sigma}, \widetilde{Y}_{\sigma^\prime})=  |\widetilde{Y}_{\sigma}|+|\widetilde{Y}_{\sigma^\prime}| - 2 \!\!\! 
\min_{\varsigma \in \lgeo \sigma , \sigma^\prime \rgeo} \!\!\! |\widetilde{Y}_{\varsigma}| \; . 
\end{equation}

Let us rewrite this formula in terms of the contour function $C(t)$ of $t$, the spatial contour 
$V$, its endpoint process $\widehat{V}$ as defined in (\ref{kotschach}), and its corresponding snake 
$\gzW$ as defined in (\ref{contsnk}). To that end, first recall from (\ref{Globasnitz}) the definition of the spatial contour pseudo-distance $d_{\mathtt{cont}} $ and recall from (\ref{snkmeff}) the notation $M_{C(t), \gzW} (\cdot , \cdot ) $ and the definition of the snake distance $d_{C(t), \gzW}$  associated with the snake 
$(C(t), \gzW)$. Then, (\ref{spfreeran}) translates into the following:  $\bP$-a.s.~for all $s_1, s_2 \ino [0, 2 \# t]$,  
 \begin{equation}
\label{spfreeran_bis}
 d_{\mathtt{cont}} (s_1,s_2) =\widetilde{d} (\widehat{V}_{s_1}, \widehat{V}_{s_2})=  
\widehat{\gzW}_{s_1} + \widehat{\gzW}_{s_2}- 2 M_{C(t), \gzW} (s_1, s_2) = d_{C(t), \gzW} (s_1, s_2) \; .
\end{equation}
Recall from (\ref{Gloggnitz}) the definition of the occupation measure $\widetilde{\bm}_{\mathtt{cont}}$ induced by the spatial contour. Thus, (\ref{Troepolach}) and (\ref{spfreeran_bis}) 
imply that 
\begin{equation}
\label{Klagenfurt}
\textrm{$\bP$-a.s.} \quad  \big( \widetilde{\cR}, \widetilde{d}, \varnothing, \widetilde{\bm}_{\mathtt{cont}} \big) \; \textrm{and} \; \big( T_{C(t), \gzW}, d_{C(t), \gzW}, r_{C(t), \gzW}, \mu_{C(t), \gzW} \big) \; \textrm{are isometric.}
\end{equation}
where we recall from (\ref{notasnkm}) that 
$(T_{C(t), \gzW}, d_{C(t), \gzW}, r_{C(t), \gzW}, \mu_{C(t), \gzW})$ stands for the pointed measured compact real tree coded by the pseudo-metric $d_{C(t), \gzW}$.

We next couple the free branching random walk $(Y_v)_{t\in t}$ with the $\bbW_\bgg$-valued one via the $\bgg$-contraction application $\Phi_\bgg$ as in Definition \ref{tracage}: namely we set 
$$ \forall v\in t, \quad  Y^{\bgg}_v\! = \! \Phi_\bgg (Y_v), \quad \textrm{and} \quad \cR_\bgg \! = \! 
\big\{ Y^\bgg_v; v\ino t \big\} \; .$$  
By Remark \ref{ooobvvv}, $Y^\bgg$ has law $Q^{+\bgg}_{\mathbf{t}}$ as defined in Definition \ref{defimodel} $(i)$. 
Observe that $\cR_\bgg$ is a subtree of $\bbW_{\! \bgg}$ as in Definition \ref{Wsubtree}. 
We next denote by $d^{\bgg}$ the graph distance in $\bbW_{\! \bgg}$ and we denote by
 $(\widetilde{\bbW}_{\! \bgg}, \widetilde{d}^{_{ \, \bgg}}_{^{\! }} )$ the real tree spanned by $\bbW_{\! \bgg}$; with a slight abuse of notation, we suppose that $\bbW_{\! \bgg}\! \subset \! \widetilde{\bbW}_{\! \bgg}$. 
 We extend $\Phi_\bgg$ as a map from $\widetilde{\bbW}_{\! [0, 1]} $ to $\widetilde{\bbW}_{\! \bgg}$ explained in Remark \ref{adjadja} and we set 
$$ \forall s\ino [0, 2\# t], \, \forall r\ino [0, C_s(t)], \quad V^\bgg_s(r)= \Phi_\bgg (V_s(r)) \quad \textrm{and} \quad \widehat{V}_s^\bgg=   V^\bgg_s (C_s(t)) \; ,$$
where we recall that $(V_s (\cdot))_{s\in [0, 2\#t]}$ stands for the spatial contour of $Y$ as defined in (\ref{spacctr}). 
Note that $\widehat{V}_s^\bgg\! = \! \Phi_\bgg (\widehat{V}_s)$. By Remark \ref{adjadja}, it turns out that 
$V^\bgg$ is the spatial contour associated with $Y^\bgg$. 
Then, denote by 
$\widetilde{\cR}_\bgg$ the compact 
real tree spanned by $\cR_\bgg$ in $\widetilde{\bbW}_{\! \bgg}$. 
Namely,  
$$\widetilde{\cR}_\bgg\! = \! \{ \widehat{V}_s^\bgg; s\ino [0, 2\# t ] \} \! = \!  \Phi_\bgg (\widetilde{\cR}) \; .$$ 
We next extend Lemma \ref{contrace} as follows: since  
$(\widehat{V}^\bgg_s)_{s\in [0, 2\# t]}$ is the continuous interpolation of the path $(Y^\bgg_{v(k)})_{0\leq k\leq 2\# t}$ in $\cR_\bgg$, we easily derive the following from (\ref{zip}). 

\smallskip

\begin{compactenum}
\item[] \textit{For all $s_1, s_2 \ino [0, 2\# t] $, 
there exists a nonnegative r.v.~$G_{s_1,s_2}$ such that}
\end{compactenum}
\begin{equation}
\label{rezip}
2G_{s_1,s_2}\! = \! \widetilde{d} \big( \widehat{V}_{s_1}, \widehat{V}_{s_2} \big)  - \widetilde{d}^{_{\, \bgg}}  \big(\widehat{V}^\bgg_{s_1}, \widehat{V}^\bgg_{s_2}\big) \quad \textit{and} \quad \bP \big( G_{s_1,s_2} 
 \! \geq x \big) \leq \bgg^{2-x}, \; x\ino [0, \infty)  
\end{equation} 
This combined with (\ref{Klagenfurt}) shows that $T_{C(t), W} $ 
and $\cR_\bgg$ are close in a rough sense. It turns out that it is sufficient for the proof of Theorem \ref{main2}.

\subsection{Invariance principle for discrete snakes.}
\label{invsnksc}  
In this section we recall one important result due to Marzouk (Theorem 1 \cite{Mar20}) that is an invariance principle for {real valued} endpoint processes of discrete snakes. Here we concentrate on the one-dimensional {case} but we shall need actually a slightly stronger version that holds for path-valued snakes and not only for endpoint processes.

More precisely, let us fix $\gamma \ino (1, 2]$ and 
let $\tau$ be a Galton--Watson tree (as in Definition \ref{GWtreedef}) whose offspring distribution $\mu$ satisfies $(\mathbf{H})$ as in (\ref{hyposta}). Conditionally given $\tau$, 
let $(Y_v)_{v\in \tau}$ 
be a $\bbZ$-valued branching random walk whose transition kernel $q(y, dy^\prime)$ is that of the simple symmetric random walk on $\bbZ$.  
Then recall from (\ref{contsnk}) the definition of $(\gzW_s (\cdot) )_{s\in [0, 2\#\tau]}$, the continuous interpolation of the discrete snake (here the tree $T$ is $\bbZ$). Then,  
\begin{equation}
 \label{Lienz}
\forall k\ino \{ 0, \ldots, 2\# t\} , \quad \widehat{\gzW}_k = Y_{v(k)} \; , 
\end{equation} 
where $v(k) \ino \tau$, $0\! \leq \! k \! \leq \! 2\#\tau $, stands for 
the sequence of vertices of $\tau$ visited by the contour (or the depth-first exploration) of $\tau$. 
Then, the following result is a special case of Theorem 1 in Marzouk \cite{Mar20}: 
\begin{theorem}[Theorem 1 \cite{Mar20}] 
\label{cvdissnk} Let $\gamma$, $\mu$, $\tau$ and $(a_n)_{n\in \bbN}$ be as in Theorem \ref{rappcvtr}.   
Conditionally given $\tau$, let $(Y_v)_{v\in \tau}$ be defined as above 
and let $(\widehat{\gzW}_k)_{0\leq k\leq 2\# \tau}$ be the endpoint process of the snake associated with $(Y_v)_{v\in \tau}$ as in (\ref{Lienz}). Then 
\begin{equation*}
\Big( \tfrac{1}{a_n} C_{2ns}(\tau) ,  \tfrac{1}{\sqrt{a_n}}  \widehat{\gzW}_{\lfloor 2ns \rfloor} \Big)_{\! \! s\in [0, 1]} \; \textrm{under $\; \bP \big(  \cdot  \big| \, \# \tau \! = \! n \big) $ }  \overset{\textrm{(law)}}{\underset{n\rightarrow \infty}{-\!\!\! -\!\!\!-\!\!\! -\!\!\!  -\!\! \! \longrightarrow}} \; (H_s, \widehat{W}_s)_{s\in [0, 1]} 
\end{equation*}
where $(H, W)$ is the normalised one dimensional Brownian snake with $\gamma$-stable branching mechanism and where $\widehat{W}_s \! = \! W_s (H_s)$, $s\ino [0, 1]$, stands for its endpoint process. 
\end{theorem}
As already mentioned this result actually holds for fairly more general spatial motions in $\bbR^d$. It extends earlier results by Janson and Marckert \cite{JanMar05} who considered snakes whose genealogical tree are in the domain of attraction of the Brownian tree. We refer to Marzouk \cite{Mar20} for more details. 
We next derive from Theorem \ref{cvdissnk} the following proposition that holds for 
the path-valued continuous interpolation of the discrete snake $(\gzW_s (\cdot))_{s\in [0, 2\# \tau]}$. \begin{proposition} 
\label{contsnlim} We keep the same notation and the same assumption as in Theorem \ref{cvdissnk}. 
Recall from Definition \ref{rePopol} $(a)$ the metric $\Delta$ that makes the space $\bC([0, 1], \bbR) \! \times \! \bC([0, 1], \bC_0)$ Polish. Then, weakly on that space, the following convergence holds. 
\begin{equation}
\label{strsnkk}
\Big( \big( \tfrac{1}{a_n} C_{2ns}(\tau) )_{\! s\in [0, 1] }\, ,  \big( \tfrac{1}{\sqrt{a_n}} \gzW_{ 2ns }  (\tfrac{\cdot}{a_n}) \big)_{\!  s\in [0, 1]} \Big)\; \textrm{under $\; \bP \big(  \cdot  \big| \, \# \tau \! = \! n \big) $ }  \overset{\textrm{(law)}}{\underset{n\rightarrow \infty}{-\!\!\! -\!\!\!-\!\!\! -\!\!\!  -\!\! \! \longrightarrow}}\big(H, W \big). 
\end{equation}
By taking the absolute value, this results holds true for the snake reflected at $0$. 
\end{proposition}
\noi
\textbf{Proof.} Let $k_0\! = \! 0\! \leq \! k_1 \! \leq \! \ldots \! \leq \! k_p \! \leq \! 2\# \tau$ be integers. Then for all $j\ino \{0, \ldots, p\! -\! 1\}$, set $b_j\! = \!  \min_{k_j \leq \ell \leq k_{j+1}} C_\ell (\tau)$ and $b^\prime_j\! = \! C_{k_{j+1}} (\tau)\! -\! b_j$. 
By definition of the branching random walk and its associated snakes, the following holds true conditionally given $\tau$ (or $C(\tau)$). 

\smallskip

\begin{compactenum} 

\item[$(i)$] The paths $S(j)\! := \! (\gzW_{k_{j+1}} \big( (b_j+\ell) \! \wedge \! C_{k_{j+1}} (\tau) \big) \!- \! \gzW_{k_{j+1}} (b_j)  \big)_{\ell \in \bbN}$ are (conditionally) independent. 

\smallskip

\item[$(ii)$] The conditional law of $S(j)$ is that of a $\bbZ$-valued symmetric random walk, starting at $0$ and stopped at time $b_{j+1}^\prime$.
\end{compactenum} 

\noi
Next observe that $(b,b^\prime\! , w, w^\prime) \ino [0, \infty)^2 \! \times \! \bC_0^2 \! \mapsto \! (w(r\wedge b)+ w^\prime((r\! -\! b)_+ \! \wedge \! b^\prime ) \! -\! w^\prime (0))_{r\in [0, \infty)}$ is continuous (recall from Definition \ref{rePopol} that $\bC_0$ is equipped with the Polish topology of the convergence on all compact subsets of times, which corresponds for instance to the metric $\delta_{\mathtt{u}}$ given in (\ref{uttcpct})). This combined with $(i)$, $(ii)$ and easy arguments on linear interpolation imply that for all real numbers $0\! \leq \! s_1 \! \leq \! \ldots \! \leq \! s_p \! \leq  \! 1$, 
$(\mathbf{w}^{_{(n)}}_{{s_j}})_{1\leq j\leq p}$ under $\bP(  \cdot | \# \tau \! = \! n)$ converges  weakly to $(W_{s_j})_{1\leq j\leq p}$ on $\bC_0^p$, where we have set 
$\mathbf{w}^{_{(n)}}_{s} (\cdot)\! =\! a_{^n}^{_{-1/2}}\gzW_{2ns} (\cdot/a_n)$, for all $s\ino [0, 1]$. 

We next prove the tightness of the rescaled snakes. A standard argument on linear interpolation shows that Theorem \ref{cvdissnk} implies that the continuous endpoint processes $\widehat{\mathbf{w}}^{_{(n)}}$ under $\bP(  \cdot \, | \# \tau \! = \! n)$ converge  weakly to the endpoint process $\widehat{W}$. This implies that 
$$\forall \epp \ino (0, \infty), \quad   \lim_{\eta \rightarrow +} \sup_{n\geq 0}\bP ( \omega_\eta (\widehat{\mathbf{w}}^{_{(n)}}) \! \geq \! \epp \, | \# \tau \! = \! n )\! = \! 0, $$
where $\omega_\eta (\widehat{\mathbf{w}}^{_{(n)}})$ stands for the $\eta$-uniform modulus of continuity of $\widehat{\mathbf{w}}^{_{(n)}}$ over $[0, 1]$. We next denote by $\omega_\eta (\mathbf{w}^{_{(n)}})$ the $\eta$-uniform modulus of continuity of $\mathbf{w}^{_{(n)}}$ with respect to the metric $\delta_{\mathtt{u}}$ on $\bC_0$ as introduced in Lemma \ref{Innsbruck}; this lemma  asserts that $\omega_\eta (\mathbf{w}^{_{(n)}}) \! \leq \!  2\omega_\eta (\widehat{\mathbf{w}}^{_{(n)}})$. This, combined with the weak convergence of the finite dimensional marginals of $\mathbf{w}^{_{(n)}}$, {entails} that the laws of $\mathbf{w}^{_{(n)}}$ are tight on $\bC([0, 1] , \bC_0)$ and we easily get (\ref{strsnkk}). \cqfd

\subsection{Proof of Theorem \ref{main2}.}
\label{Pfmain2}
Let us fix $\gamma \ino (1, 2]$ and 
let $\tau$ be a Galton--Watson tree (as in Definition \ref{GWtreedef}) whose offspring distribution $\mu$ satisfies $(\mathbf{H})$ as in (\ref{hyposta}). 
Then recall from (\ref{Crigo}) the definition of the contour process 
$(C_s (\tau))_{s\in [0, 2\# \tau]} $ of $\tau$. 
Conditionally given $\tau$, let $(Y_v)_{v\in \tau}$ 
be a $\bbW_{\! [0, 1]}$-valued branching random walk with law $Q^+_{\ftau}$ 
as in Definition \ref{defimodel} $(ii)$. 
Recall from (\ref{spacctr}) 
the spatial contour $(V_s (\cdot))_{s\in [0, 2 \# \tau]}$ associated with $Y$, and recall from 
(\ref{contsnk}) the definition of $(\gzW_s (\cdot) )_{s\in [0, 2\#\tau]}$, the continuous interpolation of the discrete snake. 
We denote by $(\widetilde{\tau}, \widetilde{d}^{{\, \tau}}_{^{\mathtt{gr}}})$, 
$(\widetilde{\bbW}_{\! [0, 1]}, \widetilde{d})$ and 
$(\widetilde{\bbW}_{\! \bgg }, \widetilde{d}^{ \bgg})$ 
be the real trees spanned by respectively $\tau$, $\bbW_{\! [0, 1]}$ and $\bbW_{\! \bgg}$ 
we assume that $\tau\! \subset \! \widetilde{\tau}$, that $\bbW_{\! [0, 1]}\! \subset \! \widetilde{\bbW}_{\! [0, 1]}$ and that $\bbW_{\! \bgg}\! \subset \! \widetilde{\bbW}_{\! \bgg}$. 
We next couple the free branching random walk $(Y_v)_{t\in t}$ with the $\bbW_\bgg$-valued one via the $\bgg$-contraction application $\Phi_\bgg$ as in Definition \ref{tracage}: namely we set 
$Y^{\bgg}_v\! = \!  \Phi_\bgg (Y_v)$ for all $v\ino \tau$. 
By Remark \ref{ooobvvv}, $Y^\bgg$ has law $Q^{+\bgg}_{\mathbf{t}}$ as in Definition \ref{defimodel} $(i)$. 
 We also set $V^\bgg\! = \! \Phi_\bgg (V)$ that is the spatial contour of $Y^\bgg$, according to Remark \ref{adjadja}. We recall the following notation $\cR \! = \! 
\big\{ Y_v; v\ino \tau \big\}$, $\cR_\bgg \! = \! 
\big\{ Y^\bgg_v; v\ino \tau \big\}\! = \! \Phi_d (\cR)$, $\widetilde{\cR} \! = \! \big\{ \widehat{V}_s ; s\ino [0, 2\# \tau] \big\}$ and $\widetilde{\cR}_\bgg\! = \! \big\{ \widehat{V}^\bgg_s ; s\ino [0, 2\# \tau] \big\}\! = \! \Phi_\bgg (\widetilde{\cR})$. 

We trivially extend $C(\tau)$ and $\gzW$ on $[0, \infty)$ by taking $C_s(\tau)$ equal to $0$ and by taking $\gzW_s$ equal to the null function for all $s\ino [2\# \tau, \infty)$.   
Let $(a_n)_{n\in \bbN}$ be as in Theorem \ref{rappcvtr}.  
To simplify notation we set for all $s\ino [0, 1]$, 
$$  \mathbf{h}_n (s) \! = \! \tfrac{1}{a_n} C_{2ns}(\tau),  \quad  \mathbf{w}^{(n)}_s (r) 
= \tfrac{1}{\sqrt{a_n}} \gzW_{ 2ns }  \big( \tfrac{r}{a_n} \big) , \; r\ino [0, \infty)  \quad \textrm{and} \quad 
\widehat{\mathbf{w}}^{(n)}_s  \! = \!  \tfrac{1}{\sqrt{a_n}} \widehat{\gzW}_{ 2ns }. $$ 
Observe that $\bP (\, \cdot \, | \, \# \tau \! = \! n)$-a.s.~$(\mathbf{h}_n,  \mathbf{w}^{(n)})\ino \fSigma ([0, 1])$ as in Definition \ref{rePopol} and we denote by $d_{\mathbf{h}_n }$ and  $d_{\mathbf{h}_n,  \mathbf{w}^{(n)}}$ the corresponding tree and snake pseudo-distances as defined in resp.~(\ref{djhczkk}) and (\ref{snkmeff}). In particular, observe that 
\begin{equation*}
\textrm{$\bP$-a.s.~for all $s_1, s_2\ino [0, \infty)$,} \quad  
\tfrac{1}{\sqrt{a_n}} \widetilde{d}  \big(\widehat{V}_{2ns_1} , \widehat{V}_{2ns_2} \big)= d_{\mathbf{h}_n,  \mathbf{w}^{(n)}}(s_1, s_2) \; 
\end{equation*}
Proposition \ref{contsnlim} combined with Lemma \ref{Hermagor} implies that weakly on $(\bC([0, 1]^2,\bbR), \lVert \cdot \rVert)^2$ 
\begin{equation*}
(d_{\mathbf{h}_n }, d_{\mathbf{h}_n,  \mathbf{w}^{(n)}}) \; \textrm{under $\bP \big(  \cdot  \big| \, \# \tau \! = \! n \big) $ } 
\underset{n\rightarrow \infty}{-\!\!\! -\!\!\!-\!\!\! -\!\!\!  \longrightarrow} \; \big(d_H, d_{H, W} \big), 
\end{equation*}  
where $(H, W)$ is the normalised one dimensional reflected Brownian snake with $\gamma$-stable branching mechanism as in Definition \ref{brosnadef} and where $d_H$ and $d_{H, W}$ are the corresponding tree and snake pseudo-distances. 

  We next set 
$$ \forall s_1, s_2 \ino [0, \infty), \quad  \mathbf{d}_n^*\! (s_1, s_2)=  \tfrac{1}{\sqrt{a_n}} 
\widetilde{d}^{\bgg}  \big(\widehat{V}^{\bgg}  _{2ns_1} , \widehat{V}^{\bgg}  _{2ns_2} \big) \; .$$
By (\ref{rezip}), conditionally given $\tau$, $ \mathbf{d}_n^* (s_1, s_2) \! \leq \!  d_{\mathbf{h}_n,  \mathbf{w}^{(n)}}(s_1, s_2)$ and for all $\epp \ino (0, 1)$, a.s.~
$$ \bP \big(d_{\mathbf{h}_n,  \mathbf{w}^{(n)}}(s_1, s_2) \! -\!  \mathbf{d}_n^* (s_1, s_2) \! > \! 2\epp \, \big| \, \tau   \big)  \leq \bgg^{3- \epp \sqrt{a_n}} \; , $$
which implies that $\lim_{n\rightarrow \infty} \bP \big(d_{\mathbf{h}_n,  \mathbf{w}^{(n)}}(s_1, s_2) \! -\!  \mathbf{d}_n^*(s_1, s_2) \! >\!  2\epp \, \big| \, \# \tau \! = \! n  \big)  \! = \! 0$. Therefore, Proposition \ref{submetric} applies (with $\mathbf{d}_n\! = \! d_{\mathbf{h}_n,  \mathbf{w}^{(n)}}$) 
to show that 
$ (\mathbf{d}_n^*,  d_{\mathbf{h}_n,  \mathbf{w}^{(n)}})$ under $\bP (\, \cdot \, | \# \tau \! = \! n)$ converge to $(d_{H, W}, d_{H,W})$, weakly on $(\bC([0, 1]^2,\bbR), \lVert \cdot \rVert)^2$. Actually it is easy to see that 
$(d_{\mathbf{h}_n},  \mathbf{d}_n^*,  d_{\mathbf{h}_n,  \mathbf{w}^{(n)}})$ under $\bP (\, \cdot \, | \# \tau \! = \! n)$ converges to $(d_H, d_{H, W}, d_{H,W})$, weakly on $(\bC([0, 1]^2,\bbR), \lVert \cdot \rVert)^3$. 

We next recall that $\bP (\, \cdot \, | \# \tau \! = \! n)$-a.s.~the real tree coded by $d_{\mathbf{h}_n}$ is isometric to 
$(\widetilde{\tau}, \tfrac{_1}{^{a_n}} \widetilde{d}^{_{\, \tau}}_{^{\mathtt{gr}}}, \varnothing, \tfrac{_1}{^{n}}\widetilde{\mathtt{m}})$ where $ \widetilde{\mathtt{m}}\! = \! \delta_{\varnothing} + \mathtt{Length}$. Also recall from (\ref{Troepolach}) that $\bP (\, \cdot \, | \# \tau \! = \! n)$-a.s.~the real tree coded by $\mathbf{d}_n^*$ is isometric to $(\widetilde{\cR}_\bgg, \tfrac{_1}{^{\sqrt{a_n}}} \widetilde{d}^{_{\, \bgg}}_{^{\! }}, \varnothing, \tfrac{_1}{^{2n}}\widetilde{\bm}_{\mathtt{cont}}^\bgg )$, where $\widetilde{\bm}_{\mathtt{cont}}^\bgg$ is the occupation measure of $(\widehat{V}^\bgg_s)_{s\in [0, 2n]}$, as defined in general in (\ref{Gloggnitz}). 
Since $(d_{\mathbf{h}_n},  \mathbf{d}_n^*)$ under $\bP (\, \cdot \, | \# \tau \! = \! n)$ converges to $(d_H, d_{H, W})$, weakly on $(\bC([0, 1]^2,\bbR), \lVert \cdot \rVert)^2$, Proposition \ref{groweak} implies that the pointed measured compact spaces 
$\big( \widetilde{\tau}, \tfrac{1}{a_n} \widetilde{d}^{_{\, \tau}}_{^{\mathtt{gr}}}  , \varnothing, \tfrac{1}{n} \widetilde{\mathtt{m}} \big)$ and $\big( \widetilde{\cR}_\bgg, \tfrac{1}{\sqrt{a_n}}\widetilde{d}^{\, \bgg}, \varnothing, \tfrac{1}{2n} \widetilde{\bm}^\bgg_{\mathtt{cont}} \big)$ under $\bP (\, \cdot \, | \# \tau \! = \! n)$ jointly converge to resp. $(T_H, d_H, r_H, \mu_H)$ and $(T_{H, W}, d_{H, W}, r_{H, W}, \mu_{H, W})$ weakly on $(\bbM, \bdelta_{\mathtt{GHP}})^2$, where $T_H$ and $T_{H,W}$ here stand for the real trees coded by resp.~$H$ and $(H, W)$.

We next recall the following notation. 
$$ \mathtt{m}\! = \! \sum_{v\in \tau} \delta_v, \; \bm_{\mathtt{occ}}^\bgg\! = \! \! 
\sum_{v\in \tau} \delta_{Y^\bgg_v} , \;  \overline{\bm}_{\mathtt{occ}}^\bgg\! = \!  \tfrac{1}{n}\bm_{\mathtt{occ}}^\bgg,  \;   \bm_{\mathtt{count}}^\bgg\! \! = \! \!
\sum_{x\in \cR_\bgg} \! \delta_{x}  \; \, \textrm{and} \; \, \overline{\bm}_{\mathtt{count}}^\bgg\! \! = \! \frac{_1}{^{\# \cR_\bgg}} \bm_{\mathtt{count}}^\bgg .$$ 
By the inequalities specified in Example \ref{rempli2}, we get $d_{\mathtt{Haus}} (\tau, \widetilde{\tau})\! \leq \! 1$ and 
$  d_{\mathtt{Prok}} (\widetilde{\mathtt{m}},\mathtt{m})\! \leq \! 3$. Similarly, by 
Lemma \ref{Tobiach} we get $d_{\mathtt{Haus}} (\cR_\bgg, \widetilde{\cR}_\bgg) \! \leq \! 1$ and 
$d_{\mathtt{Prok}} ( 2\bm_{\mathtt{occ}}^\bgg, \widetilde{\bm}^\bgg_{\mathtt{cont}}) \! \leq \! 2$. Thus, 
\begin{eqnarray*}
\lefteqn{ 
\Big( \big(\tau , \tfrac{1}{a_n} d_{{\mathtt{gr}}}  , \varnothing, \tfrac{1}{n} \mathtt{m} \big)\, , \,  \big( \cR_\bgg, \tfrac{1}{\sqrt{a_n}}d_{\mathtt{gr}}, \varnothing, \tfrac{1}{n} 
\bm^\bgg_{\mathtt{occ}} \big) \Big) 
 \quad  \textrm{under $\; \bP \big(  \cdot  \big| \, \# \tau \! = \! n \big) $ } } \hspace{40mm} 
                    \nonumber\\
 & &\underset{n\rightarrow \infty}{-\!\!\!-\!\!\! -\!\!\!  \longrightarrow} \; \big( (T_H, d_H, r_H, \mu_H)\, ,\,  (T_{H, W}, d_{H, W}, r_{H, W}, \mu_{H, W})\big), 
\end{eqnarray*}
where $d_{\mathtt{gr}}$ stands for both graph-tree distances in $\tau$ and $\bbW_{\! \bgg}$ to simplify notation. 

We next control $\bm_{\mathtt{count}}^\bgg$ in terms of $ \bm_{\mathtt{occ}}^\bgg$ thanks to Theorem \ref{main1} and Lemma \ref{Loibach}. To that end, we set $c\! = \! c_{\mu, \bgg} /2$ where $c_{\mu, \bgg}$ is as in Theorem \ref{main1}; we denote by $u_0, u_1, \ldots, u_{\# \tau -1}$ the vertices of $\tau$ listed in the lexicographical order and we also set 
$ \alpha \! =\!  \max_{1\leq i\leq \# \tau }  | \# \{ Y^\bgg_{u_j}; 0\! \leq \! j \! < \! i \} \! - \! 2c i | $ and $\beta \! = \!  \max_{s \in (0, 2\# \tau] } | \#\{ \widehat{V}^\bgg_{l} \, ; \; 0\! \leq \! l \! < \! \! \lceil  s \rceil \} \! -\!  c s |$. 
Theorem \ref{main1} asserts that for all $\epp \ino (0, \infty)$, $\bP (\alpha /n\! >\! \epp \, | \, \# \tau \! = \! n)$ tends to $0$ as $n\! \rightarrow \! \infty$. By (\ref{SktAndrae}), $\bP (\, \cdot \, | \, \# \tau \! = \! n)$-a.s.~
$$\frac{\beta}{n}\! \leq \! \frac{1}{n}(\alpha+ 3c)+ \frac{{c a_n}}{n} \max_{s\in [0, 1]} h_n(s)\; .$$
Since $(a_n)$ is $\tfrac{\gamma-1}{\gamma}$-regularly varying, $\lim_{n\rightarrow \infty} 
a_n/n\! =  \! 0$, which implies that 
$\bP (\beta /n\! >\! \epp \, | \, \# \tau \! = \! n)$ tends to $0$ as $n\! \rightarrow \! \infty$. 

Next, we denote by $d_{\mathtt{Prok}}$ the Prokhorov distance on the space of 
finite measures on $(\cR_\bgg, d_{\mathtt{gr}})$ and we denote by $d^{(n)}_{\mathtt{Prok}}$ the Prokhorov distance on the space of $(\cR_\bgg, \tfrac{_1}{^{\sqrt{a_n}}}d_{\mathtt{gr}} )$. 
Since $n^{-1}\! \leq \! a_{n}^{_{-1/2}}$ for all sufficiently large $n$, Remark \ref{Kremsmunster} combined with (\ref{SktVeitanderGlan}) 
implies that 
$\bP (\, \cdot \, | \, \# \tau \! = \! n)$-a.s.
$$d^{(n)}_{\mathtt{Prok}} \big(
\overline{\bm}^\bgg_{\mathtt{occ}} , \overline{\bm}^\bgg_{\mathtt{count}},  \big) \!\! \!\!  \!\!  \!\! \overset{\textrm{{Remark \ref{Kremsmunster}}}}{\leq}   \!\! \!\! \!\! \!\!\tfrac{_1}{^{\sqrt{a_n}}}
d_{\mathtt{Prok}} \big(\bm^\bgg_{\mathtt{occ}} , \tfrac{n}{\# \cR_\bgg} \bm^\bgg_{\mathtt{count}} \big) \!\!  \!\! \overset{\textrm{{by (\ref{SktVeitanderGlan})}}}{\leq}   \!\! \!\! \tfrac{_1}{^{\sqrt{a_n}}}+ 2q (\mathbf{d}_n^*, \tfrac{4\beta+1}{2nc} ) , $$
where $q (\mathbf{d}_n^*,\eta )\! = \! \max \{\mathbf{d}_n^*(s,s^\prime) ; 
s,s^\prime \ino [0, 1] : |s\! -\! s^\prime | \! \leq \! \eta  \}$. Since $\mathbf{d}^*_n\! \rightarrow \! d_{H, W}$ weakly on $(\bC([0, 1]^2, \bbR), \lVert \cdot \rVert)$ and since for all $\epp \ino (0, \infty)$, 
$\lim_{n\rightarrow \infty}\bP (\beta /n\! >\! \epp \, | \, \# \tau \! = \! n)\! = \!  0$, 
Proposition \ref{tightness} implies that $\lim_{n\rightarrow \infty} \bP ( q (\mathbf{d}_n^*, \tfrac{4\beta+1}{2nc} ) \! >\! \epp \, | \, \# \tau \! = \! n)\! = \! 0$ and thus 
$$ \forall \epp \ino (0, 1) , \quad \lim_{n\rightarrow \infty} \bP \big(   d^{(n)}_{\mathtt{Prok}} \big(\overline{\bm}^\bgg_{\mathtt{occ}} , \overline{\bm}^\bgg_{\mathtt{count}}  \big)  \! >\! \epp \, \big| \, \# \tau \! = \! n \big)\! = \! 0\; .$$ 
Then, observe that $\# \cR_\bgg \! \leq \! n$ and that 
$\tfrac{\# \cR_\bgg}{n} \overline{\bm}^\bgg_{\mathtt{count}} \! = \! 
\tfrac{1}{n} \bm^\bgg_{\mathtt{count}}$. 
By Remark \ref{Kremsmunster}, we get 
\begin{eqnarray*}
d^{(n)}_{\mathtt{Prok}} \big( c_{\mu, \bgg} \overline{\bm}^\bgg_{\mathtt{occ}} , \tfrac{1}{n} \bm^\bgg_{\mathtt{count}} \big)  &\leq  & d^{(n)}_{\mathtt{Prok}} \big( c_{\mu, \bgg} \overline{\bm}^\bgg_{\mathtt{occ}} ,  \tfrac{\# \cR_\bgg}{n}  \overline{\bm}^\bgg_{\mathtt{occ}}   \big)  + d^{(n)}_{\mathtt{Prok}} \big( \tfrac{\# \cR_\bgg}{n}  \overline{\bm}^\bgg_{\mathtt{occ}}  , 
\tfrac{\# \cR_\bgg}{n}  \overline{\bm}^\bgg_{\mathtt{count}}   \big) \\
& \leq &   \Big|\tfrac{\# \cR_\bgg}{n} \! -\! c_{\mu, \bgg} \Big| + d^{(n)}_{\mathtt{Prok}} \big(\overline{\bm}^\bgg_{\mathtt{occ}} , \overline{\bm}^\bgg_{\mathtt{count}}  \big)  . 
\end{eqnarray*}
Now, recall  that $c_{\mu, \bgg}\! = \! 2c$ and that $|\tfrac{\# \cR_\bgg}{n} \! -\! c_{\mu, \bgg}| \leq \alpha/n$. 
It finally implies that 
$$ \forall \epp \ino (0, 1) , \quad \lim_{n\rightarrow \infty} \bP \big(   d^{(n)}_{\mathtt{Prok}} \big(c_{\mu, \bgg} \overline{\bm}^\bgg_{\mathtt{occ}} , \tfrac{1}{n} \bm^\bgg_{\mathtt{count}}  \big)  \! >\! \epp \, \big| \, \# \tau \! = \! n \big)\! = \! 0\; ,$$
which completes the proof of Theorem \ref{main2}. \cqfd

{\small
\setlength{\bibsep}{.3em}
\bibliographystyle{acm}
\bibliography{trace}
}

\end{document}